\input amstex
\documentstyle{amsppt}
\magnification=\magstep1
\vsize =21 true cm
\hsize =16 true cm
\loadmsbm
\topmatter

\centerline{\bf Exotic arithmetic structure on the first Hurwitz
                triplet}
\author{\smc Lei Yang}\endauthor
\endtopmatter
\document

\centerline{\bf Abstract}

\vskip 0.5 cm

  We find that the first Hurwitz triplet possesses two distinct
arithmetic structures. As Shimura curves $X_1$, $X_2$, $X_3$,
whose levels are with norm $13$. As non-congruence modular curves
$Y_1$, $Y_2$, $Y_3$, whose levels are $7$. Both of them are
defined over ${\Bbb Q}(\cos \frac{2 \pi}{7})$. However, for the
third non-congruence modular curve $Y_3$, there exist an
``exotic'' duality between the associated non-congruence modular
forms and the Hilbert modular forms, both of them are related to
${\Bbb Q}(e^{\frac{2 \pi i}{13}})$! Our results have relations and
applications to modular equations of degree fourteen (including
Jacobian modular equation and ``exotic'' modular equation),
``triality'' of the representation of $PSL(2, 13)$, Haagerup
subfactor, geometry of the exceptional Lie group $G_2$, and even
the Monster finite simple group ${\Bbb M}$!

\vskip 0.5 cm

\centerline{\bf Contents}
$$\aligned
 &\text{1. Introduction}\\
 &\text{2. Standard arithmetic structure on the first Hurwitz triplet}\\
 &\text{3. Six-dimensional representations of $PSL(2, 13)$ and}\\
 &\text{\quad exotic arithmetic structure on the first Hurwitz triplet}\\
 &\text{4. Seven-dimensional representations of $PSL(2, 13)$,}\\
 &\text{\quad Jacobian modular equation of degree fourteen and exotic duality theorem}\\
 &\text{5. Haagerup subfactor and exceptional Lie group $G_2$}\\
 &\text{6. Fourteen-dimensional representations of $PSL(2, 13)$,}\\
 &\text{\quad exotic modular equation of degree fourteen and Monster simple group ${\Bbb M}$}
\endaligned$$

\vskip 0.5 cm

\centerline{\bf 1. Introduction}

\vskip 0.5 cm

  A classical theorem of Hurwitz asserts that a Riemann surface
$S$ of genus $g>1$ can have at most $84(g-1)$ automorphisms, and a
group of order $84(g-1)$ is the automorphism group of some Riemann
surface of genus $g$ if and only if it is generated by an element
of order two and one of order three such that their product has
order seven. In that case the quotient of $S$ by the group is the
Riemann sphere, and the quotient map $S \to {\Bbb C} {\Bbb P}^1$
is ramified above only three points of ${\Bbb C} {\Bbb P}^1$, with
the automorphisms of orders two, three, seven of $S$ appearing as
the deck transformations lifted from cycles around the three
branch points (see \cite{El2}).

  A Riemann surface with the maximal number $84(g-1)$ of automorphisms,
regarded as an algebraic curve over ${\Bbb C}$, is called a
Hurwitz curve of genus $g$. Hurwitz curves can be characterized in
terms of their uniformization by the hyperbolic plane ${\Bbb H}$.
Any Riemann surface $S$ of genus greater than one can be
identified with ${\Bbb H}/\pi_1(S)$; conversely, any discrete
co-compact subgroup $\Gamma \subset \text{Aut}({\Bbb H}) \cong
PSL(2, {\Bbb R})$ that acts freely on ${\Bbb H}$ (that is, every
point has trivial stabilizer) yields a Riemann surface ${\Bbb
H}/\Gamma$ of genus greater than one whose fundamental group is
$\Gamma$. The automorphism group of ${\Bbb H}/\Gamma$ is
$N(\Gamma)/\Gamma$, where $N(\Gamma)$ is the normalizer of
$\Gamma$ in $\text{Aut}({\Bbb H})$. It follows that ${\Bbb
H}/\Gamma$ is a Hurwitz curve if and only if $N(\Gamma)$ is the
triangle group $G_{2, 3, 7}$ of orientation-preserving
transformations generated by reflections in the sides of a given
hyperbolic triangle with angles $\frac{\pi}{2}$, $\frac{\pi}{3}$,
$\frac{\pi}{7}$ in ${\Bbb H}$. Equivalently, $\Gamma$ is to be a
normal subgroup of $G_{2, 3, 7}$. Since $G_{2, 3, 7}$ has the
presentation
$$G_{2, 3, 7}=\langle \sigma_2, \sigma_3, \sigma_7: \sigma_2^2=
  \sigma_3^3=\sigma_7^7=\sigma_2 \sigma_3 \sigma_7=1 \rangle$$
with $\sigma_j$ being a $2 \pi/j$ rotation about the $\pi/j$
vertex of the triangle, this yields a characterization of the
groups that can occur as $\text{Aut}(S)=G_{2, 3, 7}/\Gamma$. We
can identify the Hurwitz curves $S$ with a Shimura curve by
recognizing $G_{2, 3, 7}$ as an arithmetic group in $PSL(2, {\Bbb
R})$, and $\pi_1(S)$ with a congruence subgroup of $G_{2, 3, 7}$
(see \cite{El2}).

  In the theory of Riemann surfaces, the first Hurwitz triplet is a
triple of distinct Hurwitz curves with the identical automorphism
group of the lowest possible genus, namely $14$ (genera $3$ and
$7$ admit a unique Hurwitz curve, respectively the Klein quartic
curve (see \cite{K2}) and the Fricke-Macbeath curve (see \cite{F}
and \cite{M1})). It was studied by Shimura (see \cite{Sh}) and
Macbeath (see \cite{M2}). The explanation for this phenomenon is
arithmetic. Namely, in the ring of integers of the real subfield
${\Bbb Q}(\cos \frac{2 \pi}{7})$ of the cyclotomic field ${\Bbb
Q}(e^{\frac{2 \pi i}{7}})$, the rational prime $13$ splits as a
product of three distinct prime ideals. The principal congruence
subgroups defined by the triplet of primes produce Fuchsian groups
corresponding to the triplet of Riemann surfaces. Each of the
three Riemann surfaces in the first Hurwitz triplet can be formed
as a Fuchsian model, the quotient of the hyperbolic plane by one
of these three Fuchsian groups, which is just the model of Shimura
curves. We call it the standard arithmetic structure on the first
Hurwitz triplet.

  In the present paper, we find that the first Hurwitz triplet
possesses an exotic arithmetic structure. Consequently, we give a
different arithmetic explanation for this phenomenon. Our main
result is as follows:

{\bf Theorem 1.1. (Main Theorem).} {\it The first Hurwitz triplet
possesses two distinct arithmetic structures. The standard
arithmetic structure$:$ as Shimura curves $X_1$, $X_2$, $X_3$,
whose levels are with norm $13$. The exotic arithmetic
structure$:$ as non-congruence modular curves $Y_1$, $Y_2$, $Y_3$,
whose levels are $7$. Both of them are defined over ${\Bbb Q}(\cos
\frac{2 \pi}{7})$. However, for the third non-congruence modular
curve $Y_3$, there exist an ``exotic'' duality between the
associated non-congruence modular forms and the Hilbert modular
forms, both of them are related to ${\Bbb Q}(e^{\frac{2 \pi
i}{13}})!$}

  In order to explain Theorem 1.1, let us recall that the list of
all arithmetic subgroups of $PSL(2, {\Bbb R})$ is exhausted up to
commensurability by Fuchsian groups derived from quaternion
algebras over totally real number fields (see \cite{Ka}). If this
field is ${\Bbb Q}$ and the quaternion algebra is isomorphic to
$M_2({\Bbb Q})$, the full matrix algebra over ${\Bbb Q}$, then the
quotient space ${\Bbb H}/\Gamma$ is not compact but has finite
volume, and $\Gamma$ is commensurable with the modular group; in
all other cases ${\Bbb H}/\Gamma$ is compact. Now, we need the
following theorem, which is a hyperbolic and slightly non-standard
version of arithmeticity criteria which in essence are due to Weil
(see \cite{L}). Namely:

{\smc Proposition 1.1.} {\it Let $X$ be a smooth complex curve
which can be written as a quotient $X \simeq {\Bbb H}/G$, where $G
\subset PSL(2, {\Bbb R})$ is commensurable with a Fuchsian
triangular group. Then $X$ can be defined over a number field.}

  It may be useful to add Belyi theorem, which has to do with
the converse of the above statement. We state for clarity the
corresponding hyperbolic version.

{\smc Proposition 1.2. (Hyperbolic unramified version of Belyi
theorem).} {\it A smooth complex curve $X$ can be defined over a
number field if and only if there exists a finite set $Z \subset
X$ such that the affine curve $\check{X}=X \backslash Z$ is
uniformized by a Fuchsian group $G \subset PSL(2, {\Bbb R})$ with
$G$ commensurable to $PSL(2, {\Bbb Z})$.}

  Later on, we will use the other version of Belyi theorem.
Let us recall that one particular aspect of the modular group
$\Gamma=PSL(2, {\Bbb Z})$ is the balance (or rather the lack of
it) between its congruence and non-congruence subgroups (see
\cite{J}). Among the arithmetic subgroups (those of finite index)
in $\Gamma$, the congruence subgroups have proved to be the most
important and the most widely studied. Nevertheless, it has been
known for some time that, in a certain sense, most of the
arithmetic subgroups of $\Gamma$ are non-congruence subgroups. One
tends to regard congruence subgroups as, in some vague sense,
``known'', and it would make the study of $\Gamma$ a great deal
simpler if these were the only arithmetic subgroups. However, as
is often the case, low-dimensional behavior is subtle, the modular
group $\Gamma=PSL(2, {\Bbb Z})$ is truly exceptional, and it does
indeed possess non-congruence subgroups. By comparison, all finite
index subgroups of $PSL(n, {\Bbb Z})$, for $n \geq 3$, are
congruence! The arithmetic of non-congruence subgroups is very
interesting, because it is not encompassed in the Langlands
program per se. We will study the extraordinary phenomena
exhibited by the non-congruence subgroups and the connections with
representation theory and modular forms. On the other hand,
Shimura curves and their higher-dimensional analogue-Shimura
varieties are of central importance in the Langlands program. It
should be pointed out that when one pass from congruence modular
curves to Shimura curves, new difficulties arise due to the
absence of natural modular forms. In contrast with this fact, we
can associate non-congruence modular curves with non-congruence
modular forms. The importance of finite index subgroups of
$\Gamma$ comes from the following:

{\smc Proposition 1.3. (The other version of Belyi theorem).} {\it
Any smooth compact complex projective curve defined over
$\overline{{\Bbb Q}}$ can be realized $($in many ways$)$ as a
modular curve for a finite index subgroup of $\Gamma$.}

  Belyi's theorem tells us that, viewed simply as algebraic curves,
non-congruence modular curves are very general, and that we should
not expect them to have any special arithmetic properties. On the
other hand, the uniformization of non-congruence modular curves by
the upper half-plane is quite special, and leads to surprising
consequences. By Belyi theorem, Shimura curves over totally real
fields are also modular curves of some finite index, often
non-congruence, subgroups of $\Gamma$. In fact, for the Klein
quartic curve of genus three, Klein (see \cite{K2}) found that the
action of its automorphism group $PSL(2, 7)$ is defined over the
field ${\Bbb Q}(e^{\frac{2 \pi i}{7}})$. In his paper \cite{Sh},
Shimura proved that the model of Shimura curves and the model of
congruence modular curves are complex analytically isomorphic.
Moreover, he found that both models have ${\Bbb Q}(e^{\frac{2 \pi
i}{7}})$ as their natural field of definition. For the
Fricke-Macbeath curve of genus seven, besides the realization of
Shimura curves (see \cite{Sh}), Wohlfahrt (see \cite{W2}) gave the
realization of non-congruence modular curves with level seven.
Macbeath (see \cite{M1}) obtained the canonical model for this
curve, which is defined over the field ${\Bbb Q}(e^{\frac{2 \pi
i}{7}})$.

  In contrast to the Hurwitz curves with genera three and seven, we
find that the Hurwitz curves of genus fourteen possesses two
distinct arithmetic structures. This is the background of Theorem
1.1.

  In order to prove Theorem 1.1. we construct a six-dimensional
representation of the simple group $PSL(2, 13)$ of order $1092$ on
the five-dimensional projective space
$${\Bbb P}^5=\{ (z_1, z_2, z_3, z_4, z_5, z_6): z_i \in {\Bbb C} \quad (i=1, 2, 3, 4, 5, 6) \}.$$
In particular, this representation is defined over the cyclotomic
field ${\Bbb Q}(e^{\frac{2 \pi i}{13}})$. More precisely, put
$$S=-\frac{1}{\sqrt{13}} \left(\matrix
  \zeta^{12}-\zeta & \zeta^{10}-\zeta^3 & \zeta^4-\zeta^9 & \zeta^5-\zeta^8 & \zeta^2-\zeta^{11} & \zeta^6-\zeta^7\\
  \zeta^{10}-\zeta^3 & \zeta^4-\zeta^9 & \zeta^{12}-\zeta & \zeta^2-\zeta^{11} & \zeta^6-\zeta^7 & \zeta^5-\zeta^8\\
  \zeta^4-\zeta^9 & \zeta^{12}-\zeta & \zeta^{10}-\zeta^3 & \zeta^6-\zeta^7 & \zeta^5-\zeta^8 & \zeta^2-\zeta^{11}\\
  \zeta^5-\zeta^8 & \zeta^2-\zeta^{11} & \zeta^6-\zeta^7 & \zeta-\zeta^{12} & \zeta^3-\zeta^{10} & \zeta^9-\zeta^4\\
  \zeta^2-\zeta^{11} & \zeta^6-\zeta^7 & \zeta^5-\zeta^8 & \zeta^3-\zeta^{10} & \zeta^9-\zeta^4 & \zeta-\zeta^{12}\\
  \zeta^6-\zeta^7 & \zeta^5-\zeta^8 & \zeta^2-\zeta^{11} & \zeta^9-\zeta^4 & \zeta-\zeta^{12} & \zeta^3-\zeta^{10}
\endmatrix\right)\eqno{(1.1)}$$
and
$$T=\left(\matrix
  \zeta^7 &   &   &   &   &    \\
  & \zeta^{11} &  &   &   &    \\
  &    & \zeta^8  &   &   &    \\
  &    &     & \zeta^6    &    \\
  &    &     &  & \zeta^2 &    \\
  &    &     &        &   & \zeta^5
\endmatrix\right),\eqno{(1.2)}$$
where $\zeta=\exp(2 \pi i/13)$. Let $G=\langle S, T \rangle$. We
prove that $G \cong PSL(2, 13)$. Furthermore, we find three
different kinds of representations of $G$ as $(2, 3, 7)$-generated
groups, which correspond to the conjugacy classes $7A$, $7B$ and
$7C$ of $G$, respectively. Put
$$(2, 3, n; p):=\langle u, v: u^3=v^2=(uv)^n=(u^{-1} v^{-1} uv)^p=1 \rangle.$$
Now, we can give the results:

  The first model: $x_1^3=y_1^2=(x_1 y_1)^7=1$, where
$$x_1=-\frac{1}{\sqrt{13}} \left(\matrix
  \zeta^9-\zeta^{12} & \zeta^{11}-\zeta^7 & \zeta^8-\zeta^9 & \zeta^7-\zeta^5 & \zeta^4-\zeta^{11} & \zeta^4-\zeta^{12}\\
  \zeta^7-\zeta^3 & \zeta^3-\zeta^4 & \zeta^8-\zeta^{11} & \zeta^{10}-\zeta^4 & \zeta^{11}-\zeta^6 & \zeta^{10}-\zeta^8\\
  \zeta^7-\zeta^8 & \zeta^{11}-\zeta & \zeta-\zeta^{10} & \zeta^{12}-\zeta^7 & \zeta^{12}-\zeta^{10} & \zeta^8-\zeta^2\\
  \zeta^8-\zeta^6 & \zeta^2-\zeta^9 & \zeta-\zeta^9 & \zeta^4-\zeta & \zeta^2-\zeta^6 & \zeta^5-\zeta^4\\
  \zeta^9-\zeta^3 & \zeta^7-\zeta^2 & \zeta^5-\zeta^3 & \zeta^6-\zeta^{10} & \zeta^{10}-\zeta^9 & \zeta^5-\zeta^2\\
  \zeta^6-\zeta & \zeta^3-\zeta & \zeta^{11}-\zeta^5 & \zeta^6-\zeta^5 & \zeta^2-\zeta^{12} & \zeta^{12}-\zeta^3
\endmatrix\right),\eqno{(1.3)}$$
$$y_1=\left(\matrix
  0 & 0 & 0 & -\zeta & 0 & 0\\
  0 & 0 & 0 & 0 & -\zeta^9 & 0\\
  0 & 0 & 0 & 0 & 0 & -\zeta^3\\
  \zeta^{12} & 0 & 0 & 0 & 0 & 0\\
  0 & \zeta^4 & 0 & 0 & 0 & 0\\
  0 & 0 & \zeta^{10} & 0 & 0 & 0
\endmatrix\right).\eqno{(1.4)}$$
$$\langle x_1, y_1 \rangle \cong (2, 3, 7; 7).$$

  The second model: $x_2^3=y_2^2=(x_2 y_2)^7=1$, where
$$x_2=-\frac{1}{\sqrt{13}} \left(\matrix
  \zeta^9-\zeta^{10} & \zeta^5-\zeta^8 & \zeta-\zeta^{10} & \zeta^3-\zeta^{11} & \zeta^{11}-\zeta^9 & \zeta-\zeta^8\\
  \zeta^9-\zeta^{12} & \zeta^3-\zeta^{12} & \zeta^6-\zeta^7 & \zeta^9-\zeta^7 & \zeta-\zeta^8 & \zeta^8-\zeta^3\\
  \zeta^2-\zeta^{11} & \zeta^3-\zeta^4 & \zeta-\zeta^4 & \zeta^7-\zeta & \zeta^3-\zeta^{11} & \zeta^9-\zeta^7\\
  \zeta^2-\zeta^{10} & \zeta^4-\zeta^2 & \zeta^5-\zeta^{12} & \zeta^4-\zeta^3 & \zeta^8-\zeta^5 & \zeta^{12}-\zeta^3\\
  \zeta^6-\zeta^4 & \zeta^5-\zeta^{12} & \zeta^{10}-\zeta^5 & \zeta^4-\zeta & \zeta^{10}-\zeta & \zeta^7-\zeta^6\\
  \zeta^{12}-\zeta^6 & \zeta^2-\zeta^{10} & \zeta^6-\zeta^4 & \zeta^{11}-\zeta^2 & \zeta^{10}-\zeta^9 & \zeta^{12}-\zeta^9
\endmatrix\right),\eqno{(1.5)}$$
$$y_2=-\frac{1}{\sqrt{13}} \left(\matrix
  \zeta^7-\zeta^6 & \zeta^8-\zeta^5 & \zeta^{11}-\zeta^2 & \zeta^4-\zeta^9 & \zeta^{12}-\zeta & \zeta^{10}-\zeta^3\\
  \zeta^8-\zeta^5 & \zeta^{11}-\zeta^2 & \zeta^7-\zeta^6 & \zeta^{12}-\zeta & \zeta^{10}-\zeta^3 & \zeta^4-\zeta^9\\
  \zeta^{11}-\zeta^2 & \zeta^7-\zeta^6 & \zeta^8-\zeta^5 & \zeta^{10}-\zeta^3 & \zeta^4-\zeta^9 & \zeta^{12}-\zeta\\
  \zeta^4-\zeta^9 & \zeta^{12}-\zeta & \zeta^{10}-\zeta^3 & \zeta^6-\zeta^7 & \zeta^5-\zeta^8 & \zeta^2-\zeta^{11}\\
  \zeta^{12}-\zeta & \zeta^{10}-\zeta^3 & \zeta^4-\zeta^9 & \zeta^5-\zeta^8 & \zeta^2-\zeta^{11} & \zeta^6-\zeta^7\\
  \zeta^{10}-\zeta^3 & \zeta^4-\zeta^9 & \zeta^{12}-\zeta & \zeta^2-\zeta^{11} & \zeta^6-\zeta^7 & \zeta^5-\zeta^8
\endmatrix\right).\eqno{(1.6)}$$
$$\langle x_2, y_2 \rangle \cong (2, 3, 7; 6).$$

  The third model: $x_3^3=y_3^2=(x_3 y_3)^7=1$, where
$$x_3=-\frac{1}{\sqrt{13}} \left(\matrix
  \zeta^{12}-\zeta^3 & \zeta^6-\zeta^5 & \zeta^2-\zeta^{12} & \zeta^5-\zeta^{11} & \zeta-\zeta^6 & \zeta-\zeta^3\\
  \zeta^5-\zeta^4 & \zeta^4-\zeta & \zeta^2-\zeta^6 & \zeta^9-\zeta & \zeta^6-\zeta^8 & \zeta^9-\zeta^2\\
  \zeta^5-\zeta^2 & \zeta^6-\zeta^{10} & \zeta^{10}-\zeta^9 & \zeta^3-\zeta^5 & \zeta^3-\zeta^9 & \zeta^2-\zeta^7\\
  \zeta^2-\zeta^8 & \zeta^7-\zeta^{12} & \zeta^{10}-\zeta^{12} & \zeta-\zeta^{10} & \zeta^7-\zeta^8 & \zeta^{11}-\zeta\\
  \zeta^{12}-\zeta^4 & \zeta^5-\zeta^7 & \zeta^{11}-\zeta^4 & \zeta^8-\zeta^9 & \zeta^9-\zeta^{12} & \zeta^{11}-\zeta^7\\
  \zeta^8-\zeta^{10} & \zeta^4-\zeta^{10} & \zeta^6-\zeta^{11} & \zeta^8-\zeta^{11} & \zeta^7-\zeta^3 & \zeta^3-\zeta^4
\endmatrix\right),\eqno{(1.7)}$$
$$y_3=-\frac{1}{\sqrt{13}} \left(\matrix
  \zeta^8-\zeta^5 & \zeta^4-\zeta^8 & \zeta^2-\zeta & \zeta^4-\zeta^6 & \zeta^9-\zeta^2 & \zeta-\zeta^6\\
  \zeta^5-\zeta^9 & \zeta^7-\zeta^6 & \zeta^{10}-\zeta^7 & \zeta^9-\zeta^2 & \zeta^{10}-\zeta^2 & \zeta^3-\zeta^5\\
  \zeta^{12}-\zeta^{11} & \zeta^6-\zeta^3 & \zeta^{11}-\zeta^2 & \zeta-\zeta^6 & \zeta^3-\zeta^5 & \zeta^{12}-\zeta^5\\
  \zeta^7-\zeta^9 & \zeta^{11}-\zeta^4 & \zeta^7-\zeta^{12} & \zeta^5-\zeta^8 & \zeta^9-\zeta^5 & \zeta^{11}-\zeta^{12}\\
  \zeta^{11}-\zeta^4 & \zeta^{11}-\zeta^3 & \zeta^8-\zeta^{10} & \zeta^8-\zeta^4 & \zeta^6-\zeta^7 & \zeta^3-\zeta^6\\
  \zeta^7-\zeta^{12} & \zeta^8-\zeta^{10} & \zeta^8-\zeta & \zeta-\zeta^2 & \zeta^7-\zeta^{10} & \zeta^2-\zeta^{11}
\endmatrix\right).\eqno{(1.8)}$$
$$\langle x_3, y_3 \rangle \cong (2, 3, 7; 13) \quad \text{with an extra condition (see Theorem 3.3)}.$$
Here, $z_1=x_1^{-1} y_1^{-1}$, $z_2=x_2^{-1} y_2^{-1}$ and
$z_3=x_3^{-1} y_3^{-1}$ correspond to the conjugacy classes $7A$,
$7B$ and $7C$, respectively. Moreover,
$$\langle x_1, y_1 \rangle=\langle x_2, y_2 \rangle=\langle x_3, y_3 \rangle=G.\eqno{(1.9)}$$

  It is well-known that the modular group $\Gamma=PSL(2, {\Bbb Z})$
is generated by the following linear fractional transformations:
$T \tau=\tau+1$, $S \tau=-1/\tau$. Let $P=ST$. Then $P
\tau=-1/(\tau+1)$. Here, $S^2=P^3=1$. Let $\phi_i: \Gamma \to
PGL(6, {\Bbb C})$ be three representations where
$$\phi_i: S \mapsto y_i, \quad P \mapsto x_i, \quad T^{-1} \mapsto z_i, \quad (i=1, 2, 3).\eqno{(1.10)}$$
Let $Y_i=\overline{{\Bbb H}/G_i}$ be the compactification of
${\Bbb H}/G_i$ where $G_i=\text{ker} \phi_i$. Then
$$\Gamma/G_1 \cong \Gamma/G_2 \cong \Gamma/G_3 \cong PSL(2, 13).\eqno{(1.11)}$$
$G_1$, $G_2$ and $G_3$ are non-congruence normal subgroups of
level $7$ of $\Gamma$. As a $(2, 3, 7)$-generated group, by
Riemann-Hurwitz formula, we have
$$2g-2=1092 \left(1-\frac{1}{2}-\frac{1}{3}-\frac{1}{7}\right).\eqno{(1.12)}$$
Hence $g=14$, which is of genus of ${\Bbb H}/G_i$ for $i=1, 2, 3$.
Therefore, $Y_1$, $Y_2$ and $Y_3$ must be Hurwitz curves. By
Macbeath's theorem (see Theorem 3.7), there are only three Hurwitz
curves with genus $14$. Hence, $Y_1, Y_2, Y_3$ must be complex
analytically isomorphic to $X_1, X_2, X_3$. Note that $X_1, X_2,
X_3$ correspond to the conjugacy classes $7A$, $7B$, $7C$,
respectively, and $Y_1, Y_2, Y_3$ also correspond to the conjugacy
classes $7A$, $7B$, $7C$, respectively. This implies that $Y_i$ is
complex analytically isomorphic to $X_i$ for $i=1, 2, 3$. By the
result of Streit in \cite{St}, we show that the minimal fields of
definition for $Y_1$, $Y_2$ and $Y_3$ are ${\Bbb Q}(\cos \frac{2
\pi}{7})$.

  Let
$$Q=-\frac{1}{\sqrt{13}} \left(\matrix
  \zeta^7-\zeta^9 & \zeta^4-\zeta^{10} & \zeta^2-\zeta^7 & \zeta^{10}-1 & \zeta^8-\zeta^4 & \zeta^8-\zeta^9\\
  \zeta^5-\zeta^{11} & \zeta^{11}-\zeta^3 & \zeta^{10}-\zeta^{12} & \zeta^7-\zeta^3 & \zeta^{12}-1 & \zeta^7-\zeta^{10}\\
  \zeta^{12}-\zeta^4 & \zeta^6-\zeta^8 & \zeta^8-\zeta & \zeta^{11}-\zeta^{12} & \zeta^{11}-\zeta & \zeta^4-1\\
  1-\zeta^3 & \zeta^9-\zeta^5 & \zeta^4-\zeta^5 & \zeta^6-\zeta^4 & \zeta^9-\zeta^3 & \zeta^{11}-\zeta^6\\
  \zeta^{10}-\zeta^6 & 1-\zeta & \zeta^3-\zeta^6 & \zeta^8-\zeta^2 & \zeta^2-\zeta^{10} & \zeta^3-\zeta\\
  \zeta-\zeta^2 & \zeta^{12}-\zeta^2 & 1-\zeta^9 & \zeta-\zeta^9 & \zeta^7-\zeta^5 & \zeta^5-\zeta^{12}
 \endmatrix\right).\eqno{(1.13)}$$
We find that $Q^7=1$ and $x_3 y_3=Q^3$. Put
$${\Bbb A}_0(z_1, z_2, z_3, z_4, z_5, z_6)=z_1 z_4+z_2 z_5+z_3 z_6.\eqno{(1.14)}$$
Let $H=y_2 \cdot S$. Then $\langle H, T \rangle \cong {\Bbb
Z}_{13} \rtimes {\Bbb Z}_6$. Hence, it is the maximal subgroup of
order $78$ of $G$ with index $14$. We find that ${\Bbb A}_0^2$ is
invariant under the action of the maximal subgroup $\langle H, T
\rangle$.

{\bf Theorem 1.2. (Exotic duality theorem)}. {\it There exists an
``exotic'' duality between non-congruence modular forms and
Hilbert modular forms, both of them are related to ${\Bbb
Q}(\zeta):$
$$\left\{\matrix
  Q^{\nu}({\Bbb A}_0)^2, y_3 Q^{\nu}({\Bbb A}_0)^2\\
  \text{$\nu$ mod $7$}
  \endmatrix\right\}=
  \left\{\matrix
  T^{\nu}({\Bbb A}_0)^2, S T^{\nu}({\Bbb A}_0)^2\\
  \text{$\nu$ mod $13$}
\endmatrix\right\}\eqno{(1.15)}$$}

  The left hand side of (1.15) corresponds to $(2, 3, 7; 13)$,
the right hand side of (1.15) corresponds to $(2, 3, 13; 7)$. When
we exchange the position of $7$ and $13$, we will get the duality
between the two representations of $G$. $T^{\nu}({\Bbb A}_0)^2, S
T^{\nu}({\Bbb A}_0)^2$ for $\nu$ mod $13$ give a coset
decomposition of $G$ with respect to the maximal subgroup $\langle
H, T \rangle$. Therefore, any symmetric polynomial of ${\Bbb
A}_0^2$, $ST^{\nu}({\Bbb A}_0)^2$ for $\nu$ mod $13$ gives an
invariant polynomial under the action of $PSL(2, 13)$. By (1.15),
we have that any symmetric polynomial of $Q^{\nu}({\Bbb A}_0)^2,
y_3 Q^{\nu}({\Bbb A}_0)^2$ for $\nu$ mod $7$ also gives an
invariant polynomial under the action of $PSL(2, 13)$. It is known
that invariant polynomials correspond to modular forms. The right
hand side of (1.15) corresponds to the congruence subgroup
$\Gamma_{13, 1}$ with type $(1, 13)$, the associated modular forms
are Hilbert modular forms. However, the left hand side of (1.15)
corresponds to the non-congruence subgroup $\Gamma_{7, 7}$ with
type $(7, 7)$, the associated modular forms are non-congruence
modular forms! Therefore, (1.15) gives a connection between
congruence modular forms (Hilbert modular forms) and
non-congruence modular forms, both of them are related to the
cyclotomic field ${\Bbb Q}(\zeta)$! This complete the proof of
Theorem 1.1.

  Our main results have relations and applications to modular
equations of degree fourteen (including Jacobian modular equation
and ``exotic'' modular equation), ``triality'' of the
representation of $PSL(2, 13)$, Haagerup subfactor, geometry of
the exceptional Lie group $G_2$, and even the Monster finite
simple group ${\Bbb M}$!

  It is well-known that the genus of the modular curve
$X_0(p)={\Bbb H}/\Gamma_0(p)$ for prime $p$ is zero if and only if
$p=2, 3, 5, 7, 13$. In his paper \cite{K1}, Klein studied the
modular equations of orders $2$, $3$, $5$, $7$, $13$ with degrees
$3$, $4$, $6$, $8$, $14$, respectively. They are the so-called
Hauptmoduln (principal moduli). In \cite{K}, Klein investigated
the modular equation of degree six in connection with the Jacobian
equation of degree six, which can be used to solve the general
quintic equation. In \cite{K2} and \cite{K3}, he studied the
modular equation of degree eight in connection with the Jacobian
equation of degree eight, which can be used to solve the algebraic
equation of degree eight with Galois group $PSL(2, 7)$. However,
he did not investigate the last case, the modular equation of
degree fourteen. In this paper, we find the Jacobian equation of
degree fourteen, which corresponds to the modular equation of
degree fourteen. Surprisingly, we find the other ``exotic''
modular equation of degree fourteen, which does not correspond to
the Jacobian equation! As an application, we obtain the following
quartic four-fold
$$(z_3 z_4^3+z_1 z_5^3+z_2 z_6^3)-(z_6 z_1^3+z_4 z_2^3+z_5 z_3^3)+
  3(z_1 z_2 z_4 z_5+z_2 z_3 z_5 z_6+z_3 z_1 z_6 z_4)=0,\eqno{(1.16)}$$
which is invariant under the action of the simple group $G$. It is
a higher-dimensional counterpart of the Klein quartic curve (see
\cite{K2}) and the Klein cubic threefold (see \cite{K4}).

  Note that the Hurwitz curves of genus $14$ are are non-hyperelliptic.
This leads us to study their canonical models in ${\Bbb P}^{13}$
and the corresponding fourteen-dimensional representation. We
construct such a representation which is induced from our
six-dimensional representation.

  The Jones polynomials arise from the study of subfactors.
Subfactor theory is among operator algebra theory one of the most
influential subjects in the sense that it has relations and
applications to topology, representation theory, and quantum
physics. Some subfactors give rise to some quantum invariants of
knots and links. Haagerup discovered subfactors that are not
associated to any other known objects like quantum groups (see
\cite{AH}). According to \cite{J2} and \cite{J1}:`` These sporadic
subfactors discovered by Haagerup and others remain exotic
creatures in the zoo, untouched by quantum groups or any non-von
Neumann algebras approach. Hence, it is a major challenge in the
theory to come up with an interpretation of these subfactors as
members of a family related to some other mathematical objects.''
As an application, we give a connection between the Haagerup
modular data and our seven-dimensional representation of $PSL(2,
13)$ which has the ``triality''!

  The second application is concerned with geometry of the exceptional
Lie group $G_2$. According to \cite{B2}: `` There are very
beautiful conjectural relations that the geometry of $G_2$ have
ties in automorphic forms and all kinds of special geometry that,
if the physicists were right, are just waiting there to be seen.
But at the moment our techniques are far too weak to be able to
verify that claim''. We find that our invariant quadric
$${\Bbb A}_0^2+{\Bbb A}_1 {\Bbb A}_5+{\Bbb A}_2 {\Bbb A}_3+{\Bbb A}_4 {\Bbb A}_6=0\eqno{(1.17)}$$
is just Cartan's ${\Bbb Q}_5 \subset {\Bbb P}^6$, the space of
$J$-null lines in ${\Bbb C}^7$, the quotient variety $G_2/P_1$
which is isomorphic to a $5$-dimensional quadric $Q^5 \subset
{\Bbb P}^6$ and the split Cayley hexagon over ${\Bbb C}$! This
shows that these three geometric objects are modular,
parameterized by both Hilbert modular forms and con-congruence
modular forms related to ${\Bbb Q}(\zeta)$.

  The third application is the relation between genus-zero subgroups (both
congruence subgroups and non-congruence subgroups) of the modular
group and subgroups of the Monster finite simple group ${\Bbb M}$.
In particular, we are interested in $\Gamma_0(13)$ and $\Gamma_{7,
7}$, the quotient groups $\Gamma/\Gamma(13) \cong \Gamma/G_i \cong
PSL(2, 13)$ for $i=1, 2, 3$. We construct the following variety
(see section six):
$$M:=7 \cdot 13^2 {\Bbb G}_0^2+{\Bbb G}_1 {\Bbb G}_{12}+{\Bbb G}_ 2
     {\Bbb G}_{11}+\cdots+{\Bbb G}_6 {\Bbb G}_7=0.\eqno{(1.18)}$$
It is a twelve-dimensional $PSL(2, 13)$-invariant complex
algebraic variety (i.e., twenty-four-dimensional manifold) of
degree four in the projective space
$${\Bbb P}^{13}=\{({\Bbb D}_0, {\Bbb D}_1, \cdots, {\Bbb D}_{12}, {\Bbb D}_{\infty})\},$$
which is related to $N(13B^2)$, a maximal $p$-local subgroup of
${\Bbb M}$.

  This paper consists of six sections. Section two, section
three and section four are devoted to the proof of Theorem 1.1.
Namely, in section two, we give the standard arithmetic structure
on the first Hurwitz triplet. In section three, we give a
six-dimensional representation of $PSL(2, 13)$ defined over ${\Bbb
Q}(e^{\frac{2 \pi i}{13}})$ and the exotic arithmetic structure on
the first Hurwitz triplet. In section four, we give a
seven-dimensional representation of $PSL(2, 13)$ which induces
from our six-dimensional representation and find the ``triality''
of the representation of $PSL(2, 13)$. Consequently, we obtain
Jacobian modular equation of degree fourteen and find the
``exotic'' duality between non-congruence modular forms and
Hilbert modular forms, both of them are related to ${\Bbb
Q}(\zeta)$. In section five, we give applications of our results
to the Haagerup subfactor and the geometry associated to the
exceptional Lie group $G_2$. In section six, we give a
fourteen-dimensional representation of $PSL(2, 13)$ associated to
the canonical model of the first Hurwitz triplet in ${\Bbb
P}^{13}$, which also induces from our six-dimensional
representation. By this representation, we find the ``exotic''
modular equation of degree fourteen. Consequently, we construct a
twelve-dimensional $PSL(2, 13)$-invariant complex algebraic
variety (i.e., $24$-dimensional manifold) of degree four which is
related to a maximal $p$-local subgroup of the Monster simple
group ${\Bbb M}$.

\vskip 0.5 cm

\centerline{\bf 2. Standard arithmetic structure on the first
                   Hurwitz triplet}

\vskip 0.5 cm

  The existence of a quaternion algebra presentation for Hurwitz
curves is due to Shimura \cite{Sh}. An explicit order was briefly
described by Elkies in \cite{El1} and in \cite{El2}. We follow the
concrete realization of $G_{2, 3, 7}$ in terms of the group of
elements of norm one in an order of a quaternion algebra, given by
Elkies in \cite{El1} and \cite{El2} (see \cite{KSV}).

  In this section, let $K$ denote the real subfield of ${\Bbb Q}(\rho)$,
where $\rho$ is a primitive seventh root of unity. Thus $K={\Bbb
Q}(\eta)$, where the element $\eta=\rho+\rho^{-1}$ satisfies the
relation
$$\eta^3+\eta^2-2 \eta-1=0.$$
There are three embeddings of $K$ into ${\Bbb R}$, defined by
sending $\eta$ to any of the three real roots of the above
equation, namely
$$2 \cos \frac{2 \pi}{7}, \quad 2 \cos \frac{4 \pi}{7}, \quad
  2 \cos \frac{6 \pi}{7}.$$
We view the first embedding as the natural one $K \hookrightarrow
{\Bbb R}$, and denote the others by $\sigma_1, \sigma_2: K \to
{\Bbb R}$. Note that $2 \cos \frac{2 \pi}{7}$ is a positive root,
while the other two are negative.

  Let $D$ be the quaternion $K$-algebra $K(i, j)$ with $i^2=j^2=\eta$, $ji=-ij$.
Let ${\Cal O} \subset D$ be the order defined by ${\Cal O}={\Cal
O}_K[i, j]$, where ${\Cal O}_K$ is the ring of integers in $K$.
Then ${\Cal O} \cong {\Bbb Z}[\eta][i, j]$. Fix the element
$\tau=1+\eta+\eta^2$, and define an element $j^{\prime} \in D$ by
setting $j^{\prime}=\frac{1}{2}(1+\eta i+\tau j)$. We define a new
order ${\Cal Q}_{\text{Hur}} \subset D$ by setting
$${\Cal Q}_{\text{Hur}}={\Bbb Z}[\eta][i, j, j^{\prime}].$$
The group of elements of norm $1$ in the order ${\Cal
Q}_{\text{Hur}}$, modulo the center $\{ \pm 1 \}$, is isomorphic
to the $(2, 3, 7)$ group. Indeed, Elkies gave the elements
$$\aligned
  g_2 &=\frac{1}{\eta} ij,\\
  g_3 &=\frac{1}{2} \left[1+(\eta^2-2)j+(3-\eta^2)ij\right],\\
  g_7 &=\frac{1}{2} \left[(\tau-2)+(2-\eta^2)i+(\tau-3)ij\right],
\endaligned$$
satisfying the relations $g_2^2=g_3^3=g_7^7=-1$ and $g_2=g_7 g_3$,
which therefore project to generators of $G_{2, 3, 7} \subset
PSL(2, {\Bbb R})$.

  Now, we can give the model of Shimura curves for the first Hurwitz
triplet (see \cite{KSV}). In fact, one has
$$13=\eta (\eta+2) (2 \eta-1)(3-2 \eta)(\eta+3),$$
where $\eta (\eta+2)$ is invertible. Hence,
$$13 {\Cal O}_K=(2 \eta-1) {\Cal O}_K \cdot (3-2 \eta) {\Cal O}_K
                \cdot (\eta+3) {\Cal O}_K. \eqno{(2.1)}$$
The three prime ideals define a triplet of principal congruence
subgroups as follows:
$${\Cal Q}_{\text{Hur}}^1(I)=\{ x \in {\Cal Q}_{\text{Hur}}^1:
  x \equiv 1 (\text{mod $I {\Cal Q}_{\text{Hur}}$}) \},$$
where $I \subset {\Cal O}_K$ is an ideal and ${\Cal
Q}_{\text{Hur}}^1$ is the group of elements of norm $1$ in ${\Cal
Q}_{\text{Hur}}$. One therefore obtains a triplet of distinct
Hurwitz curves of genus $14$. All of them have the identical
automorphism group $PSL(2, 13)$ realized as the quotient
$${\Cal Q}_{\text{Hur}}^1/{\Cal Q}_{\text{Hur}}^1(I) \cong PSL(2, 13),$$
whose actions are defined over the real subfield ${\Bbb Q}(\cos
\frac{2 \pi}{7})$ of the cyclotomic field ${\Bbb Q}(e^{\frac{2 \pi
i}{7}})$.

\vskip 0.5 cm

\centerline{\bf 3. Six-dimensional representations of $PSL(2, 13)$
                   and}

\centerline{\bf exotic arithmetic structure on the first
                Hurwitz triplet}

\vskip 0.5 cm

  Let us recall the Weil representation of $SL_2({\Bbb F}_p)$
(see \cite{W}) which is modelled on the space $L^2({\Bbb F}_p)$ of
all complex valued functions on ${\Bbb F}_p$ where $p$ is a prime.
Denote by $L^2({\Bbb F}_p)$ the $p$-dimensional complex vector
space of all (square-integrable) complex valued functions on
${\Bbb F}_p$ with respect to counting measure, that is, all
functions from ${\Bbb F}_p$ to the complex numbers. We can
decompose $L^2({\Bbb F}_p)$ as the direct sum of the space $V^{+}$
of even functions and the space $V^{-}$ of odd functions. The
space $V^{-}$ has dimension $(p-1)/2$ and its associated
projective space ${\Bbb P}(V^{-})$ has dimension $(p-3)/2$. In
particular, for $p=13$, the Weil representation for $SL(2, 13)$ is
given as follows:
$$s=\left(\matrix 0 & -1\\ 1 & 0 \endmatrix\right), \quad
  t=\left(\matrix 1 & 1\\ 0 & 1 \endmatrix\right), \quad
  h=\left(\matrix 7 & 0\\ 0 & 2 \endmatrix\right).$$

  In this section, we will study the six-dimensional representation
of the simple group $PSL(2, 13)$ of order $1092$, which acts on
the five-dimensional projective space
$${\Bbb P}^5=\{ (z_1, z_2, z_3, z_4, z_5, z_6): z_i \in {\Bbb C} \quad (i=1, 2, 3, 4, 5, 6) \}.$$

  Let $\zeta=\exp(2 \pi i/13)$ and
$$\left\{\aligned
  &\theta_1=\zeta+\zeta^3+\zeta^9,\\
  &\theta_2=\zeta^2+\zeta^6+\zeta^5,\\
  &\theta_3=\zeta^4+\zeta^{12}+\zeta^{10},\\
  &\theta_4=\zeta^8+\zeta^{11}+\zeta^7.
\endaligned\right.\eqno{(3.1)}$$
We find that
$$\left\{\aligned
  &\theta_1+\theta_2+\theta_3+\theta_4=-1,\\
  &\theta_1 \theta_2+\theta_1 \theta_3+\theta_1 \theta_4+\theta_2 \theta_3+\theta_2 \theta_4+\theta_3 \theta_4=2,\\
  &\theta_1 \theta_2 \theta_3+\theta_1 \theta_2 \theta_4+\theta_1 \theta_3 \theta_4+\theta_2 \theta_3 \theta_4=4,\\
  &\theta_1 \theta_2 \theta_3 \theta_4=3.
\endaligned\right.$$
Hence, $\theta_1$, $\theta_2$, $\theta_3$ and $\theta_4$ satisfy
the quartic equation
$$z^4+z^3+2 z^2-4z+3=0,$$
which can be decomposed as two quadratic equations
$$\left(z^2+\frac{1+\sqrt{13}}{2} z+\frac{5+\sqrt{13}}{2}\right)
  \left(z^2+\frac{1-\sqrt{13}}{2} z+\frac{5-\sqrt{13}}{2}\right)=0$$
over the real quadratic field ${\Bbb Q}(\sqrt{13})$. Therefore,
the four roots are given as follows:
$$z_{1, 2}=\frac{1}{2}\left(-\frac{1+\sqrt{13}}{2} \pm \sqrt{\frac{-13-3 \sqrt{13}}{2}}\right)
          =\frac{1}{4} \left(-1-\sqrt{13} \pm \sqrt{-26-6 \sqrt{13}}\right),$$
$$z_{3, 4}=\frac{1}{2}\left(-\frac{1-\sqrt{13}}{2} \pm \sqrt{\frac{-13+3 \sqrt{13}}{2}}\right)
          =\frac{1}{4} \left(-1+\sqrt{13} \pm \sqrt{-26+6
\sqrt{13}}\right).$$ Note that
$$\text{Re}(\theta_1)=\cos \frac{2 \pi}{13}+\cos \frac{6 \pi}{13}-\cos \frac{5 \pi}{13}>0, \quad
  \text{Im}(\theta_1)=\sin \frac{2 \pi}{13}+\sin \frac{6 \pi}{13}-\sin \frac{5 \pi}{13}>0.$$
We have
$$\theta_1=\frac{1}{4} \left(-1+\sqrt{13}+\sqrt{-26+6 \sqrt{13}}\right).$$
Similarly,
$$\text{Re}(\theta_2)<0, \quad \text{Im}(\theta_2)>0, \quad
  \theta_2=\frac{1}{4} \left(-1-\sqrt{13}+\sqrt{-26-6 \sqrt{13}}\right).$$
$$\text{Re}(\theta_3)>0, \quad \text{Im}(\theta_3)<0, \quad
  \theta_3=\frac{1}{4} \left(-1+\sqrt{13}-\sqrt{-26+6 \sqrt{13}}\right).$$
$$\text{Re}(\theta_4)<0, \quad \text{Im}(\theta_4)<0, \quad
  \theta_4=\frac{1}{4} \left(-1-\sqrt{13}-\sqrt{-26-6 \sqrt{13}}\right).$$
Moreover, we find that
$$\left\{\aligned
  \theta_1+\theta_3+\theta_2+\theta_4 &=-1,\\
  \theta_1+\theta_3-\theta_2-\theta_4 &=\sqrt{13},\\
  \theta_1-\theta_3-\theta_2+\theta_4 &=-\sqrt{-13+2 \sqrt{13}},\\
  \theta_1-\theta_3+\theta_2-\theta_4 &=\sqrt{-13-2 \sqrt{13}}.
\endaligned\right.\eqno{(3.2)}$$

  Let
$$M=\left(\matrix
    \zeta-\zeta^{12} & \zeta^3-\zeta^{10} & \zeta^9-\zeta^4\\
    \zeta^3-\zeta^{10} & \zeta^9-\zeta^4 & \zeta-\zeta^{12}\\
    \zeta^9-\zeta^4 & \zeta-\zeta^{12} & \zeta^3-\zeta^{10}
   \endmatrix\right), \quad
  N=\left(\matrix
    \zeta^5-\zeta^8 & \zeta^2-\zeta^{11} & \zeta^6-\zeta^7\\
    \zeta^2-\zeta^{11} & \zeta^6-\zeta^7 & \zeta^5-\zeta^8\\
    \zeta^6-\zeta^7 & \zeta^5-\zeta^8 & \zeta^2-\zeta^{11}
   \endmatrix\right).\eqno{(3.3)}$$
Then $MN=NM=-\sqrt{13} I$ and $M^2+N^2=-13 I$. Put
$$S=-\frac{1}{\sqrt{13}} \left(\matrix
    -M & N\\
     N & M
    \endmatrix\right).\eqno{(3.4)}$$
Then $S^2=I$. In fact,
$$S=-\frac{2 i}{\sqrt{13}} \left(\matrix
    -\sin \frac{2 \pi}{13} & -\sin \frac{6 \pi}{13} &  \sin \frac{5 \pi}{13} &
     \sin \frac{3 \pi}{13} &  \sin \frac{4 \pi}{13} &  \sin \frac{\pi}{13}\\
    -\sin \frac{6 \pi}{13} &  \sin \frac{5 \pi}{13} & -\sin \frac{2 \pi}{13} &
     \sin \frac{4 \pi}{13} &  \sin \frac{\pi}{13}   &  \sin \frac{3 \pi}{13}\\
     \sin \frac{5 \pi}{13} & -\sin \frac{2 \pi}{13} & -\sin \frac{6 \pi}{13} &
     \sin \frac{\pi}{13}   &  \sin \frac{3 \pi}{13} &  \sin \frac{4 \pi}{13}\\
     \sin \frac{3 \pi}{13} &  \sin \frac{4 \pi}{13} &  \sin \frac{\pi}{13} &
     \sin \frac{2 \pi}{13} &  \sin \frac{6 \pi}{13} & -\sin \frac{5 \pi}{13}\\
     \sin \frac{4 \pi}{13} &  \sin \frac{\pi}{13}   &  \sin \frac{3 \pi}{13} &
     \sin \frac{6 \pi}{13} & -\sin \frac{5 \pi}{13} &  \sin \frac{2 \pi}{13}\\
     \sin \frac{\pi}{13}   &  \sin \frac{3 \pi}{13} &  \sin \frac{4 \pi}{13} &
    -\sin \frac{5 \pi}{13} &  \sin \frac{2 \pi}{13} &  \sin \frac{6 \pi}{13}
   \endmatrix\right).\eqno{(3.5)}$$
Note that
$$\frac{\sin \frac{2 \pi}{13} \sin \frac{5 \pi}{13} \sin \frac{6 \pi}{13}}
       {\sin \frac{\pi}{13} \sin \frac{3 \pi}{13} \sin \frac{4 \pi}{13}}
 =\frac{-\sqrt{\frac{-13-3 \sqrt{13}}{2}}}{-\sqrt{\frac{-13+3 \sqrt{13}}{2}}}
 =\frac{3+\sqrt{13}}{2},$$
which is a fundamental unit of ${\Bbb Q}(\sqrt{13})$. Let
$$T=\text{diag}(\zeta^7, \zeta^{11}, \zeta^8, \zeta^6, \zeta^2, \zeta^5).\eqno{(3.6)}$$

{\bf Theorem 3.1.} {\it Let $G=\langle S, T \rangle$. Then $G
\cong PSL(2, 13)$.}

{\it Proof}. We have
$$ST=-\frac{1}{\sqrt{13}} \left(\matrix
  \zeta^6-\zeta^8 & \zeta^8-\zeta & \zeta^{12}-\zeta^4 & \zeta^{11}-\zeta & \zeta^4-1 & \zeta^{11}-\zeta^{12}\\
  \zeta^4-\zeta^{10} & \zeta^2-\zeta^7 & \zeta^7-\zeta^9 & \zeta^8-\zeta^4 & \zeta^8-\zeta^9 & \zeta^{10}-1\\
  \zeta^{11}-\zeta^3 & \zeta^{10}-\zeta^{12} & \zeta^5-\zeta^{11} & \zeta^{12}-1 & \zeta^7-\zeta^{10} & \zeta^7-\zeta^3\\
  \zeta^{12}-\zeta^2 & 1-\zeta^9 & \zeta-\zeta^2 & \zeta^7-\zeta^5 & \zeta^5-\zeta^{12} & \zeta-\zeta^9\\
  \zeta^9-\zeta^5 & \zeta^4-\zeta^5 & 1-\zeta^3 & \zeta^9-\zeta^3 & \zeta^{11}-\zeta^6 & \zeta^6-\zeta^4\\
  1-\zeta & \zeta^3-\zeta^6 & \zeta^{10}-\zeta^6 & \zeta^2-\zeta^{10} & \zeta^3-\zeta & \zeta^8-\zeta^2
\endmatrix\right).\eqno{(3.7)}$$
On the other hand,
$$\aligned
 &(ST)^{-1}=T^{-1} S\\
=&-\frac{1}{\sqrt{13}} \left(\matrix
  \zeta^5-\zeta^7 & \zeta^3-\zeta^9 & \zeta^{10}-\zeta^2 & \zeta^{11}-\zeta & \zeta^8-\zeta^4 & \zeta^{12}-1\\
  \zeta^{12}-\zeta^5 & \zeta^6-\zeta^{11} & \zeta-\zeta^3 & \zeta^4-1 & \zeta^8-\zeta^9 & \zeta^7-\zeta^{10}\\
  \zeta^9-\zeta & \zeta^4-\zeta^6 & \zeta^2-\zeta^8 & \zeta^{11}-\zeta^{12} & \zeta^{10}-1 & \zeta^7-\zeta^3\\
  \zeta^{12}-\zeta^2 & \zeta^9-\zeta^5 & 1-\zeta & \zeta^8-\zeta^6 & \zeta^{10}-\zeta^4 & \zeta^3-\zeta^{11}\\
  1-\zeta^9 & \zeta^4-\zeta^5 & \zeta^3-\zeta^6 & \zeta-\zeta^8 & \zeta^7-\zeta^2 & \zeta^{12}-\zeta^{10}\\
  \zeta-\zeta^2 & 1-\zeta^3 & \zeta^{10}-\zeta^6 & \zeta^4-\zeta^{12} & \zeta^9-\zeta^7 & \zeta^{11}-\zeta^5
\endmatrix\right).
\endaligned$$
We will calculate $(ST)^2=\frac{1}{13} (a_{ij})$. Without loss of
generality, we will compute $a_{1i}$, $i=1, 2, 3, 4, 5, 6$. Here,
$$a_{11}=-2 \zeta-2 \zeta^2+2 \zeta^3-2 \zeta^9+2 \zeta^{10}+2\zeta^{11}-\zeta^5+\zeta^7.$$
By the identity
$$\sqrt{13}=\zeta+\zeta^{12}+\zeta^3+\zeta^{10}+\zeta^9+\zeta^4-\zeta^5-\zeta^8-\zeta^2-\zeta^{11}-\zeta^6-\zeta^7,\eqno{(3.8)}$$
we have
$$(\zeta^5-\zeta^7) \sqrt{13}
 =2 \zeta+2 \zeta^2-2 \zeta^3+2 \zeta^9-2 \zeta^{10}-2 \zeta^{11}+\zeta^5-\zeta^7.$$
Hence, $a_{11}=-(\zeta^5-\zeta^7) \sqrt{13}$. Similarly,
$$\aligned
  a_{12}&=-2 \zeta^2-2 \zeta^4+2 \zeta^5-2 \zeta^7+2 \zeta^8+2 \zeta^{10}-\zeta^3+\zeta^9
 =-(\zeta^3-\zeta^9) \sqrt{13},\\
  a_{13}&=-2+2 \zeta^3+2 \zeta^5-2 \zeta^7-2 \zeta^9+2 \zeta^{12}+\zeta^2-\zeta^{10}
 =-(\zeta^{10}-\zeta^2) \sqrt{13},\\
  a_{14}&=2+2 \zeta^4+2 \zeta^5-2 \zeta^7-2 \zeta^8-2 \zeta^{12}+\zeta^{11}-\zeta
 =-(\zeta^{11}-\zeta) \sqrt{13},\\
  a_{15}&=2+2 \zeta+2 \zeta^3-2 \zeta^9-2 \zeta^{11}-2 \zeta^{12}+\zeta^8-\zeta^4
 =-(\zeta^8-\zeta^4) \sqrt{13},\\
  a_{16}&=2 \zeta-2 \zeta^2+2 \zeta^4-2 \zeta^8+2 \zeta^{10}-2 \zeta^{11}+\zeta^{12}-1
 =-(\zeta^{12}-1) \sqrt{13}.
\endaligned$$
The other terms can be calculated in the same way. In conclusion,
we find that $(ST)^2=(ST)^{-1}$, i.e., $(ST)^3=1$. Hence, we have
$$S^2=T^{13}=(ST)^3=1.\eqno{(3.9)}$$
Let $u=ST$ and $v=S$. Then $uv=STS$. Hence,
$$u^3=v^2=(uv)^{13}=1.\eqno{(3.10)}$$
According to \cite{S}, put
$$(2, 3, n; p):=\langle u, v: u^3=v^2=(uv)^n=(u^{-1} v^{-1} uv)^p=1 \rangle.$$
Let $P=(uv)^{-1}$ and $Q=(uv)^2 u$. Then $u=P^2 Q$ and $v=P^3 Q$.
In \cite{S}, Sinkov proved the following:

{\bf Theorem 3.2.} (see \cite{S}). {\it Two operators of periods
$2$ and $3$ generate the simple group of order $1092$ if and only
if they satisfy one of the following sets of independent
relations:
$$\aligned
  A: \quad &(2, 3, 7; 6),\\
  B: \quad &(2, 3, 7; 7),\\
  C: \quad &(2, 3, 7); (Q^2 P^6)^3=1,\\
  D: \quad &(2, 3, 13); (Q^3 P^4)^3=1.
\endaligned$$}

  In our case, $P=S T^{-1} S$ and $Q=S T^3$.
We have
$$Q=-\frac{1}{\sqrt{13}} \left(\matrix
  \zeta^7-\zeta^9 & \zeta^4-\zeta^{10} & \zeta^2-\zeta^7 & \zeta^{10}-1 & \zeta^8-\zeta^4 & \zeta^8-\zeta^9\\
  \zeta^5-\zeta^{11} & \zeta^{11}-\zeta^3 & \zeta^{10}-\zeta^{12} & \zeta^7-\zeta^3 & \zeta^{12}-1 & \zeta^7-\zeta^{10}\\
  \zeta^{12}-\zeta^4 & \zeta^6-\zeta^8 & \zeta^8-\zeta & \zeta^{11}-\zeta^{12} & \zeta^{11}-\zeta & \zeta^4-1\\
  1-\zeta^3 & \zeta^9-\zeta^5 & \zeta^4-\zeta^5 & \zeta^6-\zeta^4 & \zeta^9-\zeta^3 & \zeta^{11}-\zeta^6\\
  \zeta^{10}-\zeta^6 & 1-\zeta & \zeta^3-\zeta^6 & \zeta^8-\zeta^2 & \zeta^2-\zeta^{10} & \zeta^3-\zeta\\
  \zeta-\zeta^2 & \zeta^{12}-\zeta^2 & 1-\zeta^9 & \zeta-\zeta^9 & \zeta^7-\zeta^5 & \zeta^5-\zeta^{12}
 \endmatrix\right).\eqno{(3.11)}$$
$$Q^2=-\frac{1}{\sqrt{13}} \left(\matrix
  \zeta^3-\zeta^8 & 1-\zeta^2 & \zeta^3-\zeta^9 & \zeta^{11}-\zeta^{12} & \zeta^8-\zeta^{11} & 1-\zeta^9\\
  \zeta-\zeta^3 & \zeta-\zeta^7 & 1-\zeta^5 & 1-\zeta^3 & \zeta^8-\zeta^4 & \zeta^7-\zeta^8\\
  1-\zeta^6 & \zeta^9-\zeta & \zeta^9-\zeta^{11} & \zeta^{11}-\zeta^7 & 1-\zeta & \zeta^7-\zeta^{10}\\
  \zeta-\zeta^2 & \zeta^2-\zeta^5 & \zeta^4-1 & \zeta^{10}-\zeta^5 & 1-\zeta^{11} & \zeta^{10}-\zeta^4\\
  \zeta^{10}-1 & \zeta^9-\zeta^5 & \zeta^5-\zeta^6 & \zeta^{12}-\zeta^{10} & \zeta^{12}-\zeta^6 & 1-\zeta^8\\
  \zeta^6-\zeta^2 & \zeta^{12}-1 & \zeta^3-\zeta^6 & 1-\zeta^7 & \zeta^4-\zeta^{12} & \zeta^4-\zeta^2
 \endmatrix\right).\eqno{(3.12)}$$
$$Q^3=-\frac{1}{\sqrt{13}} \left(\matrix
  \zeta^{11}-\zeta & \zeta^{12}-\zeta^8 & 1-\zeta & \zeta^6-\zeta^4 & \zeta^4-\zeta^{11} & 1-\zeta^8\\
  1-\zeta^9 & \zeta^8-\zeta^9 & \zeta^4-\zeta^7 & 1-\zeta^7 & \zeta^2-\zeta^{10} & \zeta^{10}-\zeta^8\\
  \zeta^{10}-\zeta^{11} & 1-\zeta^3 & \zeta^7-\zeta^3 & \zeta^{12}-\zeta^7 & 1-\zeta^{11} & \zeta^5-\zeta^{12}\\
  \zeta^9-\zeta^7 & \zeta^2-\zeta^9 & \zeta^5-1 & \zeta^2-\zeta^{12} & \zeta-\zeta^5 & 1-\zeta^{12}\\
  \zeta^6-1 & \zeta^3-\zeta^{11} & \zeta^5-\zeta^3 & 1-\zeta^4 & \zeta^5-\zeta^4 & \zeta^9-\zeta^6\\
  \zeta^6-\zeta & \zeta^2-1 & \zeta-\zeta^8 & \zeta^3-\zeta^2 & 1-\zeta^{10} & \zeta^6-\zeta^{10}
\endmatrix\right).\eqno{(3.13)}$$
$$Q^4=-\frac{1}{\sqrt{13}} \left(\matrix
  \zeta^{12}-\zeta^2 & \zeta^4-1 & \zeta^2-\zeta^3 & \zeta^6-\zeta^4 & 1-\zeta^7 & \zeta^{12}-\zeta^7\\
  \zeta^5-\zeta & \zeta^4-\zeta^5 & \zeta^{10}-1 & \zeta^4-\zeta^{11} & \zeta^2-\zeta^{10} & 1-\zeta^{11}\\
  \zeta^{12}-1 & \zeta^6-\zeta^9 & \zeta^{10}-\zeta^6 & 1-\zeta^8 & \zeta^{10}-\zeta^8 & \zeta^5-\zeta^{12}\\
  \zeta^9-\zeta^7 & \zeta^6-1 & \zeta^6-\zeta & \zeta-\zeta^{11} & \zeta^9-1 & \zeta^{11}-\zeta^{10}\\
  \zeta^2-\zeta^9 & \zeta^3-\zeta^{11} & \zeta^2-1 & \zeta^8-\zeta^{12} & \zeta^9-\zeta^8 & \zeta^3-1\\
  \zeta^5-1 & \zeta^5-\zeta^3 & \zeta-\zeta^8 & \zeta-1 & \zeta^7-\zeta^4 & \zeta^3-\zeta^7
\endmatrix\right).\eqno{(3.14)}$$
$$Q^5=-\frac{1}{\sqrt{13}} \left(\matrix
  \zeta^5-\zeta^{10} & \zeta^{10}-\zeta^{12} & \zeta^7-1 & \zeta^{11}-\zeta^{12} & 1-\zeta^3 & \zeta^{11}-\zeta^7\\
  \zeta^{11}-1 & \zeta^6-\zeta^{12} & \zeta^{12}-\zeta^4 & \zeta^8-\zeta^{11} & \zeta^8-\zeta^4 & 1-\zeta\\
  \zeta^4-\zeta^{10} & \zeta^8-1 & \zeta^2-\zeta^4 & 1-\zeta^9 & \zeta^7-\zeta^8 & \zeta^7-\zeta^{10}\\
  \zeta-\zeta^2 & \zeta^{10}-1 & \zeta^6-\zeta^2 & \zeta^8-\zeta^3 & \zeta^3-\zeta & \zeta^6-1\\
  \zeta^2-\zeta^5 & \zeta^9-\zeta^5 & \zeta^{12}-1 & \zeta^2-1 & \zeta^7-\zeta & \zeta-\zeta^9\\
  \zeta^4-1 & \zeta^5-\zeta^6 & \zeta^3-\zeta^6 & \zeta^9-\zeta^3 & \zeta^5-1 & \zeta^{11}-\zeta^9
\endmatrix\right).\eqno{(3.15)}$$
$$Q^6=-\frac{1}{\sqrt{13}} \left(\matrix
  \zeta^4-\zeta^6 & \zeta^2-\zeta^8 & \zeta^9-\zeta & \zeta^{10}-1 & \zeta^7-\zeta^3 & \zeta^{11}-\zeta^{12}\\
  \zeta^3-\zeta^9 & \zeta^{10}-\zeta^2 & \zeta^5-\zeta^7 & \zeta^8-\zeta^4 & \zeta^{12}-1 & \zeta^{11}-\zeta\\
  \zeta^6-\zeta^{11} & \zeta-\zeta^3 & \zeta^{12}-\zeta^5 & \zeta^8-\zeta^9 & \zeta^7-\zeta^{10} & \zeta^4-1\\
  1-\zeta^3 & \zeta^{10}-\zeta^6 & \zeta-\zeta^2 & \zeta^9-\zeta^7 & \zeta^{11}-\zeta^5 & \zeta^4-\zeta^{12}\\
  \zeta^9-\zeta^5 & 1-\zeta & \zeta^{12}-\zeta^2 & \zeta^{10}-\zeta^4 & \zeta^3-\zeta^{11} & \zeta^8-\zeta^6\\
  \zeta^4-\zeta^5 & \zeta^3-\zeta^6 & 1-\zeta^9 & \zeta^7-\zeta^2 & \zeta^{12}-\zeta^{10} & \zeta-\zeta^8
\endmatrix\right).\eqno{(3.16)}$$
$$Q^7=1.\eqno{(3.17)}$$
On the other hand,
$$P^4=-\frac{1}{\sqrt{13}} \left(\matrix
  \zeta^7-1 & \zeta^2-\zeta^7 & \zeta^6-\zeta^8 & \zeta^2-\zeta^{11} & \zeta^5-\zeta^6 & \zeta^8-\zeta^{11}\\
  \zeta^2-\zeta^7 & \zeta^{11}-1 & \zeta^5-\zeta^{11} & \zeta^7-\zeta^8 & \zeta^5-\zeta^8 & \zeta^6-\zeta^2\\
  \zeta^6-\zeta^8 & \zeta^5-\zeta^{11} & \zeta^8-1 & \zeta^2-\zeta^5 & \zeta^{11}-\zeta^7 & \zeta^6-\zeta^7\\
  \zeta^2-\zeta^{11} & \zeta^7-\zeta^8 & \zeta^2-\zeta^5 & \zeta^6-1 & \zeta^{11}-\zeta^6 & \zeta^7-\zeta^5\\
  \zeta^5-\zeta^6 & \zeta^5-\zeta^8 & \zeta^{11}-\zeta^7 & \zeta^{11}-\zeta^6 & \zeta^2-1 & \zeta^8-\zeta^2\\
  \zeta^8-\zeta^{11} & \zeta^6-\zeta^2 & \zeta^6-\zeta^7 & \zeta^7-\zeta^5 & \zeta^8-\zeta^2 & \zeta^5-1
\endmatrix\right).$$
We have
$$Q^3 P^4=-\frac{1}{\sqrt{13}} \left(\matrix
  \zeta^7-\zeta^5 & \zeta^2-\zeta^9 & \zeta^{10}-\zeta^5 & \zeta^6-\zeta^3 & \zeta^3-\zeta^7 & \zeta^{10}-\zeta^9\\
  \zeta^{12}-\zeta^6 & \zeta^{11}-\zeta^6 & \zeta^5-\zeta^3 & \zeta^{12}-\zeta^3 & \zeta^2-\zeta & \zeta-\zeta^{11}\\
  \zeta^6-\zeta & \zeta^4-\zeta^2 & \zeta^8-\zeta^2 & \zeta^9-\zeta^8 & \zeta^4-\zeta & \zeta^5-\zeta^9\\
  \zeta^{10}-\zeta^7 & \zeta^6-\zeta^{10} & \zeta^4-\zeta^3 & \zeta^6-\zeta^8 & \zeta^{11}-\zeta^4 & \zeta^3-\zeta^8\\
  \zeta^{10}-\zeta & \zeta^{12}-\zeta^{11} & \zeta^2-\zeta^{12} & \zeta-\zeta^7 & \zeta^2-\zeta^7 & \zeta^8-\zeta^{10}\\
  \zeta^5-\zeta^4 & \zeta^{12}-\zeta^9 & \zeta^4-\zeta^8 & \zeta^7-\zeta^{12} & \zeta^9-\zeta^{11} & \zeta^5-\zeta^{11}
\endmatrix\right),$$
$$(Q^3 P^4)^2=-\frac{1}{\sqrt{13}} \left(\matrix
  \zeta^8-\zeta^6 &  \zeta^7-\zeta & \zeta^{12}-\zeta^7 & \zeta^6-\zeta^3 & \zeta^{12}-\zeta^3 & \zeta^9-\zeta^8\\
  \zeta^4-\zeta^{11} & \zeta^7-\zeta^2 & \zeta^{11}-\zeta^9 & \zeta^3-\zeta^7 & \zeta^2-\zeta & \zeta^4-\zeta\\
  \zeta^8-\zeta^3 & \zeta^{10}-\zeta^8 & \zeta^{11}-\zeta^5 & \zeta^{10}-\zeta^9 & \zeta-\zeta^{11} & \zeta^5-\zeta^9\\
  \zeta^{10}-\zeta^7 & \zeta^{10}-\zeta & \zeta^5-\zeta^4 & \zeta^5-\zeta^7 & \zeta^6-\zeta^{12} & \zeta-\zeta^6\\
  \zeta^6-\zeta^{10} & \zeta^{12}-\zeta^{11} & \zeta^{12}-\zeta^9 & \zeta^9-\zeta^2 & \zeta^6-\zeta^{11} & \zeta^2-\zeta^4\\
  \zeta^4-\zeta^3 & \zeta^2-\zeta^{12} & \zeta^4-\zeta^8 & \zeta^5-\zeta^{10} & \zeta^3-\zeta^5 & \zeta^2-\zeta^8
\endmatrix\right),$$
and $(Q^3 P^4)^3=-I$. Note that in the projective coordinates,
this means that $(Q^3 P^4)^3=1$. Hence, we prove that the elements
$u$ and $v$ above satisfy the following relations: $(2, 3, 13)$
and $(Q^3 P^4)^3=1$, which is a presentation for the simple group
$PSL(2, 13)$ of order $1092$ by Theorem 3.2. Since the group is
simple and the generating matrices are non-trivial, we must have
$G=\langle u, v \rangle \cong PSL(2, 13)$. Hence, $\langle P, Q
\rangle=\langle S, T \rangle=G$. This completes the proof of
Theorem 3.1.

\flushpar $\qquad \qquad \qquad \qquad \qquad \qquad \qquad \qquad
\qquad \qquad \qquad \qquad \qquad \qquad \qquad \qquad \qquad
\qquad \quad \boxed{}$

{\bf Theorem 3.3.} {\it Let $x_3=QP^2$ and $y_3=Q^5 P^2$. Then
$\langle x_3, y_3 \rangle=G$.}

{\it Proof}. We have
$$P^2=-\frac{1}{\sqrt{13}} \left(\matrix
  1-\zeta & \zeta-\zeta^4 & \zeta^3-\zeta^{12} & \zeta^9-\zeta^4 & \zeta^{12}-\zeta^{10} & \zeta^9-\zeta^3\\
  \zeta-\zeta^4 & 1-\zeta^9 & \zeta^9-\zeta^{10} & \zeta^3-\zeta & \zeta^3-\zeta^{10} & \zeta^4-\zeta^{12}\\
  \zeta^3-\zeta^{12} & \zeta^9-\zeta^{10} & 1-\zeta^3 & \zeta^{10}-\zeta^4 & \zeta-\zeta^9 & \zeta-\zeta^{12}\\
  \zeta^9-\zeta^4 & \zeta^3-\zeta & \zeta^{10}-\zeta^4 & 1-\zeta^{12} & \zeta^{12}-\zeta^9 & \zeta^{10}-\zeta\\
  \zeta^{12}-\zeta^{10} & \zeta^3-\zeta^{10} & \zeta-\zeta^9 & \zeta^{12}-\zeta^9 & 1-\zeta^4 & \zeta^4-\zeta^3\\
  \zeta^9-\zeta^3 & \zeta^4-\zeta^{12} & \zeta-\zeta^{12} & \zeta^{10}-\zeta & \zeta^4-\zeta^3 & 1-\zeta^{10}
\endmatrix\right).$$
$$Q P^2=-\frac{1}{\sqrt{13}} \left(\matrix
  \zeta^{12}-\zeta^3 & \zeta^6-\zeta^5 & \zeta^2-\zeta^{12} & \zeta^5-\zeta^{11} & \zeta-\zeta^6 & \zeta-\zeta^3\\
  \zeta^5-\zeta^4 & \zeta^4-\zeta & \zeta^2-\zeta^6 & \zeta^9-\zeta & \zeta^6-\zeta^8 & \zeta^9-\zeta^2\\
  \zeta^5-\zeta^2 & \zeta^6-\zeta^{10} & \zeta^{10}-\zeta^9 & \zeta^3-\zeta^5 & \zeta^3-\zeta^9 & \zeta^2-\zeta^7\\
  \zeta^2-\zeta^8 & \zeta^7-\zeta^{12} & \zeta^{10}-\zeta^{12} & \zeta-\zeta^{10} & \zeta^7-\zeta^8 & \zeta^{11}-\zeta\\
  \zeta^{12}-\zeta^4 & \zeta^5-\zeta^7 & \zeta^{11}-\zeta^4 & \zeta^8-\zeta^9 & \zeta^9-\zeta^{12} & \zeta^{11}-\zeta^7\\
  \zeta^8-\zeta^{10} & \zeta^4-\zeta^{10} & \zeta^6-\zeta^{11} & \zeta^8-\zeta^{11} & \zeta^7-\zeta^3 & \zeta^3-\zeta^4
\endmatrix\right).\eqno{(3.18)}$$
$$(QP^2)^2=-\frac{1}{\sqrt{13}} \left(\matrix
  \zeta-\zeta^{10} & \zeta^8-\zeta^9 & \zeta^8-\zeta^{11} & \zeta^{11}-\zeta^5 & \zeta-\zeta^9 & \zeta^5-\zeta^3\\
  \zeta^7-\zeta^8 & \zeta^9-\zeta^{12} & \zeta^7-\zeta^3 & \zeta^6-\zeta & \zeta^8-\zeta^6 & \zeta^9-\zeta^3\\
  \zeta^{11}-\zeta & \zeta^{11}-\zeta^7 & \zeta^3-\zeta^4 & \zeta^3-\zeta & \zeta^2-\zeta^9 & \zeta^7-\zeta^2\\
  \zeta^8-\zeta^2 & \zeta^4-\zeta^{12} & \zeta^{10}-\zeta^8 & \zeta^{12}-\zeta^3 & \zeta^5-\zeta^4 & \zeta^5-\zeta^2\\
  \zeta^{12}-\zeta^7 & \zeta^7-\zeta^5 & \zeta^{10}-\zeta^4 & \zeta^6-\zeta^5 & \zeta^4-\zeta & \zeta^6-\zeta^{10}\\
  \zeta^{12}-\zeta^{10} & \zeta^4-\zeta^{11} & \zeta^{11}-\zeta^6 & \zeta^2-\zeta^{12} & \zeta^2-\zeta^6 & \zeta^{10}-\zeta^9
\endmatrix\right).$$
$$(QP^2)^3=I.\eqno{(3.19)}$$
$$Q^5 P^2=-\frac{1}{\sqrt{13}} \left(\matrix
  \zeta^8-\zeta^5 & \zeta^4-\zeta^8 & \zeta^2-\zeta & \zeta^4-\zeta^6 & \zeta^9-\zeta^2 & \zeta-\zeta^6\\
  \zeta^5-\zeta^9 & \zeta^7-\zeta^6 & \zeta^{10}-\zeta^7 & \zeta^9-\zeta^2 & \zeta^{10}-\zeta^2 & \zeta^3-\zeta^5\\
  \zeta^{12}-\zeta^{11} & \zeta^6-\zeta^3 & \zeta^{11}-\zeta^2 & \zeta-\zeta^6 & \zeta^3-\zeta^5 & \zeta^{12}-\zeta^5\\
  \zeta^7-\zeta^9 & \zeta^{11}-\zeta^4 & \zeta^7-\zeta^{12} & \zeta^5-\zeta^8 & \zeta^9-\zeta^5 & \zeta^{11}-\zeta^{12}\\
  \zeta^{11}-\zeta^4 & \zeta^{11}-\zeta^3 & \zeta^8-\zeta^{10} & \zeta^8-\zeta^4 & \zeta^6-\zeta^7 & \zeta^3-\zeta^6\\
  \zeta^7-\zeta^{12} & \zeta^8-\zeta^{10} & \zeta^8-\zeta & \zeta-\zeta^2 & \zeta^7-\zeta^{10} & \zeta^2-\zeta^{11}
\endmatrix\right).\eqno{(3.20)}$$
$$(Q^5 P^2)^2=-I.\eqno{(3.21)}$$
Note that in the projective coordinates, this means that $(Q^5
P^2)^2=1$.

  Let $\widetilde{P}=Q^4$ and $\widetilde{Q}=P^2$. By $P^{13}=Q^7=1$,
we have
$$\widetilde{P}^7=Q^{28}=1, \quad \widetilde{Q}^{13}=P^{26}=1.$$
Put $x_3=\widetilde{P}^2 \widetilde{Q}$ and $y_3=\widetilde{P}^3
\widetilde{Q}$. Then
$$x_3=QP^2, \quad y_3=Q^5 P^2.$$
By $(QP^2)^3=1$ and $(Q^5 P^2)^2=1$, we have $x_3^3=y_3^2=1$.
Moreover,
$$x_3 y_3=QP^2 \cdot Q^5 P^2=Q^3 \cdot (Q^5 P^2)^2=Q^3.$$
Hence, $(x_3 y_3)^7=1$. On the other hand,
$$(\widetilde{Q}^2 \widetilde{P}^6)^3=(P^4 Q^3)^3.$$
Note that $(Q^3 P^4)^3=1$. This implies that $(P^4 Q^3)^3=1$.
Thus, $x_3$ and $y_3$ satisfy the following relation: $(2, 3, 7)$;
$(\widetilde{Q}^2 \widetilde{P}^6)^3=1$. By Theorem 3.2, we have
$\langle x_3, y_3 \rangle=G$. This completes the proof of Theorem
3.3.

\flushpar $\qquad \qquad \qquad \qquad \qquad \qquad \qquad \qquad
\qquad \qquad \qquad \qquad \qquad \qquad \qquad \qquad \qquad
\qquad \quad \boxed{}$

  Let us present some details about $G$ and its character table. For
$G$, there are irreducible characters of degrees
$$1, \quad \text{$7$ (twice)}, \quad \text{$12$ ($3$ times)}, \quad
  13, \quad \text{$14$ (twice)},$$
where the sum of squares of the degrees is the group order:
$$1^2+2 \cdot 7^2+3 \cdot 12^2+13^2+2 \cdot 14^2=1092.$$
There are also irreducible characters of degrees
$$\text{$6$ (twice)}, \quad \text{$12$ ($3$ times)}, \quad \text{$14$ ($3$ times)},$$
where the sum of squares of the degrees is the group order:
$$2 \cdot 6^2+3 \cdot 12^2+3 \cdot 14^2=1092.$$
The group $G$ has order $1092$, it has standard generators $a$ and
$b$ of orders $2$ and $3$, respectively, such that $ab$ has order
$13$. In table 1, we reproduce from \cite{CC} some of the
character table of $G$. In terms of the standard generators,
representatives of the conjugacy classes are listed in Table 2
(see \cite{MM}).

$$\text{Table $1$. Some of the character table of $PSL(2, 13)$}$$
$$\matrix
            & 1A & 2A & 3A & 6A & 7A & 7B & 7C & 13A & 13B\\
     \chi_1 &  1 &  1 &  1 &  1 &  1 &  1 &  1 &   1 &   1\\
     \chi_2 &  7 & -1 &  1 & -1 &  0 &  0 &  0 & \frac{1-\sqrt{13}}{2} & \frac{1+\sqrt{13}}{2}\\
     \chi_3 &  7 & -1 &  1 & -1 &  0 &  0 &  0 & \frac{1+\sqrt{13}}{2} & \frac{1-\sqrt{13}}{2}\\
  \chi_{10} &  6 &  0 &  0 &  0 & -1 & -1 & -1 & \frac{-1+\sqrt{13}}{2} & \frac{-1-\sqrt{13}}{2}\\
  \chi_{11} &  6 &  0 &  0 &  0 & -1 & -1 & -1 & \frac{-1-\sqrt{13}}{2} & \frac{-1+\sqrt{13}}{2}\\
  \chi_{15} & 14 &  0 &  2 &  0 &  0 &  0 &  0 &   1 &   1
\endmatrix$$

$$\text{Table $2$. Representatives of the conjugacy classes of $PSL(2, 13)$}$$
$$\text{and the order of the normalizer of a representative}$$

$$\matrix
  \text{Conjugacy class} & \text{Representative} & |N_G(\langle h \rangle)|\\
  1A & \text{Identity} & 1092\\
  2A & (abababb)^3 & 12\\
  3A & (abababb)^2 & 12\\
  6A & abababb     & 12\\
  7A & ababb       & 14\\
  7B & (ababb)^2   & 14\\
  7C & (ababb)^4   & 14\\
 13A & ab          & 78\\
 13B & abab        & 78
\endmatrix$$

  We have
$$\text{Tr}(S)=0, \quad \text{Tr}(T)=\frac{-1-\sqrt{13}}{2}, \quad \text{Tr}(ST)=0.\eqno{(3.22)}$$
Hence, the above six-dimensional representation corresponds to the
character $\chi_{11}$ in Table $1$. In terms of our
six-dimensional representations of $G$, $a=S$, $b=ST$, $ab=T$.
$$abab^2=T^2 ST, \quad (abab^2)^2=(T^2 ST)^2, \quad (abab^2)^4=(T^2 ST)^4.$$
Note that $Q=ST^3$, we have
$$abab^2=T^2 Q T^{-2}, \quad (abab^2)^2=T^2 Q^2 T^{-2}, \quad (abab^2)^4=T^2 Q^4 T^{-2}.$$
We will give the other two six-dimensional representations of $G$:
$\langle x_1, y_1 \rangle$ and $\langle x_2, y_2 \rangle$, such
that $z_1=x_1^{-1} y_1^{-1}$, $z_2=x_2^{-1} y_2^{-1}$ and
$z_3=x_3^{-1} y_3^{-1}$ correspond to the conjugacy classes $7A$,
$7B$ and $7C$, respectively.

  It is well-known that the modular group $\Gamma=PSL(2, {\Bbb Z})$
is generated by the following linear fractional transformations:
$$T \tau=\tau+1, \quad S \tau=-\frac{1}{\tau}.$$
Let $P=ST$. Then
$$P \tau=-\frac{1}{\tau+1}.$$
Here, $S^2=P^3=1$. It is a discontinuous group acting on the upper
half-plane ${\Bbb H}=\{ \text{$\tau \in {\Bbb C}$: Im $\tau>0$}
\}$, and has a fundamental domain $F$ given by $|\tau|>1$,
$-\frac{1}{2} \leq \text{Re $\tau<\frac{1}{2}$}$ and $|\tau|=1$,
$-\frac{1}{2} \leq \text{Re $\tau \leq 0$}$. The only fixed points
in $F$ of elliptic transformations of $\Gamma$ are $\tau=i$ (fixed
by $S$) and $\tau=\rho=e^{2 \pi i/3}$ (fixed by $P$ and $P^{-1}$);
the local uniformizing variables are $(\tau-i)^2$ and
$(\tau-\rho)^3$ respectively. In addition $\infty$ is fixed by the
parabolic transformation $T$. The vertical sides of the boundary
of $F$ are mapped into each other by $T$ and $T^{-1}$, and the
curved side into itself by $S$.

  Now, we give some basic fact about subgroups of the modular group
and permutations (see \cite{ASD)}. Suppose that $G$ is a subgroup
of $\Gamma$ of finite index $\mu$. Then $G$ has a connected
fundamental domain $D$ consisting of $\mu$ copies of $F$, the
transforms of $F$ by a set of coset representatives for $G$ in
$\Gamma$. If the elements of $G$ conjugate in $\Gamma$ to $S$ and
$P$ form respectively $e_2$ and $e_3$ conjugacy classes in $G$,
then the boundary of $F$ will have $e_2$ and $e_3$ inequivalent
fixed point vertices of orders $2$ and $3$ respectively. Suppose
further that every element of $G$ conjugate in $\Gamma$ to a
nonzero power of $T$ is conjugate in $G$ to some power of one of
$$T^{\mu_1}, g_2 T^{\mu_2} g_2^{-1}, \cdots, g_h T^{\mu_h} g_h^{-1}, \quad \text{(in $G$)}$$
where $g_1=I$, $g_2, \cdots, g_h$ are in $\Gamma$, and no $g_i
g_j^{-1}$ is in $G$. Then the boundary of $F$ will have $h$
inequivalent parabolic fixed point vertices (called cusps) at
$\infty$ and at $g_2 \infty, \cdots, g_h \infty$ which are
rational points on the real axis. We then have
$$\mu=\mu_1+\mu_2+\cdots+\mu_h$$
and
$$g=1+\frac{\mu}{12}-\frac{h}{2}-\frac{e_2}{4}-\frac{e_3}{3},$$
where $g$ is the genus of the Riemann surface ${\Bbb H}/G$.
Following Wohlfahrt \cite{W1}, we define the level $l$ of $G$ to
be the least common multiple of $\mu_1, \mu_2, \cdots, \mu_h$. We
write $G^N$ for the intersection of the conjugates of $G$ in
$\Gamma$, so that $G^N$ is the maximal normal subgroup of $\Gamma$
contained in $G$.

  We consider permutations on $\mu$ letters named as the integers
$1$ to $\mu$, where $1$ is specially distinguished. We say that a
pair $(s, p)$ of permutations is legitimate if $s^2=p^3=I$ and the
group $\Sigma$ generated by $s$ and $p$ is transitive. If $\sigma$
is any element of $T_{\mu}$, we write $(s, p) \sim (\sigma s
\sigma^{-1}, \sigma p \sigma^{-1})$, and if $\sigma$ is any
element of $T_{\mu}$ fixing $1$, we write $(s, p) \sim_1 (\sigma s
\sigma^{-1}, \sigma p \sigma^{-1})$. Then we have

{\bf Theorem 3.4.} (see \cite{ASD}). {\it There is a one-to-one
correspondence between subgroups of index $\mu$ in the modular
group and equivalence classes of legitimate pairs of permutations
$(s, p)$ under the equivalence relation $\sim_1$. If $G$ is a
subgroup and $(s, p)$ a representative of the corresponding
equivalence class, then

\roster

\item $e_2$ and $e_3$ are the number of letters fixed by $s$ and
$p$ respectively,

\item $t=sp$ has $h$ cycles of lengths $\mu_1, \cdots, \mu_h$, and
$\mu_1$ is the length of the cycle containing $1$,

\item $\Sigma$ is isomorphic to the factor group $\Gamma/G^N$,

\item $G$ is maximal if and only if $\Gamma$ is primitive.
\endroster}

  In fact, the simple group $PSL(2, 13)$ is a primitive group
of degree $14$ (see \cite{Mi}).

{\bf Theorem 3.5.} {\it Let $x_1=Q^6 \cdot PQ^2 P^{10}$, $y_1=P
Q^2 P^{10}$, $x_2=Q^5 \cdot Q^5 P^2 \cdot P^2 Q^6 P^8 \cdot Q^5
P^2$ and $y_2=Q^5 P^2 \cdot P^2 Q^6 P^8 \cdot Q^5 P^2$. Then
$$G=\langle x_1, y_1 \rangle=\langle x_2, y_2 \rangle.$$}

{\it Proof}. In \cite{S}, Sinkov studied the group $PSL(2, 13)$ of
order $1092$ as a permutation group of degree $14$ (see \cite{S}),
which we denote by $G_{1092}$. The $91$ elements of order two
contained in $G_{1092}$ are all conjugate and it is sufficient to
consider only one of them. We choose it to be
$$s=(1, 12)(2, 11)(3, 10)(4, 9)(5, 8)(6, 7).$$
The largest subgroup of $G_{1092}$ within which $s$ is invariant
is the dihedral group of order $12$ generated by
$$(1, 4, 3, 12, 9, 10)(2, 8, 6, 11, 5, 7)$$
and
$$(1, 12)(2, 6)(3, 4)(7, 11)(9, 10)(13, 14).$$
Under this subgroup the $180$ remaining elements of order three
are divided up into $16$ conjugate sets. Of these, two contain
only six distinct elements each; the remaining $14$ sets each
contain $12$ elements. In his paper \cite{S}, Sinkov gave below
one member of each of these sets together with the order of its
product with $s$ which we denote by $\text{ord}(ps)$.
$$\aligned
  p_1&=(1, 13, 10)(2, 3, 6)(4, 9, 11)(5, 12, 7), \quad \text{ord}(p_1 s)=6,\\
  p_2&=p_1^2, \quad \text{ord}(p_2 s)=6,\\
  p_3&=(2, 10, 4)(11, 13, 5)(3, 6, 7)(8, 12, 9), \quad \text{ord}(p_3 s)=6,\\
  p_4&=p_3^2, \quad \text{ord}(p_4 s)=6.
\endaligned$$
$$\aligned
  p_5&=(3, 10, 12)(4, 6, 13)(5, 11, 8)(14, 9, 7), \quad \text{ord}(p_5 s)=3,\\
  p_6&=p_5^2, \quad \text{ord}(p_6 s)=3.
\endaligned$$
$$\aligned
  p_7&=(2, 11, 14)(3, 4, 8)(5, 9, 10)(6, 13, 7), \quad \text{ord}(p_7 s)=2,\\
  p_8&=(1, 13, 12)(9, 4, 14)(3, 8, 6)(5, 10, 7), \quad \text{ord}(p_8 s)=2,\\
  p_9&=(2, 8, 9)(4, 14, 13)(5, 10, 6)(7, 12, 11), \quad \text{ord}(p_9 s)=7,\\
  p_{10}&=(2, 7, 8)(3, 10, 11)(5, 13, 9)(6, 12, 14), \quad \text{ord}(p_{10} s)=13.
\endaligned$$
$$\aligned
  p_{11}&=(2, 3, 4)(6, 9, 11)(7, 12, 14)(8, 10, 13), \quad \text{ord}(p_{11} s)=7,\\
  p_{12}&=p_{11}^2, \quad \text{ord}(p_{12} s)=7,\\
  p_{13}&=(2, 14, 5)(3, 9, 13)(4, 7, 11)(8, 10, 12), \quad \text{ord}(p_{13} s)=7,\\
  p_{14}&=p_{13}^2, \quad \text{ord}(p_{14} s)=7.
\endaligned$$
$$\aligned
  p_{15}&=(1, 10, 6)(3, 8, 9)(4, 11, 12)(7, 13, 14), \quad \text{ord}(p_{15} s)=7,\\
  p_{16}&=(1, 10, 4)(3, 6, 14)(5, 12, 8)(9, 13, 11), \quad \text{ord}(p_{16} s)=13.
\endaligned$$

  Since no group satisfying the relations $(2, 3, 6)$ is simple, it is obvious that
none of the first eight of the above elements when coupled with
$s$ will generate the entire group. $p_{12}$ and $p_{14}$ satisfy
with $s$ the same defining relations as do $p_{11}$ and $p_{13}$,
respectively.

  $p_9$ is transformed into $p_{15}$ by the substitution
$$(1, 2)(3, 5)(4, 7)(6, 9)(8, 10)(11, 12)(13, 14)$$
which is commutative with $s$. Similarly $p_{10}$ is transformed
into $p_{16}$ by
$$(1, 7, 10, 5, 9, 11, 12, 6, 3, 8, 4, 2)$$
which is also commutative with $s$. Hence there remain only $p_9$,
$p_{10}$, $p_{11}$ and $p_{13}$ to be considered. The orders of
the commutators of $p_9$, $p_{11}$ and $p_{13}$ with $s$ are
respectively $13$, $6$, $7$.

  Let $[a, b]:=a^{-1} b^{-1} ab$. We have
$$p_{11} s=(1, 12, 14, 6, 4, 11, 7)(2, 10, 13, 5, 8, 3, 9),$$
$$[p_{11}, s]=(1, 14, 12, 4, 13, 9)(3, 7, 6, 10, 8, 5).$$
Hence,
$$\langle p_{11}, s \rangle \cong (2, 3, 7; 6).$$
$$p_{13} s=(1, 12, 5, 11, 9, 13, 10)(2, 14, 8, 3, 4, 6, 7),$$
$$[p_{13}, s]=(1, 5, 8, 12, 4, 14, 9)(2, 3, 10, 11, 7, 13, 6).$$
Hence,
$$\langle p_{13}, s \rangle \cong (2, 3, 7; 7).$$
$$p_9 s=(1, 12, 2, 5, 3, 10, 7)(4, 14, 13, 9, 11, 6, 8),$$
$$[p_9, s]=(1, 2, 14, 11, 12, 8, 6, 10, 4, 9, 3, 7, 5).$$
Hence,
$$\langle p_9, s \rangle \cong (2, 3, 7; 13).$$
$$p_{10} s=(1, 12, 14, 7, 5, 13, 4, 9, 8, 11, 10, 2, 6),$$
$$[p_{10}, s]=(1, 14, 12, 5, 9, 4, 8)(2, 13, 11, 3, 6, 7, 10).$$
Hence,
$$\langle p_{10}, s \rangle \cong (2, 3, 13; 7).$$
Let
$$\psi: p_9 \mapsto QP^2, \quad s \mapsto Q^5 P^2.$$
Then $\psi(p_9 s)=Q^3$. The map $\psi$ is an isomorphism between
the groups $\langle p_9, s \rangle$ and $\langle QP^2, Q^5 P^2
\rangle=G$. From now on, we assume that $p_9=QP^2$ and $s=Q^5
P^2$. Hence,
$$Q=(p_9 s)^5=(1, 10, 5, 12, 7, 3, 2)(4, 6, 9, 14, 8, 11, 13),$$
$$P=[(p_9 s)^{-5} p_9]^7=(1, 6, 8, 14, 13, 2, 7, 3, 11, 12, 9, 10, 5).$$

  We have
$$y_1:=P Q^2 P^{10}=(2, 8)(3, 4)(5, 12)(7, 9)(10, 14)(11, 13),\eqno{(3.23)}$$
which is of order two, and
$$x_1:=Q^6 \cdot PQ^2 P^{10}=(1, 8, 10)(2, 4, 11)(3, 9, 6)(5, 14, 7),\eqno{(3.24)}$$
which is of order three. Then $x_1 y_1=Q^6$ with order $7$, and
$$[x_1, y_1]=(1, 9, 13, 8, 2, 11, 7)(3, 4, 14, 5, 6, 12, 10),\eqno{(3.25)}$$
which is of order $7$. On the other hand,
$$P^2 Q^6 P^8=(1, 10)(2, 6)(4, 5)(7, 9)(8, 14)(11, 12).$$
$$y_2:=Q^5 P^2 \cdot P^2 Q^6 P^8 \cdot Q^5 P^2=(1, 2)(3, 12)(4, 6)(5, 14)(7, 11)(8, 9),\eqno{(3.26)}$$
which is of order two, and
$$x_2:=Q^5 \cdot Q^5 P^2 \cdot P^2 Q^6 P^8 \cdot Q^5 P^2=(1, 12, 10)(2, 11, 5)(4, 7, 14)(6, 13, 9),\eqno{(3.27)}$$
which is of order three. Then $x_2 y_2=Q^5$ with order $7$, and
$$[x_2, y_2]=(1, 2, 6, 9, 8, 4)(3, 10, 12, 7, 13, 11),\eqno{(3.28)}$$
which is of order $6$. By Theorem 3.2, we have
$$\langle x_1, y_1 \rangle=\langle x_2, y_2 \rangle=G.$$

  In the form of matrices defined over ${\Bbb Q}(\zeta)$, we have
$$P Q^2 P^{10}=\left(\matrix
  0 & 0 & 0 & -\zeta & 0 & 0\\
  0 & 0 & 0 & 0 & -\zeta^9 & 0\\
  0 & 0 & 0 & 0 & 0 & -\zeta^3\\
  \zeta^{12} & 0 & 0 & 0 & 0 & 0\\
  0 & \zeta^4 & 0 & 0 & 0 & 0\\
  0 & 0 & \zeta^{10} & 0 & 0 & 0
\endmatrix\right).\eqno{(3.29)}$$
$$\aligned
 &Q^6 \cdot PQ^2 P^{10}\\
=&-\frac{1}{\sqrt{13}} \left(\matrix
  \zeta^9-\zeta^{12} & \zeta^{11}-\zeta^7 & \zeta^8-\zeta^9 & \zeta^7-\zeta^5 & \zeta^4-\zeta^{11} & \zeta^4-\zeta^{12}\\
  \zeta^7-\zeta^3 & \zeta^3-\zeta^4 & \zeta^8-\zeta^{11} & \zeta^{10}-\zeta^4 & \zeta^{11}-\zeta^6 & \zeta^{10}-\zeta^8\\
  \zeta^7-\zeta^8 & \zeta^{11}-\zeta & \zeta-\zeta^{10} & \zeta^{12}-\zeta^7 & \zeta^{12}-\zeta^{10} & \zeta^8-\zeta^2\\
  \zeta^8-\zeta^6 & \zeta^2-\zeta^9 & \zeta-\zeta^9 & \zeta^4-\zeta & \zeta^2-\zeta^6 & \zeta^5-\zeta^4\\
  \zeta^9-\zeta^3 & \zeta^7-\zeta^2 & \zeta^5-\zeta^3 & \zeta^6-\zeta^{10} & \zeta^{10}-\zeta^9 & \zeta^5-\zeta^2\\
  \zeta^6-\zeta & \zeta^3-\zeta & \zeta^{11}-\zeta^5 & \zeta^6-\zeta^5 & \zeta^2-\zeta^{12} & \zeta^{12}-\zeta^3
\endmatrix\right).\endaligned\eqno{(3.30)}$$
$$P^2 Q^6 P^8=-\frac{1}{\sqrt{13}} \left(\matrix
  \zeta^8-\zeta^5 & \zeta^{12}-\zeta^3 & \zeta^4-\zeta^3 & \zeta^2-\zeta^4 & \zeta^{12}-\zeta^5 & \zeta^{10}-\zeta^2\\
  \zeta^{10}-\zeta & \zeta^7-\zeta^6 & \zeta^4-\zeta & \zeta^{12}-\zeta^5 & \zeta^5-\zeta^{10} & \zeta^4-\zeta^6\\
  \zeta^{10}-\zeta^9 & \zeta^{12}-\zeta^9 & \zeta^{11}-\zeta^2 & \zeta^{10}-\zeta^2 & \zeta^4-\zeta^6 & \zeta^6-\zeta^{12}\\
  \zeta^9-\zeta^{11} & \zeta^8-\zeta & \zeta^{11}-\zeta^3 & \zeta^5-\zeta^8 & \zeta-\zeta^{10} & \zeta^9-\zeta^{10}\\
  \zeta^8-\zeta & \zeta^3-\zeta^8 & \zeta^7-\zeta^9 & \zeta^3-\zeta^{12} & \zeta^6-\zeta^7 & \zeta^9-\zeta^{12}\\
  \zeta^{11}-\zeta^3 & \zeta^7-\zeta^9 & \zeta-\zeta^7 & \zeta^3-\zeta^4 & \zeta-\zeta^4 & \zeta^2-\zeta^{11}
\endmatrix\right).$$
$$\aligned
 &Q^5 P^2 \cdot P^2 Q^6 P^8 \cdot Q^5 P^2\\
=&-\frac{1}{\sqrt{13}} \left(\matrix
  \zeta^7-\zeta^6 & \zeta^8-\zeta^5 & \zeta^{11}-\zeta^2 & \zeta^4-\zeta^9 & \zeta^{12}-\zeta & \zeta^{10}-\zeta^3\\
  \zeta^8-\zeta^5 & \zeta^{11}-\zeta^2 & \zeta^7-\zeta^6 & \zeta^{12}-\zeta & \zeta^{10}-\zeta^3 & \zeta^4-\zeta^9\\
  \zeta^{11}-\zeta^2 & \zeta^7-\zeta^6 & \zeta^8-\zeta^5 & \zeta^{10}-\zeta^3 & \zeta^4-\zeta^9 & \zeta^{12}-\zeta\\
  \zeta^4-\zeta^9 & \zeta^{12}-\zeta & \zeta^{10}-\zeta^3 & \zeta^6-\zeta^7 & \zeta^5-\zeta^8 & \zeta^2-\zeta^{11}\\
  \zeta^{12}-\zeta & \zeta^{10}-\zeta^3 & \zeta^4-\zeta^9 & \zeta^5-\zeta^8 & \zeta^2-\zeta^{11} & \zeta^6-\zeta^7\\
  \zeta^{10}-\zeta^3 & \zeta^4-\zeta^9 & \zeta^{12}-\zeta & \zeta^2-\zeta^{11} & \zeta^6-\zeta^7 & \zeta^5-\zeta^8
\endmatrix\right).\endaligned\eqno{(3.31)}$$
$$\aligned
 &Q^5 \cdot Q^5 P^2 \cdot P^2 Q^6 P^8 \cdot Q^5 P^2\\
=&-\frac{1}{\sqrt{13}} \left(\matrix
  \zeta^9-\zeta^{10} & \zeta^5-\zeta^8 & \zeta-\zeta^{10} & \zeta^3-\zeta^{11} & \zeta^{11}-\zeta^9 & \zeta-\zeta^8\\
  \zeta^9-\zeta^{12} & \zeta^3-\zeta^{12} & \zeta^6-\zeta^7 & \zeta^9-\zeta^7 & \zeta-\zeta^8 & \zeta^8-\zeta^3\\
  \zeta^2-\zeta^{11} & \zeta^3-\zeta^4 & \zeta-\zeta^4 & \zeta^7-\zeta & \zeta^3-\zeta^{11} & \zeta^9-\zeta^7\\
  \zeta^2-\zeta^{10} & \zeta^4-\zeta^2 & \zeta^5-\zeta^{12} & \zeta^4-\zeta^3 & \zeta^8-\zeta^5 & \zeta^{12}-\zeta^3\\
  \zeta^6-\zeta^4 & \zeta^5-\zeta^{12} & \zeta^{10}-\zeta^5 & \zeta^4-\zeta & \zeta^{10}-\zeta & \zeta^7-\zeta^6\\
  \zeta^{12}-\zeta^6 & \zeta^2-\zeta^{10} & \zeta^6-\zeta^4 & \zeta^{11}-\zeta^2 & \zeta^{10}-\zeta^9 & \zeta^{12}-\zeta^9
\endmatrix\right).\endaligned\eqno{(3.32)}$$
This completes the proof of Theorem 3.5.

\flushpar $\qquad \qquad \qquad \qquad \qquad \qquad \qquad \qquad
\qquad \qquad \qquad \qquad \qquad \qquad \qquad \qquad \qquad
\qquad \quad \boxed{}$

  By the above argument, we find that as a permutation group of
degree $14$, $PSL(2, 13)$ has four presentations: $(2, 3, 7; 6)$,
$(2, 3, 7; 7)$, $(2, 3, 7; 13)$ and $(2, 3, 13; 7)$. In each case,
$t=sp$ has $h=2$ cycles of lengths $\mu_1$ and $\mu_2$. For the
representations of Sinkov: $s$, $p_9$, $p_{10}$, $p_{11}$ and
$p_{13}$, we have $e_2=2$ and $e_3=2$. On the other hand, for our
three representations: $\langle x_i, y_i \rangle$ ($i=1, 2, 3$),
we also have $e_2=e_3=2$. There are four cases:

(1) $PSL(2, 13) \cong (2, 3, 13; 7)$. Here, $\mu=14$, $h=2$,
$e_2=e_3=2$, $\mu_1=13$, $\mu_2=1$. Hence, $l=13$.

(2) $PSL(2, 13) \cong (2, 3, 7; 6)$. Here, $\mu=14$, $h=2$,
$e_2=e_3=2$, $\mu_1=7$, $\mu_2=7$. Hence, $l=7$.

(3) $PSL(2, 13) \cong (2, 3, 7; 7)$. Here, $\mu=14$, $h=2$,
$e_2=e_3=2$, $\mu_1=7$, $\mu_2=7$. Hence, $l=7$.

(4) $PSL(2, 13) \cong (2, 3, 7; 13)$. Here, $\mu=14$, $h=2$,
$e_2=e_3=2$, $\mu_1=7$, $\mu_2=7$. Hence, $l=7$.

  In the notation of \cite{ASD}, the first case corresponds to the
congruence subgroup $\Gamma_{13, 1}$ which is just $\Gamma_0(13)$,
the last three cases correspond to the noncongruence subgroup
$\Gamma_{7, 7}$. In all of these cases, we have
$$g=1+\frac{14}{12}-\frac{2}{2}-\frac{2}{4}-\frac{2}{3}=0,\eqno{(3.33)}$$
i.e., $g({\Bbb H}/\Gamma_{13, 1})=g({\Bbb H}/\Gamma_{7, 7})=0$. By
the above theorem, we have
$$\Gamma/\Gamma(13) \cong \Gamma/\Gamma_{7, 7}^N \cong PSL(2, 13),\eqno{(3.34)}$$
where $\Gamma(13)=\Gamma_0(13)^N$. Moreover, $\Gamma_0(13)$ and
$\Gamma_{7, 7}$ are maximal.

  In \cite{N}, Newman studied the maximal normal subgroups $G$ of the
modular group $\Gamma$; i.e. those normal subgroups $G$ such that
$\Gamma/G$ is simple. The principal congruence subgroups
$\Gamma(p)$ of prime level $p>3$ are such groups, since
$$\Gamma/\Gamma(p) \cong PSL(2, p).$$
However, these are not the only groups with quotient groups
isomorphic to $PSL(2, p)$. Newman showed that for a given $p$
there are in general many normal subgroups $G$ of different levels
such that
$$\Gamma/G \cong PSL(2, p).$$
Furthermore, those of level $\neq p$ are not congruence groups. It
is well known that $\Gamma$ contains infinitely many normal
subgroups of finite index which are not congruence groups.
However, all of these groups have the common feature that they are
lattice subgroups (in Rankin's terminology) of some normal
congruence group, and so are not maximal. Newman's result implies
that $\Gamma$ contains infinitely many maximal normal subgroups of
finite index which are not congruence groups, a somewhat
surprising fact. More precisely, Newman proved the following:

{\bf Theorem 3.6.} (see \cite{N}). {\it Suppose that the positive
integer $n$ satisfies

\roster

\item $n=p$ or $n | (p \pm 1)/2$,

\item $(n, 6)=1$,

\item $n>5$.

\endroster
Then there are elements $A, B \in PSL(2, p)$ such that $A$ is of
period $2$, $B$ is of period $3$, $AB$ is of period $n$, and
$PSL(2, p)=\langle A, B \rangle$. Suppose that $n$ satisfies
conditions $(1)$, $(2)$, $(3)$. Then there is a maximal normal
subgroup $G_n$ of $\Gamma$ such that $G_n$ is of level $n$, and
$\Gamma/G_n \cong PSL(2, p)$.}

  In fact, Newman constructed the homomorphism $\phi_n: \Gamma \to PSL(2, p)$
defined by
$$\phi_n: S \mapsto A, \quad ST \mapsto B,$$
where $S=\left(\matrix 0 & -1\\ 1 & 0 \endmatrix\right)$ and
$T=\left(\matrix 1 & 1\\ 0 & 1 \endmatrix\right)$, is actually a
homomorphism of $\Gamma$ onto $PSL(2, p)$. Let $G_n$ be the kernel
of $\phi_n$. Then $G_n$ is a normal subgroup of $\Gamma$ and
$\Gamma/G_n \cong PSL(2, p)$. Since $PSL(2, p)$ is simple, $G_n$
is a maximal normal subgroup of $\Gamma$. Furthermore the level of
$G_n$, which is the exponent of $T$ modulo $G_n$, is just $n$
since $\phi_n: T \mapsto AB$ and $AB$ is of period $n$. The groups
$G_n$ are certainly distinct, being of different levels, but are
all of index $\frac{1}{2}p(p^2-1)$ in $\Gamma$ with common
quotient group $PSL(2, p)$.

  In our case, $p=13$ and $n=7$ satisfy conditions (1), (2), (3)
and $n \neq p$. The group $G_7$ described by the above theorem is
thus a maximal normal subgroup of $\Gamma$ of index $1092$, and we
need only show that $G_7$ is not a congruence group. Suppose the
contrary. Then $G_7$, being of level $7$, would have to contain
the principal congruence subgroup $\Gamma(7)$, by Wohlfahrt's
theorem \cite{W1}. This would imply that $(\Gamma: G_7) | (\Gamma:
\Gamma(7))$, or that $1092|168$. But this is false. Thus, $G_7$ is
a noncongrunce group.

  We will show that for $\Gamma_{7, 7}^N$ or $G_7$ with $p=13$, there are
three noncongruence subgroups $G_1$, $G_2$, $G_3$ associated with
it. Let $z_i=x_i^{-1} y_i^{-1}$ ($i=1, 2, 3$). Then
$$y_i^2=x_i^3=z_i^7=y_i x_i z_i=1.\eqno{(3.35)}$$
Note that
$$\aligned
  z_1 &=PQ^2 P^{10} \cdot Q \cdot PQ^2P^{10},\\
  z_2 &=(Q^5 P^2 \cdot P^2 Q^6 P^8 \cdot Q^5 P^2) \cdot Q^2 \cdot (Q^5 P^2 \cdot P^2 Q^6 P^8 \cdot Q^5 P^2),\\
  z_3 &=Q^5 P^2 \cdot Q^4 \cdot Q^5 P^2,
\endaligned\eqno{(3.36)}$$
where $PQ^2 P^{10}$, $Q^5 P^2 \cdot P^2 Q^6 P^8 \cdot Q^5 P^2$ and
$Q^5 P^2$ are three involutions. Hence, $z_1$ is conjugate to $Q$,
$z_2$ is conjugate to $Q^2$ and $z_3$ is conjugate to $Q^4$, i.e.,
$z_1$, $z_2$ and $z_3$ are of type $7A$, $7B$ and $7C$,
respectively. Thus, we get a correspondence between the
presentations of $G$ and conjugacy classes of $G$:
$$(2, 3, 7; 7) \longleftrightarrow 7A, \quad
  (2, 3, 7; 6) \longleftrightarrow 7B, \quad
  (2, 3, 7; 13) \longleftrightarrow 7C.\eqno{(3.37)}$$

  A smooth complex projective curve of genus $g$ is a Hurwitz
curve if its automorphism group attains the maximum order
$84(g-1)$. A group that can be realized as the automorphism group
of a Hurwitz curve is called a Hurwitz group. Let $p \in {\Bbb Z}$
be a prime and let $q=p^n$. In \cite{M2}, Macbeath proved the
following theorem.

{\bf Theorem 3.7.} (see \cite{M2}). {\it The finite group $PSL(2,
q)$ is a Hurwitz group if:

\roster

\item $q=7$;

\item $q=p \equiv \pm 1$ mod $7$;

\item $q=p^3$, where $p \equiv \pm 2$ or $\pm 3$ mod $7$,
\endroster
and for no other values of $q$. In cases $(1)$ and $(3)$ there is
only one normal torsion free subgroup of the triangle group $G_{2,
3, 7}$ with quotient $PSL(2, q)$. In case $(2)$ there are three
different such subgroups leading to three non-isomorphic Riemann
surfaces $X_1$, $X_2$, $X_3$ with $PSL(2, q)$ acting as a Hurwitz
group of conformal automorphisms.}

  It is well known that there is a one-to-one correspondence
between the equivalence classes of compact Riemann surfaces and
the equivalence classes of algebraic curves. It is also true that
compact Riemann surfaces of genus $g>1$ with large automorphism
groups (in other words those Riemann surfaces with surface groups
which are normally contained in some triangle group) correspond to
algebraic curves defined over $\overline{{\Bbb Q}}$. Therefore we
have a natural action of Gal($\overline{{\Bbb Q}}/{\Bbb Q}$) on
these Riemann surfaces. For a Hurwitz curve $X$, one can consider
the moduli field of $X$, i.e. the fixed field of $\{ \sigma \in
\text{Gal}(\overline{{\Bbb Q}}/{\Bbb Q}): X \cong X^{\sigma} \}$.
In his paper \cite{St}, Streit proved the following:

{\bf Theorem 3.8.} (see \cite{St}). {\it In cases $(1)$ and $(3)$,
${\Bbb Q}$ is the moduli field and hence a minimal field of
definition of the Hurwitz curves with $PSL(2, q)$ acting as group
of conformal automorphisms. In case $(2)$ the moduli field of
$X_1$, $X_2$, $X_3$ is ${\Bbb Q}(\zeta_7+\zeta_7^{-1})$,
$\zeta_7:=\exp(2 \pi i/7)$. There are elements $\sigma, \tau \in$}
Gal{\it $(\overline{{\Bbb Q}}/{\Bbb Q})$ such that $X_1 \cong
X_2^{\sigma} \cong X_3^{\tau}$.}

  Let $G_{2, 3, 7}$ be the triangle group of orientation-preserving
transformations generated by reflections with angles $\pi/2$,
$\pi/3$ and $\pi/7$, respectively, in ${\Bbb H}$. Recall that
$G_{2, 3, 7}$ has a presentation given by
$$G_{2, 3, 7}=\langle t, u, v: t^2=u^3=v^7=tuv=1 \rangle.$$
What Macbeath proved in \cite{M2} is that there are three group
homomorphisms
$$q_i: G_{2, 3, 7} \twoheadrightarrow G$$
such that $q_i(v)$, with $i=1, 2, 3$, belong to three different
conjugacy classes of $G$ (in fact, the conjugacy classes $7A$,
$7B$ and $7C$ are fixed under the action of $\text{Aut}(G)$). As
compact Riemann surfaces, we have $X_i={\Bbb H}/N_i$ with
$N_i=\text{ker}(q_i)$ and ${\Bbb H}$ is the universal covering of
${\Bbb H}/N_i$. Thus, up to the permutation of three prime ideals,
the Shimura curve realizations of the first Hurwitz triplet, which
we denote by $X_1$, $X_2$ and $X_3$, correspond to the conjugacy
classes $7A$, $7B$ and $7C$, respectively. Now, we give the other
realizations as noncongruence modular curves.

  Let $\phi_i: \Gamma \to PGL(6, {\Bbb C})$ be three representations
where
$$\phi_i: S \mapsto y_i, \quad ST \mapsto x_i, \quad T^{-1} \mapsto z_i, \quad (i=1, 2, 3).\eqno{(3.38)}$$
Let $Y_i=\overline{{\Bbb H}/G_i}$ be the compactification of
${\Bbb H}/G_i$ where $G_i=\text{ker} \phi_i$. Then
$$\Gamma/G_1 \cong \Gamma/G_2 \cong \Gamma/G_3 \cong PSL(2, 13).\eqno{(3.39)}$$
$G_1$, $G_2$ and $G_3$ are non-congruence normal subgroups of
level $7$ of $\Gamma$. Note that, as a $(2, 3, 13)$-generated
group, by Riemann-Hurwitz formula, we have
$$2g-2=1092 \left(1-\frac{1}{2}-\frac{1}{3}-\frac{1}{13}\right).\eqno{(3.40)}$$
Hence, $g=50$, which is the genus of the modular curve $X(13)$. As
a $(2, 3, 7)$-generated group, by Riemann-Hurwitz formula, we have
$$2g-2=1092 \left(1-\frac{1}{2}-\frac{1}{3}-\frac{1}{7}\right).\eqno{(3.41)}$$
Hence $g=14$, which is of genus of ${\Bbb H}/G_i$ for $i=1, 2, 3$.
Therefore, $Y_1$, $Y_2$ and $Y_3$ must be Hurwitz curves.

  Following Wohlfahrt's paper \cite{W2}, we can give the other
approach to compute the genus of $Y_i$. Let $\Phi$ be a subgroup
of finite index $m$ in $\Gamma$ and ${\bold X}$ the associated
Riemann surface of genus $g$ and Hurwitz characteristic
$\chi=2g-2$. Then
$$6 \chi=m-6h-3 e_2-4 e_3,$$
where $h$ denotes the number of classes of parabolic, and $e_k$
the number of classes of elliptic fixed points, of order $k$, of
$\Phi$. In our case, $m=1092$, $e_2=0$, $e_3=0$, $h=156$. Hence,
$\chi=26$, $g=14$.

  By Theorem 3.7, there are only three Hurwitz curves with
genus $14$. Hence, $Y_1, Y_2, Y_3$ must be complex analytically
isomorphic to $X_1, X_2, X_3$. Note that $X_1, X_2, X_3$
correspond to the conjugacy classes $7A$, $7B$, $7C$,
respectively, and $Y_1, Y_2, Y_3$ also correspond to the conjugacy
classes $7A$, $7B$, $7C$, respectively. This implies that $Y_i$ is
complex analytically isomorphic to $X_i$ for $i=1, 2, 3$.

  It is known that (see \cite{St}) the special linear groups $PSL(2, q)$,
$q=p^n \neq 9$, prime $p>2$, may be generated by two elements
$PSL(2, q)=\langle g_0, g_1 \rangle$, such that $g_0$ is an
involution and $g_1$ is of order three. In order to be a Hurwitz
group, it is necessary and sufficient that one can choose the pair
$g_0$, $g_1$ such that $g_{\infty}:=(g_0 g_1)^{-1}$ has order $7$,
i.e. there is a surjective homomorphism $\varphi$ from the
triangle group $G_{2, 3, 7}$ to $PSL(2, q)$ preserving the orders
of the generators. We always assume this homomorphism $\varphi$ to
be realized by defining the image of three anti-clockwise
rotations $\gamma_0$, $\gamma_1$, $\gamma_{\infty}$ around the
non-accidental vertices of a canonical fundamental domain of
$G_{2, 3, 7}$, by $\varphi(\gamma_i)=g_i$, $i=0, 1, \infty$. The
kernel $\Gamma=\text{ker}(\varphi)$ defines a surface group
$\Gamma \triangleleft G_{2, 3, 7}$ and a Riemann surface $X=\Gamma
\backslash {\Bbb H}$, where $\text{Aut}(X)=PSL(2, q)$ and ${\Bbb
H}$ is the upper half plane.

{\bf Theorem 3.9.} (see \cite{St}, Theorem 3). {\it Assume $n \in
{\Bbb N}$, $n>6$ and prime $p \geq 13$. Let $p \equiv \pm 1$ mod
$n$ and $\psi: (2, 3, n) \longrightarrow PSL(2, p)$ an epimorphism
which maps $\gamma_0$, $\gamma_1$, $\gamma_{\infty}$ to $g_0$,
$g_1$, $g_{\infty}$ and which preserves the orders of the
generators of $(2, 3, n)$. Then $\text{ker}(\psi) \backslash {\Bbb
H}$ is a Riemann surface defined over ${\Bbb
Q}(\zeta_n+\zeta_n^{-1})$ as a minimal field of definition with
$\zeta_n:=\exp(2 \pi i/n)$. There are exactly $\varphi(n)/2$
non-biholomorphically equivalent but $\text{Gal}(\overline{{\Bbb
Q}}/{\Bbb Q})$-conjugate Riemann surfaces which can be constructed
in this manner, where $\varphi$ is the Euler $\varphi$-function.}

  The constructed epimorphisms lead to biholomorphically equivalent
Riemann surfaces when replacing $(2, 3, n)$ by $(2, 3, \infty)$.
Therefore the kernel of these epimorphisms are then subgroups of
the modular group. As Macbeath pointed out that they can not be
congruence subgroups. Therefore, in \cite{St} Streit constructed
the Galois orbit of quite a few normal non-congruence subgroups
and their minimal field of definition.

  In our case, $n=7$, $p=13$. The kernel of our three epimorphisms
are non-congruence subgroups of level $7$. Their minimal field of
definition is ${\Bbb Q}(\cos \frac{2 \pi}{7})$.

  In \cite{V}, Vogeler calculated the systoles for the first
Hurwitz triplet. Now, we give a connection between his results and
our three models. For genus $14A$, which corresponds to our second
model, i.e., the presentation $(2, 3, 7; 6)$, the length spectrum
$n=6$, the systole is $5.903919373$, the number of systolic loops
is $91$. For genus $14B$, which corresponds to our first model,
i.e., the presentation $(2, 3, 7; 7)$, the length spectrum $n=7$,
the systole is $6.887905935$, the number of systolic loops is
$78$. For genus $14C$, which corresponds to our third model, i.e.,
the presentation $(2, 3, 7; 13)$, the length spectrum $n=3$, the
systole is $6.393315042$, the number of systolic loops is $364$.

  Recall that the theta functions with characteristic
$\left[\matrix \epsilon\\ \epsilon^{\prime} \endmatrix\right] \in
{\Bbb R}^2$ is defined by the following series which converges
uniformly and absolutely on compact subsets of ${\Bbb C} \times
{\Bbb H}$ (see \cite{FK}, p. 73):
$$\theta \left[\matrix \epsilon\\ \epsilon^{\prime} \endmatrix\right] (z, \tau)
 =\sum_{n \in {\Bbb Z}} \exp \left\{2 \pi i \left[\frac{1}{2}
  \left(n+\frac{\epsilon}{2}\right)^2 \tau+\left(n+\frac{\epsilon}{2}\right)
  \left(z+\frac{\epsilon^{\prime}}{2}\right)\right]\right\}.$$
The modified theta constants is defined by (see \cite{FK}, p. 215)
$$\varphi_l(\tau)=\theta [\chi_l](0, k \tau),$$
where the characteristic $\chi_l=\left[\matrix \frac{2l+1}{k}\\ 1
\endmatrix\right]$, $l=0, \cdots, \frac{k-3}{2}$, for odd $k$ and
$\chi_l=\left[\matrix \frac{2l}{k}\\ 0 \endmatrix\right]$, $l=0,
\cdots, \frac{k}{2}$, for even $k$. We have the following:

{\bf Theorem 3.10.} (see \cite{FK}, p. 236). {\it For each odd
integer $k \geq 5$, the map
$$\Phi: \tau \mapsto (\varphi_0(\tau), \varphi_1(\tau), \cdots,
  \varphi_{\frac{k-5}{2}}(\tau), \varphi_{\frac{k-3}{2}}(\tau))$$
from ${\Bbb H} \cup {\Bbb Q} \cup \{ \infty \}$ to ${\Bbb
C}^{\frac{k-1}{2}}$, defines a holomorphic mapping from
$\overline{{\Bbb H}/\Gamma(k)}$ into ${\Bbb C} {\Bbb
P}^{\frac{k-3}{2}}$.}

  In our case, $k=13$, the map
$$\Phi: \tau \mapsto (\varphi_0(\tau), \varphi_1(\tau), \varphi_2(\tau),
  \varphi_3(\tau), \varphi_4(\tau), \varphi_5(\tau))$$
gives a holomorphic mapping from the modular curve
$X(13)=\overline{{\Bbb H}/\Gamma(13)}$ into ${\Bbb C} {\Bbb P}^5$,
which corresponds to our six-dimensional representation, i.e., up
to the constants, $z_1, \cdots, z_6$ are just modular forms
$\varphi_0(\tau), \cdots, \varphi_5(\tau)$.

\vskip 0.5 cm

\centerline{\bf 4. Seven-dimensional representations of $PSL(2,
                   13)$,}

\centerline{\bf Jacobian modular equation of degree fourteen and
                exotic duality theorem}

\vskip 0.5 cm

  Let us study the action of $S T^{\nu}$ on the five dimensional
projective space ${\Bbb P}^5=\{(z_1, z_2, z_3, z_4, z_5, z_6)\}$,
where $\nu=0, 1, \cdots, 12$. Put
$$\alpha=\zeta+\zeta^{12}-\zeta^5-\zeta^8, \quad
   \beta=\zeta^3+\zeta^{10}-\zeta^2-\zeta^{11}, \quad
   \gamma=\zeta^9+\zeta^4-\zeta^6-\zeta^7.$$
We find that
$$\aligned
  &13 ST^{\nu}(z_1) \cdot ST^{\nu}(z_4)\\
=&\beta z_1 z_4+\gamma z_2 z_5+\alpha z_3 z_6+\\
 &+\gamma \zeta^{\nu} z_1^2+\alpha \zeta^{9 \nu} z_2^2+\beta \zeta^{3 \nu} z_3^2
  -\gamma \zeta^{12 \nu} z_4^2-\alpha \zeta^{4 \nu} z_5^2-\beta \zeta^{10 \nu} z_6^2+\\
 &+(\alpha-\beta) \zeta^{5 \nu} z_1 z_2+(\beta-\gamma) \zeta^{6 \nu} z_2 z_3
  +(\gamma-\alpha) \zeta^{2 \nu} z_1 z_3+\\
 &+(\beta-\alpha) \zeta^{8 \nu} z_4 z_5+(\gamma-\beta) \zeta^{7 \nu} z_5 z_6
  +(\alpha-\gamma) \zeta^{11 \nu} z_4 z_6+\\
 &-(\alpha+\beta) \zeta^{\nu} z_3 z_4-(\beta+\gamma) \zeta^{9 \nu} z_1 z_5
  -(\gamma+\alpha) \zeta^{3 \nu} z_2 z_6+\\
 &-(\alpha+\beta) \zeta^{12 \nu} z_1 z_6-(\beta+\gamma) \zeta^{4 \nu} z_2 z_4
  -(\gamma+\alpha) \zeta^{10 \nu} z_3 z_5.
\endaligned$$
$$\aligned
  &13 ST^{\nu}(z_2) \cdot ST^{\nu}(z_5)\\
=&\gamma z_1 z_4+\alpha z_2 z_5+\beta z_3 z_6+\\
 &+\alpha \zeta^{\nu} z_1^2+\beta \zeta^{9 \nu} z_2^2+\gamma \zeta^{3 \nu} z_3^2
  -\alpha \zeta^{12 \nu} z_4^2-\beta \zeta^{4 \nu} z_5^2-\gamma \zeta^{10 \nu} z_6^2+\\
 &+(\beta-\gamma) \zeta^{5 \nu} z_1 z_2+(\gamma-\alpha) \zeta^{6 \nu} z_2 z_3
  +(\alpha-\beta) \zeta^{2 \nu} z_1 z_3+\\
 &+(\gamma-\beta) \zeta^{8 \nu} z_4 z_5+(\alpha-\gamma) \zeta^{7 \nu} z_5 z_6
  +(\beta-\alpha) \zeta^{11 \nu} z_4 z_6+\\
 &-(\beta+\gamma) \zeta^{\nu} z_3 z_4-(\gamma+\alpha) \zeta^{9 \nu} z_1 z_5
  -(\alpha+\beta) \zeta^{3 \nu} z_2 z_6+\\
 &-(\beta+\gamma) \zeta^{12 \nu} z_1 z_6-(\gamma+\alpha) \zeta^{4 \nu} z_2 z_4
  -(\alpha+\beta) \zeta^{10 \nu} z_3 z_5.
\endaligned$$
$$\aligned
  &13 ST^{\nu}(z_3) \cdot ST^{\nu}(z_6)\\
=&\alpha z_1 z_4+\beta z_2 z_5+\gamma z_3 z_6+\\
 &+\beta \zeta^{\nu} z_1^2+\gamma \zeta^{9 \nu} z_2^2+\alpha \zeta^{3 \nu} z_3^2
  -\beta \zeta^{12 \nu} z_4^2-\gamma \zeta^{4 \nu} z_5^2-\alpha \zeta^{10 \nu} z_6^2+\\
 &+(\gamma-\alpha) \zeta^{5 \nu} z_1 z_2+(\alpha-\beta) \zeta^{6 \nu} z_2 z_3
  +(\beta-\gamma) \zeta^{2 \nu} z_1 z_3+\\
 &+(\alpha-\gamma) \zeta^{8 \nu} z_4 z_5+(\beta-\alpha) \zeta^{7 \nu} z_5 z_6
  +(\gamma-\beta) \zeta^{11 \nu} z_4 z_6+\\
 &-(\gamma+\alpha) \zeta^{\nu} z_3 z_4-(\alpha+\beta) \zeta^{9 \nu} z_1 z_5
  -(\beta+\gamma) \zeta^{3 \nu} z_2 z_6+\\
 &-(\gamma+\alpha) \zeta^{12 \nu} z_1 z_6-(\alpha+\beta) \zeta^{4 \nu} z_2 z_4
  -(\beta+\gamma) \zeta^{10 \nu} z_3 z_5.
\endaligned$$
Note that $\alpha+\beta+\gamma=\sqrt{13}$, we find that
$$\aligned
  &\sqrt{13} \left[ST^{\nu}(z_1) \cdot ST^{\nu}(z_4)+ST^{\nu}(z_2) \cdot ST^{\nu}(z_5)+ST^{\nu}(z_3) \cdot ST^{\nu}(z_6)\right]\\
 =&(z_1 z_4+z_2 z_5+z_3 z_6)+(\zeta^{\nu} z_1^2+\zeta^{9 \nu} z_2^2+\zeta^{3 \nu} z_3^2)
  -(\zeta^{12 \nu} z_4^2+\zeta^{4 \nu} z_5^2+\zeta^{10 \nu} z_6^2)+\\
  &-2(\zeta^{\nu} z_3 z_4+\zeta^{9 \nu} z_1 z_5+\zeta^{3 \nu} z_2 z_6)
   -2(\zeta^{12 \nu} z_1 z_6+\zeta^{4 \nu} z_2 z_4+\zeta^{10 \nu} z_3 z_5).
\endaligned$$
Let
$$\varphi_{\infty}(z_1, z_2, z_3, z_4, z_5, z_6)=\sqrt{13} (z_1 z_4+z_2 z_5+z_3 z_6)\eqno{(4.1)}$$
and
$$\varphi_{\nu}(z_1, z_2, z_3, z_4, z_5, z_6)=\varphi_{\infty}(ST^{\nu}(z_1, z_2, z_3, z_4, z_5, z_6))\eqno{(4.2)}$$
for $\nu=0, 1, \cdots, 12$. Then
$$\aligned
  \varphi_{\nu}
=&(z_1 z_4+z_2 z_5+z_3 z_6)+\zeta^{\nu} (z_1^2-2 z_3 z_4)+\zeta^{4 \nu} (-z_5^2-2 z_2 z_4)+\\
 &+\zeta^{9 \nu} (z_2^2-2 z_1 z_5)+\zeta^{3 \nu} (z_3^2-2 z_2 z_6)+
   \zeta^{12 \nu} (-z_4^2-2 z_1 z_6)+\zeta^{10 \nu} (-z_6^2-2 z_3 z_5).
\endaligned\eqno{(4.3)}$$
This leads us to define the following senary quadratic forms
(quadratic forms in six variables):
$$\left\{\aligned
  {\Bbb A}_0 &=z_1 z_4+z_2 z_5+z_3 z_6,\\
  {\Bbb A}_1 &=z_1^2-2 z_3 z_4,\\
  {\Bbb A}_2 &=-z_5^2-2 z_2 z_4,\\
  {\Bbb A}_3 &=z_2^2-2 z_1 z_5,\\
  {\Bbb A}_4 &=z_3^2-2 z_2 z_6,\\
  {\Bbb A}_5 &=-z_4^2-2 z_1 z_6,\\
  {\Bbb A}_6 &=-z_6^2-2 z_3 z_5.
\endaligned\right.\eqno{(4.4)}$$
Hence,
$$\sqrt{13} ST^{\nu}({\Bbb A}_0)={\Bbb A}_0+\zeta^{\nu} {\Bbb A}_1+\zeta^{4 \nu} {\Bbb A}_2+
  \zeta^{9 \nu} {\Bbb A}_3+\zeta^{3 \nu} {\Bbb A}_4+\zeta^{12 \nu} {\Bbb A}_5+\zeta^{10 \nu} {\Bbb A}_6.\eqno{(4.5)}$$
Let
$$\left\{\aligned
  p_1 &=\zeta^2+\zeta^{11}-2+2(\zeta+\zeta^{12}-\zeta^9-\zeta^4),\\
  p_2 &=2-\zeta^9-\zeta^4+2(\zeta^5+\zeta^8-\zeta^2-\zeta^{11}),\\
  p_3 &=\zeta^6+\zeta^7-2+2(\zeta^3+\zeta^{10}-\zeta-\zeta^{12}),\\
  p_4 &=\zeta^5+\zeta^8-2+2(\zeta^9+\zeta^4-\zeta^3-\zeta^{10}),\\
  p_5 &=2-\zeta^3-\zeta^{10}+2(\zeta^6+\zeta^7-\zeta^5-\zeta^8),\\
  p_6 &=2-\zeta-\zeta^{12}+2(\zeta^2+\zeta^{11}-\zeta^6-\zeta^7).
\endaligned\right.\eqno{(4.6)}$$
We find that
$$\left\{\aligned
  13 S({\Bbb A}_1) &=2 \sqrt{13} {\Bbb A}_0+p_1 {\Bbb A}_1+p_2 {\Bbb A}_2+
                     p_3 {\Bbb A}_3+p_4 {\Bbb A}_4+p_5 {\Bbb A}_5+p_6 {\Bbb A}_6,\\
  13 S({\Bbb A}_2) &=2 \sqrt{13} {\Bbb A}_0+p_2 {\Bbb A}_1+p_4 {\Bbb A}_2+
                     p_6 {\Bbb A}_3+p_5 {\Bbb A}_4+p_3 {\Bbb A}_5+p_1 {\Bbb A}_6,\\
  13 S({\Bbb A}_3) &=2 \sqrt{13} {\Bbb A}_0+p_3 {\Bbb A}_1+p_6 {\Bbb A}_2+
                     p_4 {\Bbb A}_3+p_1 {\Bbb A}_4+p_2 {\Bbb A}_5+p_5 {\Bbb A}_6,\\
  13 S({\Bbb A}_4) &=2 \sqrt{13} {\Bbb A}_0+p_4 {\Bbb A}_1+p_5 {\Bbb A}_2+
                     p_1 {\Bbb A}_3+p_3 {\Bbb A}_4+p_6 {\Bbb A}_5+p_2 {\Bbb A}_6,\\
  13 S({\Bbb A}_5) &=2 \sqrt{13} {\Bbb A}_0+p_5 {\Bbb A}_1+p_3 {\Bbb A}_2+
                     p_2 {\Bbb A}_3+p_6 {\Bbb A}_4+p_1 {\Bbb A}_5+p_4 {\Bbb A}_6,\\
  13 S({\Bbb A}_6) &=2 \sqrt{13} {\Bbb A}_0+p_6 {\Bbb A}_1+p_1 {\Bbb A}_2+
                     p_5 {\Bbb A}_3+p_2 {\Bbb A}_4+p_4 {\Bbb A}_5+p_3 {\Bbb A}_6,\\
\endaligned\right.\eqno{(4.7)}$$
Note that
$$\left\{\aligned
  p_1 &=\sqrt{13} (\zeta^2+\zeta^{11}), \\
  p_2 &=\sqrt{13} (\zeta^9+\zeta^4),\\
  p_3 &=\sqrt{13} (\zeta^6+\zeta^7),\\
  p_4 &=\sqrt{13} (\zeta^5+\zeta^8),\\
  p_5 &=\sqrt{13} (\zeta^3+\zeta^{10}),\\
  p_6 &=\sqrt{13} (\zeta+\zeta^{12}).
\endaligned\right.\eqno{(4.8)}$$
We obtain a seven-dimensional representation of the simple group
$PSL(2, 13) \cong \langle \widetilde{S}, \widetilde{T} \rangle$
which induces from the actions of $S$ and $T$ on the basis $({\Bbb
A}_0, {\Bbb A}_1, {\Bbb A}_2, {\Bbb A}_3, {\Bbb A}_4, {\Bbb A}_5,
{\Bbb A}_6)$. Here
$$\widetilde{S}=\frac{1}{\sqrt{13}} \left(\matrix
  1 & 1                  & 1               & 1               & 1               & 1                  & 1\\
  2 & \zeta^2+\zeta^{11} & \zeta^9+\zeta^4 & \zeta^6+\zeta^7 & \zeta^5+\zeta^8 & \zeta^3+\zeta^{10} & \zeta+\zeta^{12}\\
  2 & \zeta^9+\zeta^4 & \zeta^5+\zeta^8 & \zeta+\zeta^{12} & \zeta^3+\zeta^{10} & \zeta^6+\zeta^7 & \zeta^2+\zeta^{11}\\
  2 & \zeta^6+\zeta^7 & \zeta+\zeta^{12} & \zeta^5+\zeta^8 & \zeta^2+\zeta^{11} & \zeta^9+\zeta^4 & \zeta^3+\zeta^{10}\\
  2 & \zeta^5+\zeta^8 & \zeta^3+\zeta^{10} & \zeta^2+\zeta^{11} & \zeta^6+\zeta^7 & \zeta+\zeta^{12} & \zeta^9+\zeta^4\\
  2 & \zeta^3+\zeta^{10} & \zeta^6+\zeta^7 & \zeta^9+\zeta^4 & \zeta+\zeta^{12} & \zeta^2+\zeta^{11} & \zeta^5+\zeta^8\\
  2 & \zeta+\zeta^{12} & \zeta^2+\zeta^{11} & \zeta^3+\zeta^{10} & \zeta^9+\zeta^4 & \zeta^5+\zeta^8 & \zeta^6+\zeta^7
\endmatrix\right),\eqno{(4.9)}$$
and
$$\widetilde{T}=\left(\matrix
  1 &       &         &         &         &            &           \\
    & \zeta &         &         &         &            &           \\
    &       & \zeta^4 &         &         &            &           \\
    &       &         & \zeta^9 &         &            &           \\
    &       &         &         & \zeta^3 &            &           \\
    &       &         &         &         & \zeta^{12} &           \\
    &       &         &         &         &            & \zeta^{10}
\endmatrix\right).\eqno{(4.10)}$$
We have
$$\text{Tr}(\widetilde{S})=-1, \quad
  \text{Tr}(\widetilde{T})=\frac{1+\sqrt{13}}{2}, \quad
  \text{Tr}(\widetilde{S} \widetilde{T})=1.\eqno{(4.11)}$$
Hence, our seven-dimensional representation corresponds to the
character $\chi_3$ in Table $1$. In fact, for a prime $q \equiv 1$
(mod $4$), Klein and Hecke obtained an $r=\frac{1}{2} (q+1)$ rowed
matrix representation of the finite group $\Gamma/\Gamma(q)$ by
means of matrices whose elements are in the cyclotomic field
generated by $e^{2 \pi i/q}$ (see \cite{K3} and \cite{He}). When
$q=13$, it is just our seven dimensional representation. However,
our seven dimensional representation is induced from a six
dimensional representation, which does not appear in Klein and
Hecke's papers \cite{K3} and \cite{He}. Furthermore, our
representation involves invariants ${\Bbb A}_0$, $\cdots$, ${\Bbb
A}_6$, which they did not study.

  Now, we study the ``triality'' associated with $PSL(2, 13)$,
i.e., there is an automorphism of order three as follows: besides
$\langle S, T \rangle$, there are two other generators for $PSL(2,
13)$.
$$S_1=-\frac{1}{\sqrt{13}} \left(\matrix
  \zeta^4-\zeta^9 & \zeta^{12}-\zeta & \zeta^{10}-\zeta^3 & \zeta^5-\zeta^8 & \zeta^2-\zeta^{11} & \zeta^6-\zeta^7\\
  \zeta^{12}-\zeta & \zeta^{10}-\zeta^3 & \zeta^4-\zeta^9 & \zeta^2-\zeta^{11} & \zeta^6-\zeta^7 & \zeta^5-\zeta^8\\
  \zeta^{10}-\zeta^3 & \zeta^4-\zeta^9 & \zeta^{12}-\zeta & \zeta^6-\zeta^7 & \zeta^5-\zeta^8 & \zeta^2-\zeta^{11}\\
  \zeta^5-\zeta^8 & \zeta^2-\zeta^{11} & \zeta^6-\zeta^7 & \zeta^3-\zeta^{10} & \zeta^9-\zeta^4 & \zeta-\zeta^{12}\\
  \zeta^2-\zeta^{11} & \zeta^6-\zeta^7 & \zeta^5-\zeta^8 & \zeta^9-\zeta^4 & \zeta-\zeta^{12} & \zeta^3-\zeta^{10}\\
  \zeta^6-\zeta^7 & \zeta^5-\zeta^8 & \zeta^2-\zeta^{11} & \zeta-\zeta^{12} & \zeta^3-\zeta^{10} & \zeta^9-\zeta^4
\endmatrix\right),\eqno{(4.12)}$$
$$T_1=\text{diag}(\zeta^{11}, \zeta^8, \zeta^7, \zeta^5, \zeta^6, \zeta^2),\eqno{(4.13)}$$
$$S_2=-\frac{1}{\sqrt{13}} \left(\matrix
  \zeta^{10}-\zeta^3 & \zeta^4-\zeta^9 & \zeta^{12}-\zeta & \zeta^5-\zeta^8 & \zeta^2-\zeta^{11} & \zeta^6-\zeta^7\\
  \zeta^4-\zeta^9 & \zeta^{12}-\zeta & \zeta^{10}-\zeta^3 & \zeta^2-\zeta^{11} & \zeta^6-\zeta^7 & \zeta^5-\zeta^8\\
  \zeta^{12}-\zeta & \zeta^{10}-\zeta^3 & \zeta^4-\zeta^9 & \zeta^6-\zeta^7 & \zeta^5-\zeta^8 & \zeta^2-\zeta^{11}\\
  \zeta^5-\zeta^8 & \zeta^2-\zeta^{11} & \zeta^6-\zeta^7 & \zeta^9-\zeta^4 & \zeta-\zeta^{12} & \zeta^3-\zeta^{10}\\
  \zeta^2-\zeta^{11} & \zeta^6-\zeta^7 & \zeta^5-\zeta^8 & \zeta-\zeta^{12} & \zeta^3-\zeta^{10} & \zeta^9-\zeta^4\\
  \zeta^6-\zeta^7 & \zeta^5-\zeta^8 & \zeta^2-\zeta^{11} & \zeta^3-\zeta^{10} & \zeta^9-\zeta^4 & \zeta-\zeta^{12}
\endmatrix\right)\eqno{(4.14)}$$
and
$$T_2=\text{diag}(\zeta^8, \zeta^7, \zeta^{11}, \zeta^2, \zeta^5, \zeta^6).\eqno{(4.15)}$$
Let
$$R=\left(\matrix
    0 & 0 & 1 &   &   &  \\
    1 & 0 & 0 &   &   &  \\
    0 & 1 & 0 &   &   &  \\
      &   &   & 0 & 1 & 0\\
      &   &   & 0 & 0 & 1\\
      &   &   & 1 & 0 & 0
   \endmatrix\right).\eqno{(4.16)}$$
Then $R^3=1$ and $R T_1 R^{-1}=T$, $R S_1 R^{-1}=S$, $R^{-1} T_2
R=T$ and $R^{-1} S_2 R=S$. Hence, $\langle S_1, T_1 \rangle \cong
\langle S_2, T_2 \rangle \cong PSL(2, 13)$.

  Let
$$\psi_{\infty}(z_1, z_2, z_3, z_4, z_5, z_6)=\sqrt{13} (z_1 z_6+z_2 z_4+z_3 z_5)\eqno{(4.17)}$$
and
$$\psi_{\nu}(z_1, z_2, z_3, z_4, z_5, z_6)=\psi_{\infty}(S_1 T_1^{\nu}(z_1, z_2, z_3, z_4, z_5, z_6))\eqno{(4.18)}$$
for $\nu=0, 1, \cdots, 12$. Then
$$\psi_{\nu}
 ={\Bbb B}_0+\zeta^{\nu} {\Bbb B}_1+\zeta^{4 \nu} {\Bbb B}_2
 +\zeta^{9 \nu} {\Bbb B}_3+\zeta^{3 \nu} {\Bbb B}_4+\zeta^{12 \nu} {\Bbb B}_5
 +\zeta^{10 \nu} {\Bbb B}_6,\eqno{(4.19)}$$
where
$$\left\{\aligned
  {\Bbb B}_0 &=z_1 z_6+z_2 z_4+z_3 z_5,\\
  {\Bbb B}_1 &=z_3^2-2 z_2 z_5,\\
  {\Bbb B}_2 &=-z_6^2-2 z_1 z_5,\\
  {\Bbb B}_3 &=z_1^2-2 z_3 z_6,\\
  {\Bbb B}_4 &=z_2^2-2 z_1 z_4,\\
  {\Bbb B}_5 &=-z_5^2-2 z_3 z_4,\\
  {\Bbb B}_6 &=-z_4^2-2 z_2 z_6.
\endaligned\right.\eqno{(4.20)}$$

  Let
$$\phi_{\infty}(z_1, z_2, z_3, z_4, z_5, z_6)=\sqrt{13} (z_1 z_5+z_2 z_6+z_3 z_4)\eqno{(4.21)}$$
and
$$\phi_{\nu}(z_1, z_2, z_3, z_4, z_5, z_6)=\phi_{\infty}(S_2 T_2^{\nu}(z_1, z_2, z_3, z_4, z_5, z_6))\eqno{(4.22)}$$
for $\nu=0, 1, \cdots, 12$. Then
$$\phi_{\nu}
 ={\Bbb C}_0+\zeta^{\nu} {\Bbb C}_1+\zeta^{4 \nu} {\Bbb C}_2
 +\zeta^{9 \nu} {\Bbb C}_3+\zeta^{3 \nu} {\Bbb C}_4+\zeta^{12 \nu} {\Bbb C}_5
 +\zeta^{10 \nu} {\Bbb C}_6,\eqno{(4.23)}$$
where
$$\left\{\aligned
  {\Bbb C}_0 &=z_1 z_5+z_2 z_6+z_3 z_4,\\
  {\Bbb C}_1 &=z_2^2-2 z_1 z_6,\\
  {\Bbb C}_2 &=-z_4^2-2 z_3 z_6,\\
  {\Bbb C}_3 &=z_3^2-2 z_2 z_4,\\
  {\Bbb C}_4 &=z_1^2-2 z_3 z_5,\\
  {\Bbb C}_5 &=-z_6^2-2 z_2 z_5,\\
  {\Bbb C}_6 &=-z_5^2-2 z_1 z_4.
\endaligned\right.\eqno{(4.24)}$$

  Similar as above, we obtain the other two kinds of seven-dimensional
representations of the simple group $PSL(2, 13) \cong \langle
\widetilde{S_1}, \widetilde{T_1} \rangle \cong \langle
\widetilde{S_2}, \widetilde{T_2} \rangle$ which induce from the
actions of $S_1$ and $T_1$ on the basis $({\Bbb B}_0, {\Bbb B}_1,
{\Bbb B}_2, {\Bbb B}_3, {\Bbb B}_4, {\Bbb B}_5, {\Bbb B}_6)$ and
the actions of $S_2$ and $T_2$ on the basis $({\Bbb C}_0, {\Bbb
C}_1, {\Bbb C}_2, {\Bbb C}_3, {\Bbb C}_4, {\Bbb C}_5, {\Bbb
C}_6)$, respectively. Remarkably, we find that
$$\widetilde{S_1}=\widetilde{S_2}=\widetilde{S}, \quad
  \widetilde{T_1}=\widetilde{T_2}=\widetilde{T}!\eqno{(4.25)}$$

  Now, let us recall some facts about theta functions over number
fields and Hilbert modular forms (see \cite{E} and \cite{Hi}, pp.
796--798). Consider the field $K={\Bbb Q}(\zeta)$ where
$\zeta=e^{2 \pi i/p}$. Let $k={\Bbb Q}(\zeta+\zeta^{-1})$ be the
real subfield. Let ${\frak O}$ be the ring of integers in $K$, and
let ${\frak P}$ be the principal ideal of ${\frak O}$ generated by
the element $1-\zeta$. Since $k$ is the real subfield of $K$ and
$[k: {\Bbb Q}]=\frac{p-1}{2}$, there exist exactly $\frac{p-1}{2}$
distinct real embeddings $\sigma_l: k \to {\Bbb R}$, $l=1, \cdots,
\frac{p-1}{2}$. each of these $\sigma_l$ is of the form
$\zeta+\zeta^{-1} \mapsto \zeta^a+\zeta^{-a}$ for a suitable
integer $a$. In particular, one has $\sigma_l(k)=k$. Hence the
$\sigma_l$ form a group (with respect to composition), the Galois
group of $k$ over ${\Bbb Q}$. We denote it by $G$. Consider the
product
$${\Bbb H}^{\frac{p-1}{2}}={\Bbb H} \times \cdots \times {\Bbb H}
  \quad (\text{$(p-1)/2$ times})$$
of $\frac{p-1}{2}$ upper half planes. Let $z=(z_1, \cdots,
z_{(p-1)/2})$ be a point of ${\Bbb H}^{\frac{p-1}{2}}$. We define
the theta function $\theta_j(z)$ depending on $\frac{p-1}{2}$
variables $z_l \in {\Bbb H}$ by
$$\theta_j(z):=\sum_{x \in {\frak P}+j} e^{2 \pi i
\text{Tr}_{k/{\Bbb Q}}\left(z \frac{x \overline{x}}{p}\right)},$$
for $j=0, 1, \cdots, \frac{p-1}{2}$, where
$$\text{Tr}_{k/{\Bbb Q}} \left(z \frac{x \overline{x}}{p}\right)
 :=\sum_{l=1}^{(p-1)/2} z_l \cdot \frac{\sigma_l(x \overline{x})}{p}.$$
The function $\theta_j$ is holomorphic in $z \in {\Bbb
H}^{(p-1)/2}$.

  Let ${\frak O}_k$ be the ring of integers in $k$. The group
$SL_2({\frak O}_k)$ is the group of all $2 \times 2$-matrices
$$\left(\matrix \alpha & \beta\\
                \gamma & \delta
  \endmatrix\right)$$
with entries $\alpha, \beta, \gamma, \delta \in {\frak O}_k$ and
with determinant $\alpha \delta-\beta \gamma=1$. This group
operates on ${\Bbb H}^{(p-1)/2}$ by
$$z \mapsto \frac{\alpha z+\beta}{\gamma z+\delta}, \quad
  z_l \mapsto \frac{\sigma_l(\alpha) z_l+\sigma_l(\beta)}
  {\sigma_l(\gamma) z_l+\sigma_l(\delta)}, \quad
  l=1, \cdots, \frac{p-1}{2}.$$
The norm $N_{k/{\Bbb Q}}(\alpha)$ of an element $\alpha \in k$
over ${\Bbb Q}$ is defined by
$$N_{k/{\Bbb Q}}(\alpha):=\prod_{l=1}^{(p-1)/2} \sigma_l(\alpha).$$
For $z \in {\Bbb H}^{(p-1)/2}$, $\gamma, \delta \in {\frak O}_k$
we define
$$N_{k/{\Bbb Q}}(\gamma z+\delta):=\prod_{l=1}^{(p-1)/2} (\sigma_l(\gamma) z_l+\sigma_l(\delta)).$$
If $\sigma \in G$ we set $\sigma(z)=(z_{\varepsilon(1)}, \cdots,
z_{\varepsilon\left(\frac{p-1}{2}\right)})$, where $\varepsilon$
denotes that permutation of the indices $1, \cdots, \frac{p-1}{2}$
such that $\sigma_l \circ \sigma=\sigma_{\varepsilon(l)}$ for $1
\leq l \leq \frac{p-1}{2}$. Finally let $\Gamma$ be a subgroup of
$SL_2({\frak O}_k)$.

{\it Definition} 4.1. A holomorphic function $f: {\Bbb
H}^{(p-1)/2} \to {\Bbb C}$ is called a Hilbert modular form of
weight $m$ for $\Gamma$, if
$$f \left(\frac{\alpha z+\beta}{\gamma z+\delta}\right)
 =f(z) \cdot N_{k/{\Bbb Q}}(\gamma z+\delta)^m \quad
 \text{for all $\left(\matrix \alpha & \beta\\ \gamma & \delta
 \endmatrix\right) \in \Gamma$}.$$
It is called symmetric, if $f(\sigma(z))=f(z)$ for all $\sigma \in
G$.

  Let ${\frak p}$ be the ideal ${\frak p}:={\frak P} \cap {\frak
O}_k$ of ${\frak O}_k$. Then
$${\frak p}=(\zeta+\zeta^{-1}-2)=((\zeta-1)(\zeta^{-1}-1)).$$
Moreover,
$${\frak p}^{\frac{p-1}{2}}=(p).$$
We define
$$\Gamma({\frak p}):=\left\{\left(\matrix \alpha & \beta\\
  \gamma & \delta \endmatrix\right) \in SL_2({\frak O}_k):
  \text{$\alpha \equiv \delta \equiv 1$ (mod ${\frak p}$}),
  \text{$\beta \equiv \gamma \equiv 0$ (mod ${\frak p}$})\right\},$$
$$\Gamma_0({\frak p}):=\left\{\left(\matrix \alpha & \beta\\
  \gamma & \delta \endmatrix\right) \in SL_2({\frak O}_k):
  \text{$\gamma \equiv 0$ (mod ${\frak p}$}) \right\}.$$
Then we have the following result.

{\bf Theorem 4.2.} (see \cite{E}). {\it The function $\theta_j$,
$j=0, 1, \cdots, \frac{p-1}{2}$, is a Hilbert modular form of
weight $1$ for the group $\Gamma({\frak p})$. Moreover one has
$$\theta_0 \left(\frac{\alpha z+\beta}{\gamma z+\delta}\right)
 =\theta_0(z) \cdot \left(\frac{\delta}{p}\right) \cdot
  N_{k/{\Bbb Q}}(\gamma z+\delta) \quad \text{for all
  $\left(\matrix \alpha & \beta\\
  \gamma & \delta \endmatrix\right) \in \Gamma_0({\frak p})$}.$$}

  Note that the group $SL(2, {\Bbb F}_p)$ acts on the $\frac{p+1}{2}$-dimensional
vector space over ${\Bbb C}$ generated by the $\theta_j$. For $p
\equiv 1$ mod $4$ this is an action of $PSL(2, {\Bbb F}_p)$.

  Now let $C \subset {\Bbb F}_p^n$ be a self-dual code. We have $n(p-1) \equiv 0$ (mod $8$).
Let $\Gamma_{C}$ be the lattice constructed from $C$. Then
$\Gamma_C$ is an even unimodular lattice of rank $n(p-1)$. By the
definition of the symmetric bilinear form on $\Gamma_C$, the usual
theta function of the lattice $\Gamma_C$ is the function
$$\vartheta_C(z)=\sum_{x \in \Gamma_C} e^{2 \pi i z \text{Tr}_{k/{\Bbb Q}}
  \left(\frac{x \overline{x}}{p}\right)},\quad \text{where $z \in {\Bbb H}$}.$$
This is a modular form in one variable $z \in {\Bbb H}$. Now
$\Gamma_C$ is not only a ${\Bbb Z}$-module, but also an ${\frak
O}_k$-module. As above we can define a theta function in several
variables. For $z \in {\Bbb H}^{(p-1)/2}$ define
$$\theta_C(z):=\sum_{x \in \Gamma_C} e^{2 \pi i \text{Tr}_{k/{\Bbb Q}}
  \left(z \frac{x \overline{x}}{p}\right)}.$$

{\bf Theorem 4.3.} (see \cite{E}). {\it The function $\theta_C$ is
a Hilbert modular form of weight $n$ for the whole group
$SL_2({\frak O}_k)$.}

  Note that the Hilbert modular forms $\theta_j$ and $\theta_C$ are
symmetric (in the sense of the definition above). This is due to
the fact that the lattices ${\frak P}$ and $\Gamma_C$ are
invariant under the obvious action of the Galois group of ${\Bbb
Q}(\zeta)$ over ${\Bbb Q}$.

  The Lee weight enumerator of a code $C \subset {\Bbb F}_p^n$ is
the polynomial
$$W_C \left(X_0, X_1, \cdots, X_{\frac{p-1}{2}}\right)
 :=\sum_{u \in C} X_0^{l_0(u)} X_1^{l_1(u)} \cdots X_{\frac{p-1}{2}}^{l_{(p-1)/2}(u)},$$
where $l_0(u)$ is the number of zeros in $u$, and $l_i(u)$, for
$i=1, \cdots, \frac{p-1}{2}$, is the number of $+i$ or $-i$
occurring in the codeword $u$. This is a homogeneous polynomial of
degree $n$.

  We can now formulate the main theorem of G. van der Geer and F.
Hirzebruch.

{\bf Theorem 4.4. (van der Geer and Hirzebruch).} (see \cite{E}).
{\it Let $C \subset {\Bbb F}_p^n$ be a code with $C \subset
C^{\bot}$. Then the following identity holds:
$$\theta_C=W_C \left(\theta_0, \theta_1, \cdots, \theta_{\frac{p-1}{2}}\right).$$}

  The polynomial $W_C \left(\theta_0, \theta_1, \cdots, \theta_{\frac{p-1}{2}}\right)$
is an invariant polynomial for the above mentioned representation
of dimension $\frac{p+1}{2}$ of the group $SL(2, {\Bbb F}_p)$.

  For $p=3$, $k={\Bbb Q}$. The well-known result of Brou\'{e} and
Enguehard drops out. Namely, the Hamming weight enumerator
$H_C(\theta_0, \theta_1)$ is a polynomial in the modular forms
$E_4$, $E_6^2$ of $SL(2, {\Bbb Z})$ where
$$E_4=\theta_0^2+8 \theta_0 \theta_1^3, \quad
  E_6=\theta_0^6-20 \theta_0^3 \theta_1^3-8 \theta_1^6.$$

  For $p=5$, $k={\Bbb Q}(\sqrt{5})$. In his paper \cite{Hi1},
Hirzebruch proved that the ring of symmetric Hilbert modular forms
for $SL_2({\frak O}_k)(\sqrt{5})$ equals ${\Bbb C}[A_0, A_1, A_2]$
where the Klein invariants $A_0$, $A_1$, $A_2$ (see \cite{K}) have
weight $1$. He proved that the ring of symmetric Hilbert modular
forms for $SL_2({\frak O}_k)$ of even weight equals ${\Bbb C}[A,
B, C]$ where $A$, $B$, $C$ are the Klein invariants of degrees
$2$, $6$, $10$ (see \cite{K}). One has $A_0=\theta_0$, $A_1=2
\theta_1$, $A_2=2 \theta_2$. Using the $3$-dimensional
representation of $PSL(2, {\Bbb F}_5) \cong \langle \widetilde{S},
\widetilde{T} \rangle$ constructed by Klein (see \cite{K}):
$$\widetilde{S}=-\frac{1}{\sqrt{5}} \left(\matrix
  1 & 2 & 2\\
  1 & \varepsilon+\varepsilon^4 & \varepsilon^2+\varepsilon^3\\
  1 & \varepsilon^2+\varepsilon^3 & \varepsilon+\varepsilon^4
\endmatrix\right), \quad
  \widetilde{T}=\left(\matrix
  1 & 0 & 0\\
  0 & \varepsilon^2 & 0\\
  0 & 0 & \varepsilon^3
\endmatrix\right),$$
where $\varepsilon=e^{2 \pi i/5}$, one can obtain the MacWilliams
identity for Lee weight enumerators of codes over ${\Bbb F}_5$:
$$W_{C^{\bot}}(X_0, X_1, X_2)=\text{const} \cdot
  W_C \left(\matrix X_0+2 X_1+2 X_2\\
  X_0+(\varepsilon+\varepsilon^4) X_1+(\varepsilon^2+\varepsilon^3) X_2\\
  X_0+(\varepsilon^2+\varepsilon^3) X_1+(\varepsilon+\varepsilon^4) X_2
\endmatrix\right).$$

  For $p=7$, the corresponding Hilbert modular variety of dimension $3$
was investigated by E. Thomas. The invariant theory for the above
mentioned $4$-dimensional representation of $SL(2, {\Bbb F}_7)$
enters (see \cite{MS}).

  In our case, $p=13$. Up to the constants, our invariants ${\Bbb A}_0$,
$\cdots$, ${\Bbb A}_6$ are just $\theta_0$, $\cdots$, $\theta_6$.
Using our $7$-dimensional representation of $PSL(2, {\Bbb F}_{13})
\cong \langle \widetilde{S}, \widetilde{T} \rangle$ constructed as
above, we get the MacWilliams identity for Lee weight enumerators
of codes over ${\Bbb F}_{13}$:
$$W_{C^{\bot}}(X_0, X_1, X_2, X_4, X_5, X_6)
=\text{const} \cdot W_C \left( (X_0, X_1, X_2, X_4, X_5, X_6)
 \widetilde{S}\right).\eqno{(4.26)}$$

  In \cite{K1}, Klein obtained the modular equation of degree fourteen,
which corresponds to the transformation of order thirteen:
$$\aligned
J: J-1: 1=&(\tau^2+5 \tau+13)(\tau^4+7 \tau^3+20 \tau^2+19 \tau+1)^3\\
         :&(\tau^2+6 \tau+13)(\tau^6+10 \tau^5+46 \tau^4+108 \tau^3+122 \tau^2+38 \tau-1)^2\\
         :&1728 \tau,
\endaligned$$
Note that the Hauptmodul $J$ can be defined over the real
quadratic field ${\Bbb Q}(\sqrt{13})$:
$$\aligned
 &\tau^4+7 \tau^3+20 \tau^2+19 \tau+1\\
=&\left(\tau^2+\frac{7+\sqrt{13}}{2} \tau+\frac{11+3
  \sqrt{13}}{2}\right)\left(\tau^2+\frac{7-\sqrt{13}}{2} \tau+\frac{11-3 \sqrt{13}}{2}\right),
\endaligned$$
$$\aligned
 &\tau^6+10 \tau^5+46 \tau^4+108 \tau^3+122 \tau^2+38 \tau-1\\
=&\left(\tau^3+5 \tau^2+\frac{21-\sqrt{13}}{2}
  \tau+\frac{3+\sqrt{13}}{2}\right) \left(\tau^3+5 \tau^2+\frac{21+\sqrt{13}}{2} \tau
  +\frac{3-\sqrt{13}}{2}\right).
\endaligned$$
Let us confine our thought to an especially important result which
Jacobi had established as early as $1829$ in his ``Notices sur les
fonctions elliptiques'' (see \cite{K}). Jacobi there considered,
instead of the modular equation, the so-called
multiplier-equation, together with other equations equivalent to
it, and found that their $(n+1)$ roots are composed in a simple
manner of $\frac{n+1}{2}$ elements, with the help of merely
numerical irrationalities. Namely, if we denote these elements by
${\Bbb A}_0$, ${\Bbb A}_1$, $\cdots$, ${\Bbb A}_{\frac{n-1}{2}}$,
and further, for the roots $z$ of the equation under
consideration, apply the indices employed by Galois, we have, with
appropriate determination of the square root occurring on the
left-hand side:
$$\left\{\aligned
  \sqrt{z_{\infty}} &=\sqrt{(-1)^{\frac{n-1}{2}} \cdot n} \cdot {\Bbb A}_0,\\
  \sqrt{z_{\nu}} &={\Bbb A}_0+\epsilon^{\nu} {\Bbb A}_1+\epsilon^{4 \nu} {\Bbb A}_2
  +\cdots+\epsilon^{(\frac{n-1}{2})^2 \nu} {\Bbb A}_{\frac{n-1}{2}}
\endaligned\right.$$
for $\nu=0, 1, \cdots, n-1$ and $\epsilon=e^{\frac{2 \pi i}{n}}$.
Jacobi had himself emphasized the special significance of his
result by adding to his short communication: ``C'est un
th\'{e}or\`{e}me des plus importants dans la th\'{e}orie
alg\'{e}brique de la transformation et de la division des
fonctions elliptiques.'' Now, we give the Jacobian equation of
degree fourteen which corresponds to the above modular equation of
degree fourteen.

  Let $H=y_2 \cdot S$. Then
$$H=Q^5 P^2 \cdot P^2 Q^6 P^8 \cdot Q^5 P^2 \cdot P^3 Q
 =\left(\matrix
  0 &  0 &  0 & 0 & 0 & 1\\
  0 &  0 &  0 & 1 & 0 & 0\\
  0 &  0 &  0 & 0 & 1 & 0\\
  0 &  0 & -1 & 0 & 0 & 0\\
 -1 &  0 &  0 & 0 & 0 & 0\\
  0 & -1 &  0 & 0 & 0 & 0
\endmatrix\right).\eqno{(4.27)}$$
Note that
$$H^2=\left(\matrix
  0 & -1 &  0 &  0 &  0 &  0\\
  0 &  0 & -1 &  0 &  0 &  0\\
 -1 &  0 &  0 &  0 &  0 &  0\\
  0 &  0 &  0 &  0 & -1 &  0\\
  0 &  0 &  0 &  0 &  0 & -1\\
  0 &  0 &  0 & -1 &  0 &  0
\endmatrix\right), \quad
  H^3=\left(\matrix
  0 & 0 & 0 & -1 &  0 &  0\\
  0 & 0 & 0 &  0 & -1 &  0\\
  0 & 0 & 0 &  0 &  0 & -1\\
  1 & 0 & 0 &  0 &  0 &  0\\
  0 & 1 & 0 &  0 &  0 &  0\\
  0 & 0 & 1 &  0 &  0 &  0
\endmatrix\right)$$
and $H^6=-I$. In the projective coordinates, this means that
$H^6=1$. We have
$$H^{-1} T H=-T^4.$$
Thus, $\langle H, T \rangle \cong {\Bbb Z}_{13} \rtimes {\Bbb
Z}_6$. Hence, it is the maximal subgroup of order $78$ of $G$ with
index $14$ (see \cite{CC}). We find that $\varphi_{\infty}^2$ is
invariant under the action of the maximal subgroup $\langle H, T
\rangle$. Note that
$$\varphi_{\infty}=\sqrt{13} {\Bbb A}_0, \quad
  \varphi_{\nu}={\Bbb A}_0+\zeta^{\nu} {\Bbb A}_1+\zeta^{4 \nu} {\Bbb A}_2+
  \zeta^{9 \nu} {\Bbb A}_3+\zeta^{3 \nu} {\Bbb A}_4+\zeta^{12 \nu} {\Bbb A}_5+\zeta^{10 \nu} {\Bbb A}_6$$
for $\nu=0, 1, \cdots, 12$. Let $w=\varphi^2$,
$w_{\infty}=\varphi_{\infty}^2$ and $w_{\nu}=\varphi_{\nu}^2$.
Then $w_{\infty}$, $w_{\nu}$ for $\nu=0, \cdots, 12$ form an
algebraic equation of degree fourteen, which is just the Jacobian
equation of degree fourteen, whose roots are these $w_{\nu}$ and
$w_{\infty}$:
$$w^{14}+a_1 w^{13}+\cdots+ a_{13} w+a_{14}=0.$$
In particular, the coefficients
$$a_{14}=\varphi_{\infty}^2 \cdot \prod_{\nu=0}^{12} \varphi_{\nu}^2
 =13 {\Bbb A}_0^2 \prod_{\nu=0}^{12} ({\Bbb A}_0+\zeta^{\nu} {\Bbb A}_1+\zeta^{4 \nu} {\Bbb A}_2+
  \zeta^{9 \nu} {\Bbb A}_3+\zeta^{3 \nu} {\Bbb A}_4+\zeta^{12 \nu} {\Bbb A}_5+\zeta^{10 \nu} {\Bbb A}_6)^2,$$
and
$$-a_1=w_{\infty}+\sum_{\nu=0}^{12} w_{\nu}=
  26 ({\Bbb A}_0^2+{\Bbb A}_1 {\Bbb A}_5+{\Bbb A}_2 {\Bbb A}_3+{\Bbb A}_4 {\Bbb A}_6).\eqno{(4.28)}$$

  Recall that the modular equation of degree fourteen is also given by (see \cite{Hu})
$$\aligned
  z^{14}
 &+13 [2 \Delta^{\frac{1}{12}} z^{13}+25 \Delta^{\frac{2}{12}} z^{12}
  +196 \Delta^{\frac{3}{12}} z^{11}+1064 \Delta^{\frac{4}{12}} z^{10}\\
 &+4180 \Delta^{\frac{5}{12}} z^{9}+12086 \Delta^{\frac{6}{12}} z^{8}
  +25660 \Delta^{\frac{7}{12}} z^{7}+39182 \Delta^{\frac{8}{12}} z^{6}\\
 &+41140 \Delta^{\frac{9}{12}} z^{5}+27272 \Delta^{\frac{10}{12}} z^{4}
  +9604 \Delta^{\frac{11}{12}} z^{3}+1165 \Delta z^{2}]\\
 &+[746 \Delta-(12 g_2)^3] \Delta^{\frac{1}{12}} z+13 \Delta^{\frac{14}{12}}=0.
\endaligned$$
The coefficients of $z^{13}$ and the constant are just $26$ and
$13$.

  In fact, the invariant quadric $L:={\Bbb A}_0^2+{\Bbb A}_1 {\Bbb
A}_5+{\Bbb A}_2 {\Bbb A}_3+{\Bbb A}_4 {\Bbb A}_6$ is equal to
$$2 \left[(z_3 z_4^3+z_1 z_5^3+z_2 z_6^3)-(z_6 z_1^3+z_4 z_2^3+z_5 z_3^3)+
  3(z_1 z_2 z_4 z_5+z_2 z_3 z_5 z_6+z_3 z_1 z_6 z_4) \right].\eqno{(4.29)}$$
Hence, the variety $L=0$ is a quartic four-fold, which is
invariant under the action of the simple group $G$.

  For $g \in G=\langle S, T \rangle$, set
$$g({\Bbb A}_0):={\Bbb A}_0(g(z_1, z_2, z_3, z_4, z_5, z_6)).$$
We find that
$$\left\{\aligned
  {\Bbb A}_0 &={\Bbb A}_0,\\
    \sqrt{13} Q({\Bbb A}_0) &={\Bbb A}_0+\zeta^3 {\Bbb A}_1+\zeta^{12} {\Bbb A}_2
  +\zeta {\Bbb A}_3+\zeta^9 {\Bbb A}_4+\zeta^{10} {\Bbb A}_5+\zeta^4 {\Bbb A}_6,\\
  \sqrt{13} Q^2({\Bbb A}_0) &={\Bbb A}_0+\zeta^7 {\Bbb A}_1+\zeta^2 {\Bbb A}_2
  +\zeta^{11} {\Bbb A}_3+\zeta^8 {\Bbb A}_4+\zeta^6 {\Bbb A}_5+\zeta^5 {\Bbb A}_6,\\
  \sqrt{13} Q^3({\Bbb A}_0) &=-{\Bbb A}_0-\zeta {\Bbb A}_1-\zeta^4 {\Bbb A}_2
  -\zeta^9 {\Bbb A}_3-\zeta^3 {\Bbb A}_4-\zeta^{12} {\Bbb A}_5-\zeta^{10} {\Bbb A}_6,\\
  \sqrt{13} Q^4({\Bbb A}_0) &=-{\Bbb A}_0-\zeta^2 {\Bbb A}_1-\zeta^8 {\Bbb A}_2
  -\zeta^5 {\Bbb A}_3-\zeta^6 {\Bbb A}_4-\zeta^{11} {\Bbb A}_5-\zeta^7 {\Bbb A}_6,\\
  \sqrt{13} Q^5({\Bbb A}_0) &={\Bbb A}_0+\zeta^9 {\Bbb A}_1+\zeta^{10} {\Bbb A}_2
  +\zeta^3 {\Bbb A}_3+\zeta {\Bbb A}_4+\zeta^4 {\Bbb A}_5+\zeta^{12} {\Bbb A}_6,\\
  \sqrt{13} Q^6({\Bbb A}_0) &={\Bbb A}_0+{\Bbb A}_1+{\Bbb A}_2+{\Bbb A}_3+{\Bbb A}_4
  +{\Bbb A}_5+{\Bbb A}_6.
\endaligned\right.\eqno{(4.30)}$$
For our third model, we have
$$y_3 Q=-\frac{1}{\sqrt{13}} \left(\matrix
  \zeta^3-\zeta^5 & \zeta^6-\zeta^{12} & \zeta^6-\zeta^{11} &
  \zeta^{11}-\zeta & \zeta^3-\zeta^{12} & \zeta-\zeta^2\\
  \zeta^2-\zeta^8 & \zeta-\zeta^6 & \zeta^2-\zeta^4 &
  \zeta^9-\zeta^5 & \zeta^8-\zeta^9 & \zeta-\zeta^4\\
  \zeta^5-\zeta^{10} & \zeta^5-\zeta^7 & \zeta^9-\zeta^2 &
  \zeta^9-\zeta^{10} & \zeta^3-\zeta^6 & \zeta^7-\zeta^3\\
  \zeta^{12}-\zeta^2 & \zeta-\zeta^{10} & \zeta^{11}-\zeta^{12} &
  \zeta^{10}-\zeta^8 & \zeta^7-\zeta & \zeta^7-\zeta^2\\
  \zeta^8-\zeta^4 & \zeta^4-\zeta^5 & \zeta^9-\zeta^{12} &
  \zeta^{11}-\zeta^5 & \zeta^{12}-\zeta^7 & \zeta^{11}-\zeta^9\\
  \zeta^3-\zeta^4 & \zeta^7-\zeta^{10} & \zeta^{10}-\zeta^6 &
  \zeta^8-\zeta^3 & \zeta^8-\zeta^6 & \zeta^4-\zeta^{11}
\endmatrix\right),$$
$$y_3 Q^2=-\frac{1}{\sqrt{13}} \left(\matrix
  1-\zeta & 1-\zeta^3 & \zeta^6-\zeta^2 &
  \zeta^6-\zeta & \zeta^{10}-\zeta^8 & \zeta^3-\zeta^{10}\\
  \zeta^2-\zeta^5 & 1-\zeta^9 & 1-\zeta &
  \zeta-\zeta^{12} & \zeta^2-\zeta^9 & \zeta^{12}-\zeta^7\\
  1-\zeta^9 & \zeta^5-\zeta^6 & 1-\zeta^3 &
  \zeta^4-\zeta^{11} & \zeta^9-\zeta^4 & \zeta^5-\zeta^3\\
  \zeta^{12}-\zeta^7 & \zeta^5-\zeta^3 & \zeta^3-\zeta^{10} &
  1-\zeta^{12} & 1-\zeta^{10} & \zeta^7-\zeta^{11}\\
  \zeta-\zeta^{12} & \zeta^4-\zeta^{11} & \zeta^6-\zeta &
  \zeta^{11}-\zeta^8 & 1-\zeta^4 & 1-\zeta^{12}\\
  \zeta^2-\zeta^9 & \zeta^9-\zeta^4 & \zeta^{10}-\zeta^8 &
  1-\zeta^4 & \zeta^8-\zeta^7 & 1-\zeta^{10}
\endmatrix\right).$$
$$y_3 Q^3=-\frac{1}{\sqrt{13}} \left(\matrix
  \zeta^{12}-\zeta^3 & \zeta^4-\zeta^3 & \zeta^8-\zeta^5 &
  \zeta^{12}-\zeta^5 & \zeta^{10}-\zeta^2 & \zeta^2-\zeta^4\\
  \zeta^7-\zeta^6 & \zeta^4-\zeta & \zeta^{10}-\zeta &
  \zeta^5-\zeta^{10} & \zeta^4-\zeta^6 & \zeta^{12}-\zeta^5\\
  \zeta^{12}-\zeta^9 & \zeta^{11}-\zeta^2 & \zeta^{10}-\zeta^9 &
  \zeta^4-\zeta^6 & \zeta^6-\zeta^{12} & \zeta^{10}-\zeta^2\\
  \zeta^8-\zeta & \zeta^{11}-\zeta^3 & \zeta^9-\zeta^{11} &
  \zeta-\zeta^{10} & \zeta^9-\zeta^{10} & \zeta^5-\zeta^8\\
  \zeta^3-\zeta^8 & \zeta^7-\zeta^9 & \zeta^8-\zeta &
  \zeta^6-\zeta^7 & \zeta^9-\zeta^{12} & \zeta^3-\zeta^{12}\\
  \zeta^7-\zeta^9 & \zeta-\zeta^7 & \zeta^{11}-\zeta^3 &
  \zeta-\zeta^4 & \zeta^2-\zeta^{11} & \zeta^3-\zeta^4
\endmatrix\right).$$
$$y_3 Q^4=-\frac{1}{\sqrt{13}} \left(\matrix
  \zeta-\zeta^{10} & \zeta^8-\zeta^9 & \zeta^8-\zeta^{11} &
  \zeta^{11}-\zeta^5 & \zeta-\zeta^9 & \zeta^5-\zeta^3\\
  \zeta^7-\zeta^8 & \zeta^9-\zeta^{12} & \zeta^7-\zeta^3 &
  \zeta^6-\zeta & \zeta^8-\zeta^6 & \zeta^9-\zeta^3\\
  \zeta^{11}-\zeta & \zeta^{11}-\zeta^7 & \zeta^3-\zeta^4 &
  \zeta^3-\zeta & \zeta^2-\zeta^9 & \zeta^7-\zeta^2\\
  \zeta^8-\zeta^2 & \zeta^4-\zeta^{12} & \zeta^{10}-\zeta^8 &
  \zeta^{12}-\zeta^3 & \zeta^5-\zeta^4 & \zeta^5-\zeta^2\\
  \zeta^{12}-\zeta^7 & \zeta^7-\zeta^5 & \zeta^{10}-\zeta^4 &
  \zeta^6-\zeta^5 & \zeta^4-\zeta & \zeta^6-\zeta^{10}\\
  \zeta^{12}-\zeta^{10} & \zeta^4-\zeta^{11} & \zeta^{11}-\zeta^6 &
  \zeta^2-\zeta^{12} & \zeta^2-\zeta^6 & \zeta^{10}-\zeta^9
\endmatrix\right).$$
$$y_3 Q^5=-\frac{1}{\sqrt{13}} \left(\matrix
  1-\zeta^{12} & \zeta^{12}-\zeta^9 & \zeta^{10}-\zeta &
  \zeta^4-\zeta^9 & \zeta-\zeta^3 & \zeta^4-\zeta^{10}\\
  \zeta^{12}-\zeta^9 & 1-\zeta^4 & \zeta^4-\zeta^3 &
  \zeta^{10}-\zeta^{12} & \zeta^{10}-\zeta^3 & \zeta^9-\zeta\\
  \zeta^{10}-\zeta & \zeta^4-\zeta^3 & 1-\zeta^{10} &
  \zeta^3-\zeta^9 & \zeta^{12}-\zeta^4 & \zeta^{12}-\zeta\\
  \zeta^4-\zeta^9 & \zeta^{10}-\zeta^{12} & \zeta^3-\zeta^9 &
  1-\zeta & \zeta-\zeta^4 & \zeta^3-\zeta^{12}\\
  \zeta-\zeta^3 & \zeta^{10}-\zeta^3 & \zeta^{12}-\zeta^4 &
  \zeta-\zeta^4 & 1-\zeta^9 & \zeta^9-\zeta^{10}\\
  \zeta^4-\zeta^{10} & \zeta^9-\zeta & \zeta^{12}-\zeta &
  \zeta^3-\zeta^{12} & \zeta^9-\zeta^{10} & 1-\zeta^3
\endmatrix\right).$$
$$y_3 Q^6=-\frac{1}{\sqrt{13}} \left(\matrix
  \zeta^{10}-\zeta^8 & \zeta^6-1 & \zeta^{10}-\zeta^5 &
  \zeta^{12}-\zeta^9 & \zeta^8-\zeta^{12} & \zeta^6-\zeta^5\\
  \zeta^{12}-\zeta^6 & \zeta^{12}-\zeta^7 & \zeta^2-1 &
  \zeta^2-\zeta^6 & \zeta^4-\zeta^3 & \zeta^7-\zeta^4\\
  \zeta^5-1 & \zeta^4-\zeta^2 & \zeta^4-\zeta^{11} &
  \zeta^{11}-\zeta^{10} & \zeta^5-\zeta^2 & \zeta^{10}-\zeta\\
  \zeta^4-\zeta & \zeta-\zeta^5 & \zeta^8-\zeta^7 &
  \zeta^3-\zeta^5 & \zeta^7-1 & \zeta^3-\zeta^8\\
  \zeta^7-\zeta^{11} & \zeta^{10}-\zeta^9 & \zeta^9-\zeta^6 &
  \zeta-\zeta^7 & \zeta-\zeta^6 & \zeta^{11}-1\\
  \zeta^3-\zeta^2 & \zeta^{11}-\zeta^8 & \zeta^{12}-\zeta^3 &
  \zeta^8-1 & \zeta^9-\zeta^{11} & \zeta^9-\zeta^2
\endmatrix\right).$$
Hence,
$$\left\{\aligned
  \sqrt{13} y_3({\Bbb A}_0) &=-{\Bbb A}_0-\zeta^8 {\Bbb A}_1-\zeta^6 {\Bbb A}_2-\zeta^7 {\Bbb A}_3
  -\zeta^{11} {\Bbb A}_4-\zeta^5 {\Bbb A}_5-\zeta^2 {\Bbb A}_6,\\
  \sqrt{13} y_3 Q({\Bbb A}_0) &={\Bbb A}_0+\zeta^{11} {\Bbb A}_1+\zeta^5 {\Bbb A}_2+\zeta^8 {\Bbb A}_3
  +\zeta^7 {\Bbb A}_4+\zeta^2 {\Bbb A}_5+\zeta^6 {\Bbb A}_6,\\
  \sqrt{13} y_3 Q^2({\Bbb A}_0) &=-{\Bbb A}_0-\zeta^{10} {\Bbb A}_1-\zeta {\Bbb A}_2-\zeta^{12} {\Bbb A}_3
  -\zeta^4 {\Bbb A}_4-\zeta^3 {\Bbb A}_5-\zeta^9 {\Bbb A}_6,\\
  \sqrt{13} y_3 Q^3({\Bbb A}_0) &=-{\Bbb A}_0-\zeta^{12} {\Bbb A}_1-\zeta^9 {\Bbb A}_2-\zeta^4 {\Bbb A}_3
  -\zeta^{10} {\Bbb A}_4-\zeta {\Bbb A}_5-\zeta^3 {\Bbb A}_6,\\
  \sqrt{13} y_3 Q^4({\Bbb A}_0) &=-{\Bbb A}_0-\zeta^4 {\Bbb A}_1-\zeta^3 {\Bbb A}_2-\zeta^{10} {\Bbb A}_3
  -\zeta^{12} {\Bbb A}_4-\zeta^9 {\Bbb A}_5-\zeta {\Bbb A}_6,\\
  \sqrt{13} y_3 Q^5({\Bbb A}_0) &=-{\Bbb A}_0-\zeta^6 {\Bbb A}_1-\zeta^{11} {\Bbb A}_2-\zeta^2 {\Bbb A}_3
  -\zeta^5 {\Bbb A}_4-\zeta^7 {\Bbb A}_5-\zeta^8 {\Bbb A}_6,\\
  \sqrt{13} y_3 Q^6({\Bbb A}_0) &={\Bbb A}_0+\zeta^5 {\Bbb A}_1+\zeta^7 {\Bbb A}_2+\zeta^6 {\Bbb A}_3
  +\zeta^2 {\Bbb A}_4+\zeta^8 {\Bbb A}_5+\zeta^{11} {\Bbb A}_6.
\endaligned\right.\eqno{(4.31)}$$
Moreover,
$$\left\{\aligned
  {\Bbb A}_0^2 &={\Bbb A}_0^2,\\
  Q({\Bbb A}_0)^2 &=ST^3({\Bbb A}_0)^2,\\
  Q^2({\Bbb A}_0)^2 &=ST^7({\Bbb A}_0)^2,\\
  Q^3({\Bbb A}_0)^2 &=ST({\Bbb A}_0)^2,\\
  Q^4({\Bbb A}_0)^2 &=ST^2({\Bbb A}_0)^2,\\
  Q^5({\Bbb A}_0)^2 &=ST^9({\Bbb A}_0)^2,\\
  Q^6({\Bbb A}_0)^2 &=S({\Bbb A}_0)^2.
\endaligned\right. \quad \quad
\left\{\aligned
      y_3({\Bbb A}_0)^2 &=ST^8({\Bbb A}_0)^2,\\
    y_3 Q({\Bbb A}_0)^2 &=ST^{11}({\Bbb A}_0)^2,\\
  y_3 Q^2({\Bbb A}_0)^2 &=ST^{10}({\Bbb A}_0)^2,\\
  y_3 Q^3({\Bbb A}_0)^2 &=ST^{12}({\Bbb A}_0)^2,\\
  y_3 Q^4({\Bbb A}_0)^2 &=ST^4({\Bbb A}_0)^2,\\
  y_3 Q^5({\Bbb A}_0)^2 &=ST^6({\Bbb A}_0)^2,\\
  y_3 Q^6({\Bbb A}_0)^2 &=ST^5({\Bbb A}_0)^2.
\endaligned\right.\eqno{(4.32)}$$

  On the other hand,
$$\sqrt{13} ST^{\nu}({\Bbb A}_0)={\Bbb A}_0+\zeta^{\nu} {\Bbb A}_1+\zeta^{4 \nu} {\Bbb A}_2+
  \zeta^{9 \nu} {\Bbb A}_3+\zeta^{3 \nu} {\Bbb A}_4+\zeta^{12 \nu} {\Bbb A}_5+\zeta^{10 \nu} {\Bbb A}_6, \quad
  T^{\nu}({\Bbb A}_0)={\Bbb A}_0$$
for $\nu=0, 1, \cdots, 12$. Comparing two sets $\{ Q^{\nu}({\Bbb
A}_0)^2, y_3 Q^{\nu}({\Bbb A}_0)^2 \}$ for $\nu$ mod $7$ and $\{
T^{\nu}({\Bbb A}_0)^2, S T^{\nu}({\Bbb A}_0)^2\}$ for $\nu$ mod
$13$ leads to the following:

{\bf Theorem 4.5. (Exotic duality theorem)}. {\it There exists an
``exotic'' duality between non-congruence modular forms and
Hilbert modular forms, both of them are related to ${\Bbb
Q}(\zeta):$
$$\left\{\matrix
  Q^{\nu}({\Bbb A}_0)^2, y_3 Q^{\nu}({\Bbb A}_0)^2\\
  \text{$\nu$ mod $7$}
  \endmatrix\right\}=
  \left\{\matrix
  T^{\nu}({\Bbb A}_0)^2, S T^{\nu}({\Bbb A}_0)^2\\
  \text{$\nu$ mod $13$}
\endmatrix\right\}\eqno{(4.33)}$$}

  The left hand side of (4.33) corresponds to $(2, 3, 7; 13)$,
the right hand side of (4.33) corresponds to $(2, 3, 13; 7)$. When
we exchange the position of $7$ and $13$, we will get the duality
between the two representations of $G$. Note that ${\Bbb A}_0^2$
is invariant under the action of the maximal subgroup $\langle H,
T \rangle$ of order $78$ with index $14$ of $G$. Hence,
$T^{\nu}({\Bbb A}_0)^2, S T^{\nu}({\Bbb A}_0)^2$ for $\nu$ mod
$13$ give a coset decomposition of $G$ with respect to this
maximal subgroup. Therefore, any symmetric polynomial of ${\Bbb
A}_0^2$, $ST^{\nu}({\Bbb A}_0)^2$ for $\nu$ mod $13$ gives an
invariant polynomial under the action of $PSL(2, 13)$. By (4.33),
we have that any symmetric polynomial of $Q^{\nu}({\Bbb A}_0)^2,
y_3 Q^{\nu}({\Bbb A}_0)^2$ for $\nu$ mod $7$ also gives an
invariant polynomial under the action of $PSL(2, 13)$. It is known
that invariant polynomials correspond to modular forms. The right
hand side of (4.33) corresponds to the congruence subgroup
$\Gamma_{13, 1}$ with type $(1, 13)$, the associated modular forms
are Hilbert modular forms which we study as above. However, the
left hand side of (4.33) corresponds to the non-congruence
subgroup $\Gamma_{7, 7}$ with type $(7, 7)$, the associated
modular forms are non-congruence modular forms! Therefore, (4.33)
gives a connection between congruence modular forms (Hilbert
modular forms) and non-congruence modular forms, both of them are
related to the cyclotomic field ${\Bbb Q}(\zeta)$! Thus, we
complete the proof of Theorem 1.1.

\flushpar $\qquad \qquad \qquad \qquad \qquad \qquad \qquad \qquad
\qquad \qquad \qquad \qquad \qquad \qquad \qquad \qquad \qquad
\qquad \quad \boxed{}$

  Now, we find that the third non-congruence modular curve $Y_3$
possesses exotic duality between the associated non-congruence
modular forms and the Hilbert modular forms. It is natural to ask
that what will happen for the other non-congruence modular curves
$Y_1$ and $Y_2$? For our first model which corresponds to $(2, 3,
7; 7)$, we have
$$\left\{\aligned
  y_1({\Bbb A}_0) &=-{\Bbb A}_0,\\
  \sqrt{13} y_1 Q({\Bbb A}_0) &=-{\Bbb A}_0-\zeta^3 {\Bbb A}_1-\zeta^{12} {\Bbb A}_2
  -\zeta {\Bbb A}_3-\zeta^9 {\Bbb A}_4-\zeta^{10} {\Bbb A}_5-\zeta^4 {\Bbb A}_6,\\
  \sqrt{13} y_1 Q^2({\Bbb A}_0) &=-{\Bbb A}_0-\zeta^7 {\Bbb A}_1-\zeta^2 {\Bbb A}_2
  -\zeta^{11} {\Bbb A}_3-\zeta^8 {\Bbb A}_4-\zeta^6 {\Bbb A}_5-\zeta^5 {\Bbb A}_6,\\
  \sqrt{13} y_1 Q^3({\Bbb A}_0) &={\Bbb A}_0+\zeta {\Bbb A}_1+\zeta^4 {\Bbb A}_2
  +\zeta^9 {\Bbb A}_3+\zeta^3 {\Bbb A}_4+\zeta^{12} {\Bbb A}_5+\zeta^{10} {\Bbb A}_6,\\
  \sqrt{13} y_1 Q^4({\Bbb A}_0) &={\Bbb A}_0+\zeta^2 {\Bbb A}_1+\zeta^8 {\Bbb A}_2
  +\zeta^5 {\Bbb A}_3+\zeta^6 {\Bbb A}_4+\zeta^{11} {\Bbb A}_5+\zeta^7 {\Bbb A}_6,\\
  \sqrt{13} y_1 Q^5({\Bbb A}_0) &=-{\Bbb A}_0-\zeta^9 {\Bbb A}_1-\zeta^{10} {\Bbb A}_2
  -\zeta^3 {\Bbb A}_3-\zeta {\Bbb A}_4-\zeta^4 {\Bbb A}_5-\zeta^{12} {\Bbb A}_6,\\
  \sqrt{13} y_1 Q^6({\Bbb A}_0) &=-{\Bbb A}_0-{\Bbb A}_1-{\Bbb A}_2-{\Bbb A}_3-{\Bbb A}_4
  -{\Bbb A}_5-{\Bbb A}_6.
\endaligned\right.\eqno{(4.34)}$$
For our second model which corresponds to $(2, 3, 7; 6)$, we have
$$\left\{\aligned
  \sqrt{13} y_2({\Bbb A}_0) &=-{\Bbb A}_0-{\Bbb A}_1-{\Bbb A}_2-{\Bbb A}_3-{\Bbb A}_4-{\Bbb A}_5-{\Bbb A}_6,\\
  y_2 Q({\Bbb A}_0) &=-{\Bbb A}_0,\\
  \sqrt{13} y_2 Q^2({\Bbb A}_0) &=-{\Bbb A}_0-\zeta^3 {\Bbb A}_1-\zeta^{12} {\Bbb A}_2
  -\zeta {\Bbb A}_3-\zeta^9 {\Bbb A}_4-\zeta^{10} {\Bbb A}_5-\zeta^4 {\Bbb A}_6,\\
  \sqrt{13} y_2 Q^3({\Bbb A}_0) &=-{\Bbb A}_0-\zeta^7 {\Bbb A}_1-\zeta^2 {\Bbb A}_2
  -\zeta^{11} {\Bbb A}_3-\zeta^8 {\Bbb A}_4-\zeta^6 {\Bbb A}_5-\zeta^5 {\Bbb A}_6,\\
  \sqrt{13} y_2 Q^4({\Bbb A}_0) &={\Bbb A}_0+\zeta {\Bbb A}_1+\zeta^4 {\Bbb A}_2
  +\zeta^9 {\Bbb A}_3+\zeta^3 {\Bbb A}_4+\zeta^{12} {\Bbb A}_5+\zeta^{10} {\Bbb A}_6,\\
  \sqrt{13} y_2 Q^5({\Bbb A}_0) &={\Bbb A}_0+\zeta^2 {\Bbb A}_1+\zeta^8 {\Bbb A}_2
  +\zeta^5 {\Bbb A}_3+\zeta^6 {\Bbb A}_4+\zeta^{11} {\Bbb A}_5+\zeta^7 {\Bbb A}_6,\\
  \sqrt{13} y_2 Q^6({\Bbb A}_0) &=-{\Bbb A}_0-\zeta^9 {\Bbb A}_1-\zeta^{10} {\Bbb A}_2
  -\zeta^3 {\Bbb A}_3-\zeta {\Bbb A}_4-\zeta^4 {\Bbb A}_5-\zeta^{12} {\Bbb A}_6.
\endaligned\right.\eqno{(4.35)}$$
Thus,
$$\left\{\aligned
      y_1({\Bbb A}_0)^2 &={\Bbb A}_0^2,\\
    y_1 Q({\Bbb A}_0)^2 &=Q({\Bbb A}_0)^2,\\
  y_1 Q^2({\Bbb A}_0)^2 &=Q^2({\Bbb A}_0)^2,\\
  y_1 Q^3({\Bbb A}_0)^2 &=Q^3({\Bbb A}_0)^2,\\
  y_1 Q^4({\Bbb A}_0)^2 &=Q^4({\Bbb A}_0)^2,\\
  y_1 Q^5({\Bbb A}_0)^2 &=Q^5({\Bbb A}_0)^2,\\
  y_1 Q^6({\Bbb A}_0)^2 &=Q^6({\Bbb A}_0)^2.
\endaligned\right. \quad \quad
  \left\{\aligned
      y_2({\Bbb A}_0)^2 &=Q^6({\Bbb A}_0)^2,\\
    y_2 Q({\Bbb A}_0)^2 &={\Bbb A}_0^2,\\
  y_2 Q^2({\Bbb A}_0)^2 &=Q({\Bbb A}_0)^2,\\
  y_2 Q^3({\Bbb A}_0)^2 &=Q^2({\Bbb A}_0)^2,\\
  y_2 Q^4({\Bbb A}_0)^2 &=Q^3({\Bbb A}_0)^2,\\
  y_2 Q^5({\Bbb A}_0)^2 &=Q^4({\Bbb A}_0)^2,\\
  y_2 Q^6({\Bbb A}_0)^2 &=Q^5({\Bbb A}_0)^2.
\endaligned\right.\eqno{(4.36)}$$
We find that both $Q^{\nu}({\Bbb A}_0)^2, y_1 Q^{\nu}({\Bbb
A}_0)^2$ and $Q^{\nu}({\Bbb A}_0)^2, y_2 Q^{\nu}({\Bbb A}_0)^2$
for $\nu$ mod $7$ only give half of the invariants $T^{\nu}({\Bbb
A}_0)^2, S T^{\nu}({\Bbb A}_0)^2$ for $\nu$ mod $13$.

  We have $PQP^2 \cdot Q \cdot PQP^2=Q^{-1}$ where the involution
$$PQP^2=-\frac{1}{\sqrt{13}} \left(\matrix
  \zeta^8-\zeta^5 & \zeta^{10}-\zeta & \zeta^{10}-\zeta^9 &
  \zeta^9-\zeta^{11} & \zeta^8-\zeta & \zeta^{11}-\zeta^3\\
  \zeta^{12}-\zeta^3 & \zeta^7-\zeta^6 & \zeta^{12}-\zeta^9 &
  \zeta^8-\zeta & \zeta^3-\zeta^8 & \zeta^7-\zeta^9\\
  \zeta^4-\zeta^3 & \zeta^4-\zeta & \zeta^{11}-\zeta^2 &
  \zeta^{11}-\zeta^3 & \zeta^7-\zeta^9 & \zeta-\zeta^7\\
  \zeta^2-\zeta^4 & \zeta^{12}-\zeta^5 & \zeta^{10}-\zeta^2 &
  \zeta^5-\zeta^8 & \zeta^3-\zeta^{12} & \zeta^3-\zeta^4\\
  \zeta^{12}-\zeta^5 & \zeta^5-\zeta^{10} & \zeta^4-\zeta^6 &
  \zeta-\zeta^{10} & \zeta^6-\zeta^7 & \zeta-\zeta^4\\
  \zeta^{10}-\zeta^2 & \zeta^4-\zeta^6 & \zeta^6-\zeta^{12} &
  \zeta^9-\zeta^{10} & \zeta^9-\zeta^{12} & \zeta^2-\zeta^{11}
\endmatrix\right).\eqno{(4.37)}$$
Hence, $\langle PQP^2, Q \rangle \cong {\Bbb Z}_7 \rtimes {\Bbb
Z}_2=D_{14}$, which is the maximal subgroup of order $14$ of $G$
with index $78$ (see \cite{CC}).

  Let $\rho: \Gamma \to  PGL(6, {\Bbb C})$ be a representation
given by $\rho(s)=S$ and $\rho(t)=T$. Then $\rho(h)=H$. Hence,
$\langle S, T, H \rangle$ gives the Weil representation for $SL(2,
13)$. On the other hand, we find that for the first model $\langle
y_1, Q^6, PQP^2 \rangle$:
$$\rho \left(\left(\matrix -5 & -3\\ 0 & 5 \endmatrix\right)\right)=y_1,
  \rho \left(\left(\matrix -3 & -1\\ 1 & 0 \endmatrix\right)\right)=Q^6,
  \rho \left(\left(\matrix -2 & -1\\ 5 & 2 \endmatrix\right)\right)=PQP^2.\eqno{(4.38)}$$
For the second model $\langle y_2, Q^5, PQP^2 \rangle$:
$$\rho \left(\left(\matrix 0 & -7\\ 2 & 0 \endmatrix\right)\right)=y_2,
  \rho \left(\left(\matrix 5 & 10\\ 3 & 1 \endmatrix\right)\right)=Q^5,
  \rho \left(\left(\matrix -2 & -1\\ 5 & 2 \endmatrix\right)\right)=PQP^2,\eqno{(4.39)}$$
For the third model $\langle y_3, Q^3, PQP^2 \rangle$:
$$\rho \left(\left(\matrix -1 & 10\\ 5 & 1 \endmatrix\right)\right)=y_3,
  \rho \left(\left(\matrix -3 & 5\\ 8 & 8 \endmatrix\right)\right)=Q^3,
  \rho \left(\left(\matrix -2 & -1\\ 5 & 2 \endmatrix\right)\right)=PQP^2.\eqno{(4.40)}$$
Thus, our three models give three kinds of representations which
are different from the Weil representation!

\vskip 0.5 cm

\centerline{\bf 5. Haagerup subfactor and exceptional Lie group
                   $G_2$}

\vskip 0.5 cm

  Recall that the subfactors of von Neumann algebras give $ADE$ type
classification if the Jones index $<4$. The Haagerup subfactor
(see \cite{AH}) is a finite-depth, irreducible, hyperfinite
subfactor with smallest index $\frac{5+\sqrt{13}}{2} \approx
4.3>4$.

  It is well-known that each $(2+1)$-dimensional topological quantum
field theory is conjectured to arise as a Chern-Simons-Witten
theory (see \cite{Wi}) from a compact Lie group with a certain
cohomology class. Such a Chern-Simons-Witten theory is supposed to
produce a modular tensor category. In \cite{HRW}, Hong, Rowell and
Wang studied two examples of a modular tensor category which do
not seem to arise from a Chern-Simons-Witten theory. Both examples
come from the subfactor theory of Jones. One is the quantum double
of the even part of the $E_6$ subfactor, and the other is also the
quantum double of the even part of the Haagerup subfactor. These
two examples are believed to provide counterexamples to the above
conjecture, but this is not proved here, because not all
Chern-Simons-Witten theories have been constructed in a
mathematically rigorous way. One of the significance of the
Haagerup subfactor comes from the following fact that the
$3$-sphere $S^3$ and the Poincar\'{e} homology $3$-sphere
$\Sigma(2, 3, 5)$ are distinguished by the Turaev-Viro-Ocneanu
invariant from the Haagerup subfactor (see \cite{SW}).

  On the other hand, there is a hierarchy of understanding: (1) conformal
field theory, (2) statistical mechanical models, (3) subfactors,
vertex operator algebras and twisted $K$-theory, (4) modular
tensor categories, pre-projective algebras, Calabi-Yau algebras
(see \cite{EG2}). The most basic algebraic structure here, namely
that of a modular tensor category, may arise from subfactors,
vertex operator algebras or twisted equivariant $K$-theory which
in turn may give rise to statistical mechanical models which at
criticality may produce conformal invariant field theories. That
is to say, two-dimensional conformal field theories can be
understood from the vantage point of conformal nets of subfactors
or vertex operator algebras. The most natural place to look for
exotic possibilities of subfactors and hence of conformal field
theories is with the Haagerup subfactor and its siblings. The
Haagerup modular data was computed by Izumi (see \cite{Iz}), with
$T$ being the diagonal matrix $\text{diag}(1, 1, 1, 1, \omega,
\overline{\omega}, \zeta^6, \zeta^{11}, \zeta^2, \zeta^5, \zeta^7,
\zeta^8)$, where $\omega=\exp(2 \pi i/3)$. His $S$ matrix was
complicated. In \cite{EG1}, Evans and Gannon derived an explicit
simple description for the $S$ matrix
$$S=\left(\matrix A & B\\ C & D \endmatrix\right)$$
with
$$A=\left(\matrix
  x & 1-x & 1 & 1 & 1 & 1\\
  1-x & x & 1 & 1 & 1 & 1\\
  1 &   1 & 2 & -1 & -1 & -1\\
  1 &   1 & -1 & 2 & -1 & -1\\
  1 &   1 & -1 & -1 & -1 & 2\\
  1 &   1 & -1 & -1 & 2 & -1
\endmatrix\right),
  B=\left(\matrix
  y &  y &  y &  y &  y &  y\\
  -y & -y & -y & -y & -y & -y\\
  0 & 0 & 0 & 0 & 0 & 0\\
  0 & 0 & 0 & 0 & 0 & 0\\
  0 & 0 & 0 & 0 & 0 & 0\\
  0 & 0 & 0 & 0 & 0 & 0
\endmatrix\right),$$
$$C=\left(\matrix
  y &  -y &  0 &  0 & 0 &  0\\
  y &  -y &  0 &  0 & 0 &  0\\
  y &  -y &  0 &  0 & 0 &  0\\
  y &  -y &  0 &  0 & 0 &  0\\
  y &  -y &  0 &  0 & 0 &  0\\
  y &  -y &  0 &  0 & 0 &  0
\endmatrix\right),
  D=\left(\matrix
  c(1) & c(2) & c(3) & c(4) & c(5) & c(6)\\
  c(2) & c(4) & c(6) & c(5) & c(3) & c(1)\\
  c(3) & c(6) & c(4) & c(1) & c(2) & c(5)\\
  c(4) & c(5) & c(1) & c(3) & c(6) & c(2)\\
  c(5) & c(3) & c(2) & c(6) & c(1) & c(4)\\
  c(6) & c(1) & c(5) & c(2) & c(4) & c(3)
\endmatrix\right)$$
for $x=(13-3 \sqrt{13})/26$, $y=3/\sqrt{13}$ and $c(j)=-2y \cos(2
\pi j/13)$. That this bears some relation with the double of $S_3$
may not be surprising given the relations between the Haagerup
fusion rules and those of $S_3$ and $\widehat{S}_3$. There is
however also a striking relationship with the affine algebra
modular data $B_{6, 2}$ which has central charge $c=12$ and $10$
primaries. The $T$-matrix is $\text{diag}(-1, -1; -1, -1;
-\zeta^6, -\zeta^{11}, -\zeta^2, -\zeta^5, -\zeta^7, -\zeta^8)$,
while the $S$-matrix is
$$S=\left(\matrix
  y/2 &  y/2 &  3/2 &  3/2 & y & y & y & y & y & y\\
  y/2 &  y/2 & -3/2 & -3/2 & y & y & y & y & y & y\\
  3/2 & -3/2 &  3/2 & -3/2 & 0 & 0 & 0 & 0 & 0 & 0\\
  3/2 & -3/2 & -3/2 &  3/2 & 0 & 0 & 0 & 0 & 0 & 0\\
    y & y & 0 & 0 & -c(1) & -c(2) & -c(3) & -c(4) & -c(5) & -c(6)\\
    y & y & 0 & 0 & -c(2) & -c(4) & -c(6) & -c(5) & -c(3) & -c(1)\\
    y & y & 0 & 0 & -c(3) & -c(6) & -c(4) & -c(1) & -c(2) & -c(5)\\
    y & y & 0 & 0 & -c(4) & -c(5) & -c(1) & -c(3) & -c(6) & -c(2)\\
    y & y & 0 & 0 & -c(5) & -c(3) & -c(2) & -c(6) & -c(1) & -c(4)\\
    y & y & 0 & 0 & -c(6) & -c(1) & -c(5) & -c(2) & -c(4) & -c(3)
\endmatrix\right),$$
where $y$ and $c(j)$ is as before. Ignoring the first $4$
primaries, the only difference with the Haagerup modular data are
some signs.

  According to \cite{J2} and \cite{J1}, these sporadic subfactors
discovered by Haagerup and others remain exotic creatures in the
zoo, untouched by quantum groups or any non-von Neumann algebras
approach. Hence, it is a major challenge in the theory to come up
with an interpretation of these subfactors as members of a family
related to some other mathematical objects.

  Now, we will give a connection between the Haagerup modular data
and our seven-dimensional representation of $PSL(2, 13)$. Note
that
$$p_1=\sqrt{13} c(2), p_2=\sqrt{13} c(4), p_3=\sqrt{13} c(6),
  p_4=\sqrt{13} c(5), p_5=\sqrt{13} c(3), p_6=\sqrt{13} c(1).$$
The actions of $T$ and $S$ on the basis $({\Bbb A}_0, {\Bbb A}_6,
{\Bbb A}_1, {\Bbb A}_5, {\Bbb A}_2, {\Bbb A}_4, {\Bbb A}_3)$ are
given as follows:
$$\widehat{T}=\text{diag}(1, \zeta^{10}, \zeta, \zeta^{12}, \zeta^4, \zeta^3, \zeta^9)$$ and
$$\widehat{S}=\frac{1}{\sqrt{13}} \left(\matrix
  1 &   1  &   1  &   1  &   1  &   1  &   1\\
  2 & c(1) & c(2) & c(3) & c(4) & c(5) & c(6)\\
  2 & c(2) & c(4) & c(6) & c(5) & c(3) & c(1)\\
  2 & c(3) & c(6) & c(4) & c(1) & c(2) & c(5)\\
  2 & c(4) & c(5) & c(1) & c(3) & c(6) & c(2)\\
  2 & c(5) & c(3) & c(2) & c(6) & c(1) & c(4)\\
  2 & c(6) & c(1) & c(5) & c(2) & c(4) & c(3)
 \endmatrix\right).\eqno{(5.1)}$$
Ignoring the first $1$ primary, this is just the Haagerup modular
data and the only difference with the affine algebra modular data
$B_{6, 2}$ are some signs!

  Cartan showed that there is a $14$-dimensional Lie group $G \subset SO(7, {\Bbb C})$
whose Lie algebra ${\frak g}$ has its root system of type $G_2$
(see \cite{B1}). It is interesting that, in his thesis, Cartan
does not describe this group as the subgroup preserving some
algebraic structure on ${\Bbb C}^7$. In fact, his original
description is in terms of the inner product and a set of
differential equations for curves in ${\Bbb C}^7$. What he says is
that $G$ is the subgroup of $GL(7, {\Bbb C})$ that preserves the
quadratic form
$$J=z^2+x_1 y_1+x_2 y_2+x_3 y_3\eqno{(5.2)}$$
and the system of $7$ Pfaffian equations (where $(i, j, k)$ is any
even permutation of $(1, 2, 3)$)
$$\aligned
  z dx_i-x_i dz+y_j dy_k-y_k dy_j &=0,\\
  z dy_i-y_i dz+x_j dx_k-x_k dx_j &=0,\\
  x_1 dy_1-y_1 dx_1+x_2 dy_2-y_2 dx_2+x_3 dy_3-y_3 dx_3 &=0.
\endaligned$$
Since $G$ preserves a quadratic form on ${\Bbb C}^7$, it can not
act transitively on ${\Bbb P}^6$, the space of lines in ${\Bbb
C}^7$. However, Cartan showed that $G$ does act transitively on
${\Bbb Q}_5 \subset {\Bbb P}^6$, the space of $J$-null lines in
${\Bbb C}^7$. Now $G$ does not act transitively on the space of
$J$-null $2$-planes in ${\Bbb C}^7$ (a $7$-dimensional homogeneous
space of $SO(7, {\Bbb C})$), but it does act transitively on the
$5$-dimensional space ${\Bbb N}_5$ consisting of the $J$-null
$2$-planes on which the following $2$-forms vanish (again, $(i, j,
k)$ is any even permutation of $(1, 2, 3)$):
$$\aligned
  dz \wedge dx_i+dy_j \wedge dy_k &=0,\\
  dz \wedge dy_i+dx_j \wedge dx_k &=0,\\
  dx_1 \wedge dy_1+dx_2 \wedge dy_2+dx_3 \wedge dy_3 &=0.
\endaligned$$
In fact, there exist exactly two maximal parabolic subgroups $P_1$
and $P_2$ in $G$. The quotient variety $G/P_1$ is isomorphic to a
$5$-dimensional quadric $Q^5 \subset {\Bbb P}^6$, and $G/P_2$ is a
$5$-dimensional Fano variety of index $3$.

  Generalized polygons were introduced by Tits in \cite{T} in
order to give a geometric interpretation for certain simple
algebraic groups. More generally, spherical buildings provide a
beatuiful and uniform geometric interpretation for all simple
algebraic groups, even the exceptional ones. Generalized polygons
are exactly the spherical buildings of rank two, which are also
related to the identification of simple groups of finite Morley
rank. Up to duality, there are precisely three generalized
polygons, namely the Moufang polygons associated to the simple Lie
groups $PGL_3$, $PSp_4$ and $G_2$. The corresponding geometric
objects: the Pappian projective planes, the symplectic quadrangles
and the split Cayley hexagons share a number of interesting
properties and characterizations. One of them is that they are
defined over any field, in particular over any algebraically
closed field. So in certain sense they are the prototypes of the
polygons, the universal examples.

  Let $K$ be any field and let $H(K)$ be the split Cayley hexagon
over $K$. We use the explicit description of $H(K)$ in
six-dimensional projective space (see \cite{Ma}). Let $H(K)$ be
embedded in the quadric $Q(6, K)$ with equation
$$X_0 X_4+X_1 X_5+X_2 X_6=X_3^2\eqno{(5.3)}$$
in $PG(6, K)$. The points of $H(K)$ are exactly the points of
$Q(6, K)$ and the lines of $H(K)$ are those lines on $Q(6, K)$
whose Grassmann coordinates satisfy:
$$p_{12}=p_{34}, p_{20}=p_{35}, p_{01}=p_{36}, p_{03}=p_{56},
  p_{13}=p_{64}, p_{23}=p_{45},$$
where $p_{ij}:=x_i y_j-x_j y_i$ for the line joining the points
with coordinates $x$ and $y$. This gives a complete and explicit
description for $H(K)$ on the quadric $Q(6, K)$. It is due to Tits
\cite{T}. The reason for that name is that this hexagon can also
be constructed using a split Cayley algebra over $K$ and moreover,
the corresponding simple algebraic group is also split. The split
Cayley hexagon is the most important hexagon, it is the main
example, and in fact, the only example for many fields $K$. The
full collineation group of $H(K)$ is isomorphic to the semi-direct
product $G_2(K): \text{Aut}(K)$. The split Cayley hexagons are
precisely the generalized hexagons arising from linear trialities
(i.e. where no field automorphism in the triality is involved).

  According to \cite{B2}, there are very beautiful conjectural
relations that the geometry of $G_2$ have ties in automorphic
forms and all kinds of special geometry that, if the physicists
were right, are just waiting there to be seen. But at the moment
our techniques are far too weak to be able to verify that claim.

  Now, we find that our invariant quadric
$${\Bbb A}_0^2+{\Bbb A}_1 {\Bbb A}_5+{\Bbb A}_2 {\Bbb A}_3+{\Bbb A}_4
  {\Bbb A}_6=0\eqno{(5.4)}$$
is just Cartan's ${\Bbb Q}_5 \subset {\Bbb P}^6$, the space of
$J$-null lines in ${\Bbb C}^7$, the quotient variety $G_2/P_1$
which is isomorphic to a $5$-dimensional quadric $Q^5 \subset
{\Bbb P}^6$ and the split Cayley hexagon over ${\Bbb C}$! This
shows that these three geometric objects are modular,
parameterized by both Hilbert modular forms and con-congruence
modular forms related to ${\Bbb Q}(\zeta)$. In fact, there are
five conjugacy classes of finite primitive subgroups of $G_2({\Bbb
C})$ (see \cite{CW}). Such a subgroup is isomorphic to one of
$PSL_2(13)$, $PSL_2(8)$, $PGL_2(7)$, $U_3(3)$ or $G_2(2)$.

\vskip 0.5 cm

\centerline{\bf 6. Fourteen-dimensional representations of $PSL(2,
                   13)$,}

\centerline{\bf exotic modular equation of degree fourteen and
                Monster simple group ${\Bbb M}$}

\vskip 0.5 cm

  Note that the Hurwitz curves of genus $14$ are are non-hyperelliptic.
Hence, we need to study their canonical models in ${\Bbb P}^{13}$.
There are $g(g+1)/2-(3g-3)=66$ quadratic homogeneous polynomials
in $14$ variables vanishing identically on the Hurwitz curve of
genus $14$. This leads to a fourteen-dimensional representation on
${\Bbb P}^{13}$. Now, we construct such a representation which is
induced from our six-dimensional representation.

  We have
$$\aligned
  &-13 \sqrt{13} ST^{\nu}(z_1) \cdot ST^{\nu}(z_2) \cdot ST^{\nu}(z_3)\\
 =&-\sqrt{\frac{-13-3 \sqrt{13}}{2}} (\zeta^{8 \nu} z_1^3+\zeta^{7 \nu} z_2^3+\zeta^{11 \nu} z_3^3)
   -\sqrt{\frac{-13+3 \sqrt{13}}{2}} (\zeta^{5 \nu} z_4^3+\zeta^{6 \nu} z_5^3+\zeta^{2 \nu} z_6^3)\\
  &-\sqrt{-13+2 \sqrt{13}} (\zeta^{12 \nu} z_1^2 z_2+\zeta^{4 \nu} z_2^2 z_3+\zeta^{10 \nu} z_3^2 z_1)\\
  &-\sqrt{-13-2 \sqrt{13}} (\zeta^{\nu} z_4^2 z_5+\zeta^{9 \nu} z_5^2 z_6+\zeta^{3 \nu} z_6^2 z_4)\\
  &+2 \sqrt{-13-2 \sqrt{13}} (\zeta^{3 \nu} z_1 z_2^2+\zeta^{\nu} z_2 z_3^2+\zeta^{9 \nu} z_3 z_1^2)\\
  &-2 \sqrt{-13+2 \sqrt{13}} (\zeta^{10 \nu} z_4 z_5^2+\zeta^{12 \nu} z_5 z_6^2+\zeta^{4 \nu} z_6 z_4^2)\\
  &+2 \sqrt{\frac{-13-3 \sqrt{13}}{2}} (\zeta^{7 \nu} z_1^2 z_4+\zeta^{11 \nu} z_2^2 z_5+\zeta^{8 \nu} z_3^2 z_6)+\\
  &-2 \sqrt{\frac{-13+3 \sqrt{13}}{2}} (\zeta^{6 \nu} z_1 z_4^2+\zeta^{2 \nu} z_2 z_5^2+\zeta^{5 \nu} z_3 z_6^2)+\\
  &+\sqrt{-13-2 \sqrt{13}} (\zeta^{3 \nu} z_1^2 z_5+\zeta^{\nu} z_2^2 z_6+\zeta^{9 \nu} z_3^2 z_4)+\\
  &+\sqrt{-13+2 \sqrt{13}} (\zeta^{10 \nu} z_2 z_4^2+\zeta^{12 \nu} z_3 z_5^2+\zeta^{4 \nu} z_1 z_6^2)+\\
  &+\sqrt{\frac{-13+3 \sqrt{13}}{2}} (\zeta^{6 \nu} z_1^2 z_6+\zeta^{2 \nu} z_2^2 z_4+\zeta^{5 \nu} z_3^2 z_5)+\\
  &+\sqrt{\frac{-13-3 \sqrt{13}}{2}} (\zeta^{7 \nu} z_3 z_4^2+\zeta^{11 \nu} z_1 z_5^2+\zeta^{8 \nu} z_2 z_6^2)+\\
  &+[2(\theta_1-\theta_3)-3(\theta_2-\theta_4)] z_1 z_2 z_3+[2(\theta_4-\theta_2)-3(\theta_1-\theta_3)] z_4 z_5 z_6+
\endaligned$$
$$\aligned
  &-\sqrt{\frac{-13-3 \sqrt{13}}{2}} (\zeta^{11 \nu} z_1 z_2 z_4+\zeta^{8 \nu} z_2 z_3 z_5+\zeta^{7 \nu} z_1 z_3 z_6)+\\
  &+\sqrt{\frac{-13+3 \sqrt{13}}{2}} (\zeta^{2 \nu} z_1 z_4 z_5+\zeta^{5 \nu} z_2 z_5 z_6+\zeta^{6 \nu} z_3 z_4 z_6)+\\
  &-3 \sqrt{\frac{-13-3 \sqrt{13}}{2}} (\zeta^{7 \nu} z_1 z_2 z_5+\zeta^{11 \nu} z_2 z_3 z_6+\zeta^{8 \nu} z_1 z_3 z_4)+\\
  &+3 \sqrt{\frac{-13+3 \sqrt{13}}{2}} (\zeta^{6 \nu} z_2 z_4 z_5+\zeta^{2 \nu} z_3 z_5 z_6+\zeta^{5 \nu} z_1 z_4 z_6)+\\
  &-\sqrt{-13+2 \sqrt{13}} (\zeta^{10 \nu} z_1 z_2 z_6+\zeta^{4 \nu} z_1 z_3 z_5+\zeta^{12 \nu} z_2 z_3 z_4)+\\
  &+\sqrt{-13-2 \sqrt{13}} (\zeta^{3 \nu} z_3 z_4 z_5+\zeta^{9 \nu} z_2 z_4 z_6+\zeta^{\nu} z_1 z_5 z_6).
\endaligned$$

  This leads us to define the following senary cubic forms (cubic forms in six variables):
$$\left\{\aligned
  {\Bbb D}_0 &=z_1 z_2 z_3,\\
  {\Bbb D}_1 &=2 z_2 z_3^2+z_2^2 z_6-z_4^2 z_5+z_1 z_5 z_6,\\
  {\Bbb D}_2 &=-z_6^3+z_2^2 z_4-2 z_2 z_5^2+z_1 z_4 z_5+3 z_3 z_5 z_6,\\
  {\Bbb D}_3 &=2 z_1 z_2^2+z_1^2 z_5-z_4 z_6^2+z_3 z_4 z_5,\\
  {\Bbb D}_4 &=-z_2^2 z_3+z_1 z_6^2-2 z_4^2 z_6-z_1 z_3 z_5,\\
  {\Bbb D}_5 &=-z_4^3+z_3^2 z_5-2 z_3 z_6^2+z_2 z_5 z_6+3 z_1 z_4 z_6,\\
  {\Bbb D}_6 &=-z_5^3+z_1^2 z_6-2 z_1 z_4^2+z_3 z_4 z_6+3 z_2 z_4 z_5,\\
  {\Bbb D}_7 &=-z_2^3+z_3 z_4^2-z_1 z_3 z_6-3 z_1 z_2 z_5+2 z_1^2 z_4,\\
  {\Bbb D}_8 &=-z_1^3+z_2 z_6^2-z_2 z_3 z_5-3 z_1 z_3 z_4+2 z_3^2 z_6,\\
  {\Bbb D}_9 &=2 z_1^2 z_3+z_3^2 z_4-z_5^2 z_6+z_2 z_4 z_6,\\
  {\Bbb D}_{10} &=-z_1 z_3^2+z_2 z_4^2-2 z_4 z_5^2-z_1 z_2 z_6,\\
  {\Bbb D}_{11} &=-z_3^3+z_1 z_5^2-z_1 z_2 z_4-3 z_2 z_3 z_6+2 z_2^2 z_5,\\
  {\Bbb D}_{12} &=-z_1^2 z_2+z_3 z_5^2-2 z_5 z_6^2-z_2 z_3 z_4,\\
  {\Bbb D}_{\infty}&=z_4 z_5 z_6.
\endaligned\right.\eqno{(6.1)}$$

  Let
$$r_0=2(\theta_1-\theta_3)-3(\theta_2-\theta_4), \quad
  r_{\infty}=2(\theta_4-\theta_2)-3(\theta_1-\theta_3),$$
and
$$r_1=\sqrt{-13-2 \sqrt{13}},
  r_2=\sqrt{\frac{-13+3 \sqrt{13}}{2}},
  r_3=\sqrt{-13+2 \sqrt{13}},
  r_4=\sqrt{\frac{-13-3 \sqrt{13}}{2}}.$$
Then
$$\aligned
  &-13 \sqrt{13} ST^{\nu}({\Bbb D}_0)\\
 =&r_0 {\Bbb D}_0+r_1 \zeta^{\nu} {\Bbb D}_1+r_2 \zeta^{2 \nu} {\Bbb D}_2+
   r_1 \zeta^{3 \nu} {\Bbb D}_3+r_3 \zeta^{4 \nu} {\Bbb D}_4+
   r_2 \zeta^{5 \nu} {\Bbb D}_5+r_2 \zeta^{6 \nu} {\Bbb D}_6+\\
  &+r_4 \zeta^{7 \nu} {\Bbb D}_7+r_4 \zeta^{8 \nu} {\Bbb D}_8
   +r_1 \zeta^{9 \nu} {\Bbb D}_9+r_3 \zeta^{10 \nu} {\Bbb D}_{10}+r_4 \zeta^{11 \nu} {\Bbb D}_{11}+
   r_3 \zeta^{12 \nu} {\Bbb D}_{12}+r_{\infty} {\Bbb D}_{\infty}.
\endaligned$$
$$\aligned
  -13 \sqrt{13} S({\Bbb D}_{\infty})
 =&r_{\infty} {\Bbb D}_0-r_3 {\Bbb D}_1-r_4 {\Bbb D}_2-r_3 {\Bbb D}_3+r_1 {\Bbb D}_4
  -r_4 {\Bbb D}_5-r_4 {\Bbb D}_6+\\
  &+r_2 {\Bbb D}_7+r_2 {\Bbb D}_8-r_3 {\Bbb D}_9+r_1 {\Bbb D}_{10}+r_2 {\Bbb D}_{11}+r_1 {\Bbb D}_{12}
   -r_0 {\Bbb D}_{\infty}.
\endaligned$$
Moreover, we have
$$\aligned
  -13 \sqrt{13} S({\Bbb D}_1)
 =&13 r_1 {\Bbb D}_0+q_1 {\Bbb D}_1+q_2 {\Bbb D}_2+q_3 {\Bbb D}_3+
   q_4 {\Bbb D}_4+q_5 {\Bbb D}_5+q_6 {\Bbb D}_6+\\
  &+q_7 {\Bbb D}_7+q_8 {\Bbb D}_8+q_9 {\Bbb D}_9+q_{10} {\Bbb D}_{10}
   +q_{11} {\Bbb D}_{11}+q_{12} {\Bbb D}_{12}-13 r_3 {\Bbb
   D}_{\infty},
\endaligned$$
where
$$\left\{\aligned
  q_1 &=-2 (\zeta-\zeta^{12})-2 (\zeta^5-\zeta^8)+6 (\zeta^3-\zeta^{10})
        -(\zeta^2-\zeta^{11})+4 (\zeta^9-\zeta^4)+2 (\zeta^6-\zeta^7),\\
  q_2 &=-4 (\zeta-\zeta^{12})+3 (\zeta^5-\zeta^8)+3 (\zeta^3-\zeta^{10})
        -(\zeta^2-\zeta^{11})-2 (\zeta^9-\zeta^4),\\
  q_3 &=6 (\zeta-\zeta^{12})-(\zeta^5-\zeta^8)+4 (\zeta^3-\zeta^{10})+
        2 (\zeta^2-\zeta^{11})-2 (\zeta^9-\zeta^4)-2 (\zeta^6-\zeta^7),\\
  q_4 &=-2 (\zeta-\zeta^{12})+4 (\zeta^5-\zeta^8)+2 (\zeta^3-\zeta^{10})
        -2 (\zeta^2-\zeta^{11})+(\zeta^9-\zeta^4)+6 (\zeta^6-\zeta^7),\\
  q_5 &=-2 (\zeta-\zeta^{12})-4 (\zeta^3-\zeta^{10})+3 (\zeta^2-\zeta^{11})
        +3 (\zeta^9-\zeta^4)-(\zeta^6-\zeta^7),\\
  q_6 &=3 (\zeta-\zeta^{12})-(\zeta^5-\zeta^8)-2 (\zeta^3-\zeta^{10})
        -4 (\zeta^9-\zeta^4)+3 (\zeta^6-\zeta^7),\\
  q_7 &=(\zeta-\zeta^{12})+3 (\zeta^5-\zeta^8)-2 (\zeta^2-\zeta^{11})
        -3 (\zeta^9-\zeta^4)-4 (\zeta^6-\zeta^7),\\
  q_8 &=-2 (\zeta^5-\zeta^8)-3 (\zeta^3-\zeta^{10})-4 (\zeta^2-\zeta^{11})
        +(\zeta^9-\zeta^4)+3 (\zeta^6-\zeta^7),\\
  q_9 &=4 (\zeta-\zeta^{12})+2 (\zeta^5-\zeta^8)-2 (\zeta^3-\zeta^{10})
        -2 (\zeta^2-\zeta^{11})+6 (\zeta^9-\zeta^4)-(\zeta^6-\zeta^7),\\
  q_{10} &=(\zeta-\zeta^{12})+6 (\zeta^5-\zeta^8)-2 (\zeta^3-\zeta^{10})
          +4 (\zeta^2-\zeta^{11})+2 (\zeta^9-\zeta^4)-2 (\zeta^6-\zeta^7),\\
  q_{11} &=-3 (\zeta-\zeta^{12})-4 (\zeta^5-\zeta^8)+(\zeta^3-\zeta^{10})
           +3 (\zeta^2-\zeta^{11})-2 (\zeta^6-\zeta^7),\\
  q_{12} &=2 (\zeta-\zeta^{12})-2 (\zeta^5-\zeta^8)+(\zeta^3-\zeta^{10})
           +6 (\zeta^2-\zeta^{11})-2 (\zeta^9-\zeta^4)+4 (\zeta^6-\zeta^7).
\endaligned\right.\eqno{(6.2)}$$
Similarly, we obtain
$$\aligned
   -13 \sqrt{13} S({\Bbb D}_2)
 =&26 r_2 {\Bbb D}_0+2 q_2 {\Bbb D}_1-q_4 {\Bbb D}_2+2 q_6 {\Bbb D}_3+
   2 q_8 {\Bbb D}_4-q_{10} {\Bbb D}_5-q_{12} {\Bbb D}_6+\\
  &+q_1 {\Bbb D}_7+q_3 {\Bbb D}_8+2 q_5 {\Bbb D}_9+2 q_7 {\Bbb D}_{10}+
   q_9 {\Bbb D}_{11}+2 q_{11} {\Bbb D}_{12}-26 r_4 {\Bbb D}_{\infty}.
\endaligned$$
$$\aligned
  -13 \sqrt{13} S({\Bbb D}_3)
 =&13 r_1 {\Bbb D}_0+q_3 {\Bbb D}_1+q_6 {\Bbb D}_2+q_9 {\Bbb D}_3+
   q_{12} {\Bbb D}_4+q_2 {\Bbb D}_5+q_5 {\Bbb D}_6+\\
  &+q_8 {\Bbb D}_7+q_{11} {\Bbb D}_8+q_1 {\Bbb D}_9+q_4 {\Bbb D}_{10}
   +q_7 {\Bbb D}_{11}+q_{10} {\Bbb D}_{12}-13 r_3 {\Bbb D}_{\infty},
\endaligned$$
$$\aligned
  -13 \sqrt{13} S({\Bbb D}_4)
 =&13 r_3 {\Bbb D}_0+q_4 {\Bbb D}_1+q_8 {\Bbb D}_2+q_{12} {\Bbb D}_3-
   q_3 {\Bbb D}_4+q_7 {\Bbb D}_5+q_{11} {\Bbb D}_6+\\
  &-q_2 {\Bbb D}_7-q_6 {\Bbb D}_8+q_{10} {\Bbb D}_9-q_1 {\Bbb D}_{10}
   -q_5 {\Bbb D}_{11}-q_9 {\Bbb D}_{12}+13 r_1 {\Bbb D}_{\infty},
\endaligned$$
$$\aligned
   -13 \sqrt{13} S({\Bbb D}_5)
 =&26 r_2 {\Bbb D}_0+2 q_5 {\Bbb D}_1-q_{10} {\Bbb D}_2+2 q_2 {\Bbb D}_3+
   2 q_7 {\Bbb D}_4-q_{12} {\Bbb D}_5-q_4 {\Bbb D}_6+\\
  &+q_9 {\Bbb D}_7+q_1 {\Bbb D}_8+2 q_6 {\Bbb D}_9+2 q_{11} {\Bbb D}_{10}+
   q_3 {\Bbb D}_{11}+2 q_8 {\Bbb D}_{12}-26 r_4 {\Bbb D}_{\infty}.
\endaligned$$
$$\aligned
   -13 \sqrt{13} S({\Bbb D}_6)
 =&26 r_2 {\Bbb D}_0+2 q_6 {\Bbb D}_1-q_{12} {\Bbb D}_2+2 q_5 {\Bbb D}_3+
   2 q_{11} {\Bbb D}_4-q_4 {\Bbb D}_5-q_{10} {\Bbb D}_6+\\
  &+q_3 {\Bbb D}_7+q_9 {\Bbb D}_8+2 q_2 {\Bbb D}_9+2 q_8 {\Bbb D}_{10}+
   q_1 {\Bbb D}_{11}+2 q_7 {\Bbb D}_{12}-26 r_4 {\Bbb D}_{\infty}.
\endaligned$$
$$\aligned
   -13 \sqrt{13} S({\Bbb D}_7)
 =&26 r_4 {\Bbb D}_0+2 q_7 {\Bbb D}_1+q_1 {\Bbb D}_2+2 q_8 {\Bbb D}_3-
   2 q_2 {\Bbb D}_4+q_9 {\Bbb D}_5+q_3 {\Bbb D}_6+\\
  &+q_{10} {\Bbb D}_7+q_4 {\Bbb D}_8+2 q_{11} {\Bbb D}_9-2 q_5 {\Bbb D}_{10}+
   q_{12} {\Bbb D}_{11}-2 q_6 {\Bbb D}_{12}+26 r_2 {\Bbb D}_{\infty}.
\endaligned$$
$$\aligned
   -13 \sqrt{13} S({\Bbb D}_8)
 =&26 r_4 {\Bbb D}_0+2 q_8 {\Bbb D}_1+q_3 {\Bbb D}_2+2 q_{11} {\Bbb D}_3-
   2 q_6 {\Bbb D}_4+q_1 {\Bbb D}_5+q_9 {\Bbb D}_6+\\
  &+q_4 {\Bbb D}_7+q_{12} {\Bbb D}_8+2 q_7 {\Bbb D}_9-2 q_2 {\Bbb D}_{10}+
   q_{10} {\Bbb D}_{11}-2 q_5 {\Bbb D}_{12}+26 r_2 {\Bbb D}_{\infty}.
\endaligned$$
$$\aligned
  -13 \sqrt{13} S({\Bbb D}_9)
 =&13 r_1 {\Bbb D}_0+q_9 {\Bbb D}_1+q_5 {\Bbb D}_2+q_1 {\Bbb D}_3+
   q_{10} {\Bbb D}_4+q_6 {\Bbb D}_5+q_2 {\Bbb D}_6+\\
  &+q_1 {\Bbb D}_7+q_7 {\Bbb D}_8+q_8 {\Bbb D}_9+q_{12} {\Bbb D}_{10}
   +q_8 {\Bbb D}_{11}+q_4 {\Bbb D}_{12}-13 r_3 {\Bbb D}_{\infty},
\endaligned$$
$$\aligned
  -13 \sqrt{13} S({\Bbb D}_{10})
 =&13 r_3 {\Bbb D}_0+q_{10} {\Bbb D}_1+q_7 {\Bbb D}_2+q_4 {\Bbb D}_3-
   q_1 {\Bbb D}_4+q_{11} {\Bbb D}_5+q_8 {\Bbb D}_6+\\
  &-q_5 {\Bbb D}_7-q_2 {\Bbb D}_8+q_{12} {\Bbb D}_9-q_9 {\Bbb D}_{10}
   -q_6 {\Bbb D}_{11}-q_3 {\Bbb D}_{12}+13 r_1 {\Bbb D}_{\infty},
\endaligned$$
$$\aligned
   -13 \sqrt{13} S({\Bbb D}_{11})
 =&26 r_4 {\Bbb D}_0+2 q_{11} {\Bbb D}_1+q_9 {\Bbb D}_2+2 q_7 {\Bbb D}_3-
   2 q_5 {\Bbb D}_4+q_3 {\Bbb D}_5+q_1 {\Bbb D}_6+\\
  &+q_{12} {\Bbb D}_7+q_{10} {\Bbb D}_8+2 q_8 {\Bbb D}_9-2 q_6 {\Bbb D}_{10}+
   q_4 {\Bbb D}_{11}-2 q_2 {\Bbb D}_{12}+26 r_2 {\Bbb D}_{\infty}.
\endaligned$$
$$\aligned
  -13 \sqrt{13} S({\Bbb D}_{12})
 =&13 r_3 {\Bbb D}_0+q_{12} {\Bbb D}_1+q_{11} {\Bbb D}_2+q_{10} {\Bbb D}_3-
   q_9 {\Bbb D}_4+q_8 {\Bbb D}_5+q_7 {\Bbb D}_6+\\
  &-q_6 {\Bbb D}_7-q_5 {\Bbb D}_8+q_4 {\Bbb D}_9-q_3 {\Bbb D}_{10}
   -q_2 {\Bbb D}_{11}-q_1 {\Bbb D}_{12}+13 r_1 {\Bbb D}_{\infty}.
\endaligned$$
Hence, we get the element $\widehat{S}$ which is induced from the
action of $S$ on the basis $({\Bbb D}_0, \cdots, {\Bbb
D}_{\infty})$:
$$\widehat{S}=-\frac{1}{13 \sqrt{13}} \left(\matrix S_1 & S_2\\
               S_3 & S_4 \endmatrix\right),\eqno{(6.3)}$$
with
$$S_1=\left(\matrix
  r_0 & r_1 & r_2 & r_1 & r_3 & r_2 & r_2\\
  13 r_1 & q_1 & q_2 & q_3 & q_4 & q_5 & q_6\\
  26 r_2 & 2 q_2 & -q_4 & 2 q_6 & 2 q_8 & -q_{10} & -q_{12}\\
  13 r_1 & q_3 & q_6 & q_9 & q_{12} & q_2 & q_5\\
  13 r_3 & q_4 & q_8 & q_{12} & -q_3 & q_7 & q_{11}\\
  26 r_2 & 2 q_5 & -q_{10} & 2 q_2 & 2 q_7 & -q_{12} & -q_4\\
  26 r_2 & 2 q_6 & -q_{12} & 2 q_5 & 2 q_{11} & -q_4 & -q_{10}
\endmatrix\right),\eqno{(6.4)}$$
$$S_2=\left(\matrix
  r_4 & r_4 & r_1 & r_3 & r_4 & r_3 & r_{\infty}\\
  q_7 & q_8 & q_9 & q_{10} & q_{11} & q_{12} & -13 r_3\\
  q_1 & q_3 & 2 q_5 & 2 q_7 & q_9 & 2 q_{11} & -26 r_4\\
  q_8 & q_{11} & q_1 & q_4 & q_7 & q_{10} & -13 r_3\\
  -q_2 & -q_6 & q_{10} & -q_1 & -q_5 & -q_9 & 13 r_1\\
  q_9 & q_1 & 2 q_6 & 2 q_{11} & q_3 & 2 q_8 & -26 r_4\\
  q_3 & q_9 & 2 q_2 & 2 q_8 & q_1 & 2 q_7 & -26 r_4
\endmatrix\right),\eqno{(6.5)}$$
$$S_3=\left(\matrix
  26 r_4 & 2 q_7 & q_1 & 2 q_8 & -2 q_2 & q_9 & q_3\\
  26 r_4 & 2 q_8 & q_3 & 2 q_{11} & -2 q_6 & q_1 & q_9\\
  13 r_1 & q_9 & q_5 & q_1 & q_{10} & q_6 & q_2\\
  13 r_3 & q_{10} & q_7 & q_4 & -q_1 & q_{11} & q_8\\
  26 r_4 & 2 q_{11} & q_9 & 2 q_7 & -2 q_5 & q_3 & q_1\\
  13 r_3 & q_{12} & q_{11} & q_{10} & -q_9 & q_8 & q_7\\
  r_{\infty} & -r_3 & -r_4 & -r_3 & r_1 & -r_4 & -r_4
\endmatrix\right),\eqno{(6.6)}$$
$$S_4=\left(\matrix
  q_{10} & q_4 & 2 q_{11} & -2 q_5 & q_{12} & -2 q_6 & 26 r_2\\
  q_4 & q_{12} & 2 q_7 & -2 q_2 & q_{10} & -2 q_5 & 26 r_2\\
  q_{11} & q_7 & q_3 & q_{12} & q_8 & q_4 & -13 r_3\\
  -q_5 & -q_2 & q_{12} & -q_9 & -q_6 & -q_3 & 13 r_1\\
  q_{12} & q_{10} & 2 q_8 & -2 q_6 & q_4 & -2 q_2 & 26 r_2\\
  -q_6 & -q_5 & q_4 & -q_3 & -q_2 & -q_1 & 13 r_1\\
  r_2 & r_2 & -r_3 & r_1 & r_2 & r_1 & -r_0
\endmatrix\right).\eqno{(6.7)}$$
Similarly, the element $\widehat{T}$ is induced from the action of
$T$ on the basis $({\Bbb D}_0, \cdots, {\Bbb D}_{\infty})$:
$$\widehat{T}=\text{diag}(1, \zeta, \zeta^2, \zeta^3, \zeta^4,
              \zeta^5, \zeta^6, \zeta^7, \zeta^8, \zeta^9, \zeta^{10},
              \zeta^{11}, \zeta^{12}, 1).\eqno{(6.8)}$$
We have
$$\text{Tr}(\widehat{S})=0, \quad \text{Tr}(\widehat{T})=1, \quad
  \text{Tr}(\widehat{S} \widehat{T})=-2.\eqno{(6.9)}$$
Hence, this fourteen-dimensional representation corresponds to the
character $\chi_{15}$ in Table $1$. It is defined over the
cyclotomic field ${\Bbb Q}(\zeta)$.

  Let
$$\delta_{\infty}(z_1, z_2, z_3, z_4, z_5, z_6)=13^2 (z_1^2 z_2^2 z_3^2+z_4^2 z_5^2 z_6^2)\eqno{(6.10)}$$
and
$$\delta_{\nu}(z_1, z_2, z_3, z_4, z_5, z_6)=\delta_{\infty}(ST^{\nu}(z_1, z_2, z_3, z_4, z_5, z_6))\eqno{(6.11)}$$
for $\nu=0, 1, \cdots, 12$. Then
$$\delta_{\nu}=13^2 ST^{\nu}({\Bbb G}_0)=-13 {\Bbb G}_0+\zeta^{\nu} {\Bbb G}_1
  +\zeta^{2 \nu} {\Bbb G}_2+\cdots+\zeta^{12 \nu} {\Bbb G}_{12},\eqno{(6.12)}$$
where the senary sextic forms (i.e., sextic forms in six
variables) are given as follows:
$$\left\{\aligned
  {\Bbb G}_0 &={\Bbb D}_0^2+{\Bbb D}_{\infty}^2,\\
  {\Bbb G}_1 &=-{\Bbb D}_7^2+2 {\Bbb D}_0 {\Bbb D}_1+10 {\Bbb D}_{\infty} {\Bbb D}_1
               +2 {\Bbb D}_2 {\Bbb D}_{12}-2 {\Bbb D}_3 {\Bbb D}_{11}-4 {\Bbb D}_4 {\Bbb D}_{10}
               -2 {\Bbb D}_9 {\Bbb D}_5,\\
  {\Bbb G}_2 &=-2 {\Bbb D}_1^2-4 {\Bbb D}_0 {\Bbb D}_2+6 {\Bbb D}_{\infty} {\Bbb D}_2
               -2 {\Bbb D}_4 {\Bbb D}_{11}+2 {\Bbb D}_5 {\Bbb D}_{10}-2 {\Bbb D}_6 {\Bbb D}_9
               -2 {\Bbb D}_7 {\Bbb D}_8,\\
  {\Bbb G}_3 &=-{\Bbb D}_8^2+2 {\Bbb D}_0 {\Bbb D}_3+10 {\Bbb D}_{\infty} {\Bbb D}_3
               +2 {\Bbb D}_6 {\Bbb D}_{10}-2 {\Bbb D}_9 {\Bbb D}_7-4 {\Bbb D}_{12} {\Bbb D}_4
               -2 {\Bbb D}_1 {\Bbb D}_2,\\
  {\Bbb G}_4 &=-{\Bbb D}_2^2+10 {\Bbb D}_0 {\Bbb D}_4-2 {\Bbb D}_{\infty} {\Bbb D}_4
               +2 {\Bbb D}_5 {\Bbb D}_{12}-2 {\Bbb D}_9 {\Bbb D}_8-4 {\Bbb D}_1 {\Bbb D}_3
               -2 {\Bbb D}_{10} {\Bbb D}_7,\\
  {\Bbb G}_5 &=-2 {\Bbb D}_9^2-4 {\Bbb D}_0 {\Bbb D}_5+6 {\Bbb D}_{\infty} {\Bbb D}_5
               -2 {\Bbb D}_{10} {\Bbb D}_8+2 {\Bbb D}_6 {\Bbb D}_{12}-2 {\Bbb D}_2 {\Bbb D}_3
               -2 {\Bbb D}_{11} {\Bbb D}_7,\\
  {\Bbb G}_6 &=-2 {\Bbb D}_3^2-4 {\Bbb D}_0 {\Bbb D}_6+6 {\Bbb D}_{\infty} {\Bbb D}_6
               -2 {\Bbb D}_{12} {\Bbb D}_7+2 {\Bbb D}_2 {\Bbb D}_4-2 {\Bbb D}_5 {\Bbb D}_1
               -2 {\Bbb D}_8 {\Bbb D}_{11},\\
  {\Bbb G}_7 &=-2 {\Bbb D}_{10}^2+6 {\Bbb D}_0 {\Bbb D}_7+4 {\Bbb D}_{\infty} {\Bbb D}_7
               -2 {\Bbb D}_1 {\Bbb D}_6-2 {\Bbb D}_2 {\Bbb D}_5-2 {\Bbb D}_8 {\Bbb D}_{12}
               -2 {\Bbb D}_9 {\Bbb D}_{11},\\
  {\Bbb G}_8 &=-2 {\Bbb D}_4^2+6 {\Bbb D}_0 {\Bbb D}_8+4 {\Bbb D}_{\infty} {\Bbb D}_8
               -2 {\Bbb D}_3 {\Bbb D}_5-2 {\Bbb D}_6 {\Bbb D}_2-2 {\Bbb D}_{11} {\Bbb D}_{10}
               -2 {\Bbb D}_1 {\Bbb D}_7,\\
  {\Bbb G}_9 &=-{\Bbb D}_{11}^2+2 {\Bbb D}_0 {\Bbb D}_9+10 {\Bbb D}_{\infty} {\Bbb D}_9
               +2 {\Bbb D}_5 {\Bbb D}_4-2 {\Bbb D}_1 {\Bbb D}_8-4 {\Bbb D}_{10} {\Bbb D}_{12}
               -2 {\Bbb D}_3 {\Bbb D}_6,\\
  {\Bbb G}_{10} &=-{\Bbb D}_5^2+10 {\Bbb D}_0 {\Bbb D}_{10}-2 {\Bbb D}_{\infty} {\Bbb D}_{10}
               +2 {\Bbb D}_6 {\Bbb D}_4-2 {\Bbb D}_3 {\Bbb D}_7-4 {\Bbb D}_9 {\Bbb D}_1
               -2 {\Bbb D}_{12} {\Bbb D}_{11},\\
  {\Bbb G}_{11} &=-2 {\Bbb D}_{12}^2+6 {\Bbb D}_0 {\Bbb D}_{11}+4 {\Bbb D}_{\infty} {\Bbb D}_{11}
               -2 {\Bbb D}_9 {\Bbb D}_2-2 {\Bbb D}_5 {\Bbb D}_6-2 {\Bbb D}_7 {\Bbb D}_4
               -2 {\Bbb D}_3 {\Bbb D}_8,\\
  {\Bbb G}_{12} &=-{\Bbb D}_6^2+10 {\Bbb D}_0 {\Bbb D}_{12}-2 {\Bbb D}_{\infty} {\Bbb D}_{12}
               +2 {\Bbb D}_2 {\Bbb D}_{10}-2 {\Bbb D}_1 {\Bbb D}_{11}-4 {\Bbb D}_3 {\Bbb D}_9
               -2 {\Bbb D}_4 {\Bbb D}_8.
\endaligned\right.\eqno{(6.13)}$$
We have that ${\Bbb G}_0$ is invariant under the action of
$\langle H, T \rangle$, the maximal subgroup of order $78$ of $G$
with index $14$. Note that $\delta_{\infty}$, $\delta_{\nu}$ for
$\nu=0, \cdots, 12$ form an algebraic equation of degree fourteen.
However, we have $\delta_{\infty}+\sum_{\nu=0}^{12}
\delta_{\nu}=0$. Hence, it is not the Jacobian equation of degree
fourteen! We call it exotic modular equation of degree fourteen.
We have
$$\delta_{\infty}^2+\sum_{\nu=0}^{12} \delta_{\nu}^2
 =26 (7 \cdot 13^2 {\Bbb G}_0^2+{\Bbb G}_1 {\Bbb G}_{12}+{\Bbb G}_ 2
  {\Bbb G}_{11}+\cdots+{\Bbb G}_6 {\Bbb G}_7).\eqno{(6.14)}$$

  Let us recall some basic fact about the Monster finite simple
group ${\Bbb M}$ (see \cite{CN}). For each prime $p$ with
$(p-1)|24$ there is a conjugacy class of elements of ${\Bbb M}$,
with centraliser of form $p^{1+2d}.G_p$, where $p.G_p$ is the
centraliser of a corresponding automorphism of the Leech Lattice.
The symbol $p^{1+2d}$ denotes an extraspecial $p$-group, and $2d=
24/(p-1)$. Some maximal $p$-local subgroups of ${\Bbb M}$ are
given as follows (see \cite{CC}):
$$\matrix
  p & \text{structure} & \text{specification}\\
  2 & 2_{+}^{1+24} \cdot Co_1 & N(2B)\\
  3 & 3_{+}^{1+12} \cdot 2 Suz:2 & N(3B)\\
  5 & 5_{+}^{1+6}: 4J_2 \cdot 2 & N(5B)\\
  7 & 7_{+}^{1+4}:(3 \times 2S_7) & N(7B)\\
 13 & 13_{+}^{1+2}:(3 \times 4 S_4) & N(13B)
\endmatrix$$
and
$$\matrix
  p & \text{structure} & \text{specification}\\
  2 & 2^2.2^{11}.2^{22}.(S_3 \times M_{24}) & N(2B^2)\\
  3 & 3^2.3^5.3^{10}.(M_{11} \times 2S_4) & N(3B^2)\\
  5 & 5^2.5^2.5^4.(S_3 \times GL_2(5)) & N(5B^2)\\
  7 & 7^2.7.7^2:GL_2(7) & N(7B^2)\\
 13 & 13^2:4 L_2(13) \cdot 2 & N(13B^2)
\endmatrix$$
We will study the relation between genus-zero subgroups (both
congruence subgroups and non-congruence subgroups) of the modular
group and subgroups of the Monster simple group. In particular, we
are interested in $\Gamma_0(13)$ and $\Gamma_{7, 7}$, the quotient
groups $\Gamma/\Gamma(13) \cong \Gamma/G_i \cong PSL(2, 13)$ for
$i=1, 2, 3$. Note that
$$M:=7 \cdot 13^2 {\Bbb G}_0^2+{\Bbb G}_1 {\Bbb G}_{12}+{\Bbb G}_ 2
     {\Bbb G}_{11}+\cdots+{\Bbb G}_6 {\Bbb G}_7=0\eqno{(6.15)}$$
is a $12$-dimensional $PSL(2, 13)$-invariant complex algebraic
variety (i.e., $24$-dimensional manifold) of degree $4$ in the
projective space ${\Bbb P}^{13}=\{({\Bbb D}_0, {\Bbb D}_1, \cdots,
{\Bbb D}_{12}, {\Bbb D}_{\infty})\}$, which is related to
$N(13B^2)$, a maximal $p$-local subgroup of ${\Bbb M}$. In fact,
both $PSL(2, 13)$ and ${\Bbb M}$ are Hurwitz groups.

\vskip 2.0 cm

{\smc Department of Mathematics, Peking University}

{\smc Beijing 100871, P. R. China}

{\it E-mail address}: yanglei\@math.pku.edu.cn
\vskip 1.5 cm
\Refs

\item{[AH]} {\smc M. Asaeda and U. Haagerup}, Exotic subfactors of
            finite depth with Jones indices $(5+\sqrt{13})/2$ and
            $(5+\sqrt{17})/2$, Comm. Math. Phys. {\bf 202} (1999), 1-63.

\item{[ASD]} {\smc A. O. L. Atkin and H. P. F. Swinnerton-Dyer},
             Modular forms on noncongruence subgroups, {\it
             Combinatorics}, 1-25, Proc. Sympos. Pure Math. {\bf
             19}, Amer. Math. Soc., Providence, R.I., 1971.

\item{[B1]} {\smc R. L. Bryant}, \'{E}lie Cartan and geometric
             duality, in: {\it Journ\'{e}es \'{E}lie Cartan 1998
             et 1999}, {\bf 16} (2000), 5-20, Institut \'{E}lie
             Cartan (Nancy, France).

\item{[B2]} {\smc R. L. Bryant}, Geometry of manifolds with
            special holonomy: ``100 years of holonomy'', in: {\it 150 years of
            mathematics at Washington University in St. Louis}, 29-38,
            Contemp. Math., {\bf 395}, Amer. Math. Soc., Providence, RI, 2006.

\item{[CW]} {\smc A. M. Cohen and D. B. Wales}, Finite subgroups
            of $G_2({\Bbb C})$, Comm. Algebra {\bf 11} (1983), 441-459.

\item{[CN]} {\smc J. Conway and S. Norton}, Monstrous moonshine,
            Bull. London Math. Soc. {\bf 11} (1979), 308-339.

\item{[CC]} {\smc J. H. Conway, R. T. Curtis, S. P. Norton, R.
             A. Parker and R. A. Wilson}, {\it Atlas of Finite Groups,
             Maximal Subgroups and Ordinary Characters for Simple Groups},
             Clarendon Press, Oxford, 1985.

\item{[E]} {\smc W. Ebeling}, {\it Lattices and Codes, A Course
            Partially Based on Lectures by F. Hirzebruch}, Second revised edition,
            Advanced Lectures in Mathematics, Friedr. Vieweg \& Sohn, Braunschweig, 2002.

\item{[El1]} {\smc N. D. Elkies}, Shimura curve computations,
            {\it Algorithmic number theory (Portland, OR, 1998)}, 1-47,
            Lecture Notes in Comput. Sci., {\bf 1423}, Springer, Berlin, 1998.

\item{[El2]} {\smc N. D. Elkies}, The Klein quartic in number
             theory, {\it The eightfold way}, 51-101, Math. Sci. Res. Inst.
             Publ., {\bf 35}, Cambridge Univ. Press, Cambridge, 1999.

\item{[EG1]} {\smc D. E. Evans and T. Gannon}, The exoticness and
             realisability of twisted Haagerup-Izumi modular data, Comm. Math.
             Phys. {\bf 307} (2011), 463-512.

\item{[EG2]} {\smc D. E. Evans and T. Gannon}, The search for the
            exotic-subfactors and conformal field theories,
            preprint, 2012.

\item{[FK]} {\smc H. M. Farkas and I. Kra}, {\it Theta Constants,
            Riemann Surfaces and the Modular Group, An Introduction
            with Applications to Uniformization Theorems, Partition
            Identities and Combinatorial Number Theory}, Graduate
            Studies in Mathematics, {\bf 37}, American Mathematical
            Society, Providence, RI, 2001.

\item{[F]} {\smc R. Fricke}, Ueber eine einfache Gruppe von $504$
           Operationen, Math. Ann. {\bf 52} (1899), 321-339.

\item{[He]} {\smc E. Hecke}, Grundlagen einer Theorie der
            Integralgruppen und der Integralperioden bei den Normalteilern der
            Modulgruppe, Math. Ann. {\bf 116} (1939), 469-510, in:
            {\it Mathematische Werke}, 731-772, G\"{o}ttingen, 1959.

\item{[Hi]} {\smc F. Hirzebruch}, {\it Gesammelte Abhandlungen},
            Bd. II, 1963-1987, Springer-Verlag, 1987.

\item{[Hi1]} {\smc F. Hirzebruch}, The ring of Hilbert modular
            forms for real quadratic fields of small discriminant,
            {\it Modular Functions of One Variable VI, Proceedings,
            Bonn 1976}, Edited by J. P. Serre and D. B. Zagier,
            287-323, Lecture Notes in Math. {\bf 627}, Springer-Verlag,
            1977, in: {\it Gesammelte Abhandlungen}, Bd. II, 501-536,
            Springer-Verlag, 1987.

\item{[HRW]} {\smc S.-M. Hong, E. Rowell and Z. Wang}, On exotic
             modular tensor categories, Commun. Contemp. Math. {\bf 10} (2008),
             suppl. 1, 1049-1074.

\item{[Hu]} {\smc A. Hurwitz}, Grundlagen einer independenten
            Theorie der elliptischen Modulfunktionen und Theorie
            der Multiplikator-Gleichungen erster Stufe, Math. Ann.
            {\bf 18} (1881), 528-592, in: {\it Mathematische
            Werke}, Bd. 1, 1-66, E. Birkh\"{a}user, 1932-33.

\item{[Iz]} {\smc M. Izumi}, The structure of sectors associated
            with Longo-Rehren inclusions, II. Examples, Rev. Math. Phys.
            {\bf 13} (2001), 603-674.

\item{[J]} {\smc G. A. Jones}, Congruence and noncongruence
           subgroups of the modular group: a survey, {\it
           Proceedings of groups---St. Andrews 1985}, 223-234,
           London Math. Soc. Lecture Note Series {\bf 121}, Cambridge
           Univ. Press, Cambridge, 1986.

\item{[J1]} {\smc V. F. R. Jones}, In and around the origin of
            quantum groups, in: {\it Prospects in mathematical physics},
            101-126, Contemp. Math., {\bf 437}, Amer. Math. Soc., Providence,
            RI, 2007.

\item{[J2]} {\smc V. F. R. Jones}, On the origin and development
             of subfactors and quantum topology. Bull. Amer. Math. Soc.
             {\bf 46} (2009), 309-326.

\item{[Ka]} {\smc S. Katok}, {\it Fuchsian groups}, Chicago
            Lectures in Mathematics Series, University of Chicago
            Press, Chicago, IL, 1992.

\item{[KSV]} {\smc M. G. Katz, M. Schaps and U. Vishne}, Hurwitz
             quaternion order and arithmetic Riemann surfaces,
             Geom. Dedicata {\bf 155} (2011), 151-161.

\item{[K]} {\smc F. Klein}, {\it Lectures on the Icosahedron and
             the Solution of Equations of the Fifth Degree},
             Translated by G. G. Morrice, second and revised edition,
             Dover Publications, Inc., 1956.

\item{[K1]} {\smc F. Klein}, Ueber die Transformation der
            elliptischen Functionen und die Aufl\"{o}sung der
            Gleichungen f\"{u}nften Grades, Math. Ann. {\bf 14}
            (1879), 111-172, in: {\it Gesammelte Mathematische
            Abhandlungen}, Bd. III, 13-75, Springer-Verlag, Berlin,
            1923.

\item{[K2]} {\smc F. Klein}, Ueber die Transformation siebenter
            Ordnung der elliptischen Functionen, Math. Ann. {\bf
            14} (1879), 428-471, in: {\it Gesammelte Mathematische
            Abhandlungen}, Bd. III, 90-136, Springer-Verlag, Berlin,
            1923.

\item{[K3]} {\smc F. Klein}, Ueber die Aufl\"{o}sung gewisser
            Gleichungen vom siebenten und achten Grade, Math. Ann.
            {\bf 15} (1879), 251-282, in: {\it Gesammelte Mathematische
            Abhandlungen}, Bd. II, 390-438, Springer-Verlag, Berlin,
            1922.

\item{[K4]} {\smc F. Klein}, Ueber die Transformation elfter
            Ordnung der elliptischen Functionen, Math. Ann. {\bf
            15} (1879), 533-555, in: {\it Gesammelte Mathematische
            Abhandlungen}, Bd. III, 140-168, Springer-Verlag, Berlin,
            1923.

\item{[L]} {\smc P. Lochak}, On arithmetic curves in the moduli
           spaces of curves, J. Inst. Math. Jussieu {\bf 4} (2005), 443-508.

\item{[M1]} {\smc A. M. Macbeath}, On a curve of genus $7$, Proc.
            London Math. Soc. {\bf 15} (1965), 527-542.

\item{[M2]} {\smc A. M. Macbeath}, Generators of the linear
            fractional groups, {\it Number Theory}, 14-32,
            Proc. Sympos. Pure Math., {\bf 12}, Amer. Math. Soc.,
            Providence, R.I., 1969.

\item{[Ma]} {\smc H. van Maldeghem}, {\it Generalized polygons},
            Monographs in Mathematics, {\bf 93}, Birkh\"{a}user Verlag,
            Basel, 1998.

\item{[MS]} {\smc C. L. Mallows and N. J. A. Sloane}, On the
            invariants of a linear group of order $336$, Proc.
            Cambridge Philos. Soc. {\bf 74} (1973), 435-440.

\item{[Mi]} {\smc G. A. Miller}, On the transitive substitution
            groups of degrees thirteen and fourteen, Quart. J.
            Pure Appl. Math. {\bf 29} (1898), 224-249.

\item{[MM]} {\smc I. Moreno-Mej\'{i}a}, The quadrics through the
            Hurwitz curves of genus $14$, J. London Math. Soc.
            {\bf 81} (2010), 374-388.

\item{[N]} {\smc M. Newman}, Maximal normal subgroups of the
            modular group, Proc. Amer. Math. Soc. {\bf 19} (1968),
            1138-1144.

\item{[SW]} {\smc N. Sato and M. Wakui}, Computations of
            Turaev-Viro-Ocneanu invariants of $3$-manifolds from subfactors,
            J. Knot Theory Ramifications {\bf 12} (2003), 543-574.

\item{[Sh]} {\smc G. Shimura}, Construction of class fields and
            zeta functions of algebraic curves, Ann. of Math. {\bf
            85} (1967), 58-159, in: {\it Collected Papers}, Vol.
            II, 30-131, Springer-Verlag, 2002.

\item{[S]} {\smc A. Sinkov}, Necessary and sufficient conditions
            for generating certain simple groups by two operators
            of periods two and three, Amer. J. Math. {\bf 59}
            (1937), 67-76.

\item{[St]} {\smc M. Streit}, Field of definition and Galois
            orbits for the Macbeath-Hurwitz curves, Arch. Math.
           {\bf 74} (2000), 342-349.

\item{[T]} {\smc J. Tits}, Sur la trialit\'{e} et certains groupes
           qui s'en d\'{e}duisent. Publ. Math. IHES {\bf 2} (1959), 13-60.

\item{[V]} {\smc R. Vogeler}, On the geometry of Hurwitz surfaces,
           Thesis, Florida State University (2003).

\item{[W]} {\smc A. Weil}, Sur certains groupes d'op\'{e}rateurs
           unitaires, Acta Math. {\bf 111} (1964), 143-211.

\item{[Wi]} {\smc E. Witten}, Quantum field theory and the Jones
            polynomial, Comm. Math. Phys. {\bf 121} (1989), 351-399.

\item{[W1]} {\smc K. Wohlfahrt}, An extension of F. Klein's level
           concept, Illinois J. Math. {\bf 8} (1964), 529-535.

\item{[W2]} {\smc K. Wohlfahrt}, Macbeath's curve and the modular
            group, Glasgow Math. J. {\bf 27} (1985), 239-247.

\endRefs
\end{document}